\documentclass[reqno,10pt]{amsart} \usepackage{amsmath,amsxtra,amssymb,latexsym, amscd,amsthm}
\usepackage[mathscr]{euscript}
\usepackage{mathrsfs}
\usepackage{color}
\usepackage[active]{srcltx}
\usepackage{mathtools,stmaryrd}
\usepackage{enumerate}
\usepackage{hyperref}
\usepackage{framed}
\usepackage[T1]{fontenc}
\usepackage{subfigure}
\usepackage{multirow}
\usepackage{booktabs}
\usepackage{tikz}
\usepackage{longtable}
\usepackage{tikz-qtree}
\usepackage[edges]{forest}
\usetikzlibrary{trees} 
\usepackage{stmaryrd}
\SetSymbolFont{stmry}{bold}{U}{stmry}{m}{n}
\SetSymbolFont{stmry}{bold}{U}{stmry}{b}{n}
\usepackage{graphicx}
\usepackage{subcaption}

\newcommand{\RR}{{\mathbb R}}

\newcommand{\bS}{{\mathbb S}}

\newcommand{\N}{\mathbb{N}}

\newcommand{\K}{\mathcal{K}}
\newcommand{\U}{\mathcal{U}}
\newcommand{\Y}{\mathcal{Y}}

\newcommand{\bg}{\mathcal{G}}
\def\bx{{\boldsymbol{x}}}
\def\be{{\boldsymbol{e}}}
\def\by{{\boldsymbol{y}}}
\def\bff{\mathcal{F}}
\def\bh{\mathcal{H}}
\def\bz{{\boldsymbol{z}}}
\def\bw{{\boldsymbol{w}}}
\def\bv{{\boldsymbol{v}}}
\def\bu{{\boldsymbol{u}}}

\newcommand{\supp}{\mathrm{supp}}

\def\bsa{{\boldsymbol{a}}}
\def\bd{{\boldsymbol{d}}}

\def\bu{{\boldsymbol{u}}}

\def\bv{{\boldsymbol{v}}}

\def\bgm{{\boldsymbol{\gamma}}}
\def\bgamma{{\boldsymbol{\gamma}}}
\def\bxi{{\boldsymbol{\xi}}}

\def\ny{\Gamma}

\newcommand{\ud}{\mathrm{d}}
\newcommand{\td}[1]{\tilde{#1}}
\newcommand{\suppx}[2]{\mbox{supp}_{#1}(#2)}

\newcommand{\rank}{\mbox{\upshape rank}}

\newcommand{\enpc}{\mbox{\upshape\tiny enpc}}
\newcommand{\cenpc}{\mbox{\upshape\tiny c-enpc}}
\newcommand{\nm}{\mbox{\upshape\tiny $n,m$}}
\newcommand{\subnm}[2]{\mbox{\upshape\tiny $#1$,$#2$}}
\newcommand{\sosc}{\mbox{\upshape\tiny sosc}}
\newcommand{\csosc}{\mbox{\upshape\tiny c-sosc}}

\newcommand{\psdp}{\mbox{\upshape\tiny primal}}
\newcommand{\dsdp}{\mbox{\upshape\tiny dual}}

\newcommand{\re}{\mathbb{R}}

\def\ba{{\boldsymbol{\alpha}}}
\def\bb{{\boldsymbol{\beta}}}
\def\bgamma{{\boldsymbol{\gamma}}}

\newcommand{\mX}{\mathcal{X}}

\newcommand{\mJ}{\mathfrak{J}}
\newcommand{\wt}[1]{\widetilde{#1}}

\newcommand{\mC}{\mathscr{P}}
\newcommand{\mQ}{\mathscr{Q}}
\newcommand{\Qk}{\mathscr{Q}_k}

\newcommand{\bM}{\mathbf{M}}

\newcommand{\X}{\mathcal{X}}

\newtheorem{assumption}{Assumption}[section]

\newtheorem{theorem}{Theorem}[section]

\newtheorem{corollary}{Corollary}[section]
\newtheorem{lemma}{Lemma}[section]
\newtheorem{example}{Example}[section]
\newtheorem{definition}{Definition}[section]
\newtheorem{remark}{Remark}[section]
\newtheorem{proposition}{Proposition}[section]

\title{
{Theorems of the Alternative and Nonnegativity Certificates for Systems of Jensen-type Polynomials with Universal Quantifiers}
}

\author{Feng Guo}
\address[Feng Guo]{School of Mathematical Sciences, Dalian University of Technology, Dalian 116024, Liaoning Province, China.}
\email{fguo@dlut.edu.cn}

\author{Yuhang Guo}
\address[Yuhang Guo]{School of Mathematical Sciences, Dalian University of Technology, Dalian 116024, Liaoning Province, China.}
\email{1484486158@mail.dlut.edu.cn}

\author{Guoyin Li}
\address[Guoyin Li]{Department of Applied Mathematics, University of New South Wales, 2052, Sydney, Australia}
\email{g.li@unsw.edu.au}

\thanks{The first and second authors were supported by the National Natural Science Foundation of China (grant no. 12471478). This research was also supported in part by the Australian Research Council through Discovery Project grants DP250101112 and DP250100004}

\begin{document}

\begin{abstract}
We introduce Jensen-type polynomials, a class of possibly nonconvex polynomials that extends the Jensen inequality for SOS-convex polynomials while retaining natural sums of squares (SOS) certificates for nonnegativity. We show that polynomials with essentially nonpositive coefficients belong to this class and establish a composition principle that yields further examples, including componentwise univariate compositions and compositions with structured quadratic maps motivated by source-localization loss functions. We then derive theorems of the alternative that provide SOS certificates for infeasibility of systems of strict inequalities involving Jensen-type polynomials and universally quantified variables. We illustrate our results through two applications, with numerical examples. First, we show that source-localization problems admit an exact SOS reformulation under a simple dimension condition, and that an optimal solution can be recovered by solving a semidefinite program (SDP) and applying the explicit Jensen map. Second, our SOS certificates yield two {\it convergent} SDP hierarchies for semi-infinite optimization with Jensen-type polynomials, providing upper and lower bounds on the optimal value. We also establish finite-convergence and solution-extraction results, under flat truncation rank condition and the associated Jensen maps. 
\end{abstract}
%\subjclass{}

\maketitle
\section{Introduction}
{
A fundamental problem in real algebraic geometry is to certify whether a
polynomial is nonnegative or positive on a given set. For real polynomials
$f_0,f_1,\ldots,f_J\in\RR[\bx]$ on $\mathbb{R}^n$, where $J$ is a nonnegative integer and 
$\RR[\bx]:=\RR[x_1,\ldots,x_n]$, consider the associated  semialgebraic set
\[
\K\coloneqq
\left\{\bx\in\RR^n\mid f_j(\bx)\leq 0,\ j=1,\ldots,J\right\}.
\]
Positivstellens\"atze provide algebraic certificates for the nonnegativity
or positivity of $f_0$ on $\K$. They are central to real algebraic
geometry and have important applications in optimization and control
theory
\cite{SOCAGbook,realAG,HKL2020,Lasserre2015,Laurent_sumsof,PPSOSMarshall,Nie2023,ScheidererSurvey}.

Classical Positivstellens\"atze provide certificates under different
assumptions on $\K$. For a compact polyhedron with nonempty interior,
Handelman's theorem \cite{Handelman1988} represents a polynomial strictly
positive on $\K$ as a positive linear combination of products of the
defining affine polynomials. For a compact basic semialgebraic set,
Schm\"udgen's theorem \cite{Schmugen1991} gives a representation in the
preordering generated by $-f_1,\ldots,-f_J$, while Putinar's theorem
\cite{Putinar1993} gives a representation in the corresponding quadratic
module when it is Archimedean. This condition implies compactness and, for
compact $\K$, can often be enforced by adding a redundant ball constraint.
Other results cover complementary settings: the Krivine--Stengle theorem
applies without compactness, while P\'olya's and Reznick's theorems and
extensions by Dickinson--Povh and Putinar--Vasilescu treat homogeneous
polynomials and related semialgebraic sets
\cite{DP2015,Krivine1964,Stengle1974,polya,Reznick1995,PV1999}. In
particular, Krivine--Stengle gives a general characterization of strict
positivity: $p>0$ on $\K$ if and only if
$\psi p=1+\phi$ for some $\psi$ and $\phi$ in the preordering generated by
$-f_1,\ldots,-f_J$.

These positivity certificates underpin the moment sums-of-squares (SOS) approach to
polynomial optimization, and were used successfully to  develop hierarchies of
semidefinite relaxations whose bounds converge to the global optimum under
suitable assumptions, see for example, \cite{SOCAGbook,LasserreGlobal2001,Lasserre09,Lasserre2015,Laurent_sumsof,Nie2023}.

In many applications, uncertainty gives rise to 
constraints involving infinitely many inequalities. Consider feasible sets of the form
$\U\coloneqq
\left\{\bx\in\RR^n\mid
f_j(\bx,\by)\leq 0,\ \forall\,\by\in\Y, \, j=1,\ldots,J
\right\}$, 
where $f_j\in\RR[\bx,\by]$, $\by=(y_1,\ldots,y_\ny)$, and
$\Y\subseteq\RR^\ny$ is a closed set. The universally quantified variable
$\by$ represents uncertainty in the constraints, making positivity
certificates on $\U$ useful in robust and semi-infinite optimization. Recently, Hu, Klep, and Nie \cite{HKN2024} established a Positivstellensatz for
sets of this form. Given a finite nonnegative Borel measure $\nu$ with
$\operatorname{supp}(\nu)=\Y$ satisfying the multivariate Carleman
condition, they showed that a polynomial positive on $\U$ belongs to an
associated quadratic module, provided that module is Archimedean. The
representation combines an SOS polynomial with integrals, against $\nu$,
of SOS-weighted constraint polynomials. They also developed moment
characterizations and Moment-SOS hierarchies for semi-infinite
optimization.

A key source of exactness in these moment-SOS methods is the connection between
Jensen's inequality and SOS-convexity. SOS-convex polynomials form a subclass of convex polynomials introduced by Helton and Nie
\cite{SRCSHN}, which covers convex quadratic functions or convex separable polynomials. 
An appealing feature of SOS-convexity is that it can be checked by solving a  semidefinite programming
(SDP) problem  
\cite{SRCSHN}; in contrast, checking convexity of a general polynomial
is computationally difficult \cite{Ahmadi2013}. For SOS-convex polynomials, Jensen's
inequality extends from probability measures to normalized positive linear
functionals: the value of such a functional is bounded below by evaluation
at the vector of its first moments. Under suitable assumptions, this
property yields an exact first-order semidefinite relaxation and allows an
optimizer to be recovered from an optimal moment solution
\cite{LasserreConvexity2009}. This phenomenon highlights a fundamental connection between convexity, Jensen inequalities, moment problems and
exactness of semidefinite relaxations.

This motivates the central question of this paper: {\it how far this link between Jensen inequalities and moment methods can be extended
beyond convexity?} For a general nonconvex
polynomial, evaluation at the vector of first moments need not satisfy
Jensen's inequality. We introduce \emph{Jensen classes} and their
associated \emph{Jensen maps} to capture nonconvex settings in which an
analogous principle remains valid. More precisely, 
consider an appropriate class $\mathcal{B}[\bx]_m$ of normalized positive linear
functionals on the space of real polynomials  with degree at most $m$. A map $\Phi_{\mJ}$ is called a
\emph{Jensen map} for the polynomial class $\mJ$ if
\[
L(f)\geq f\bigl(\Phi_{\mJ}(L)\bigr)
\qquad
\text{for every }L\in\mathcal{B}[\bx]_m \text{ and every }f\in\mJ.
\]
% {\color{red}for an appropriate
% class of normalized positive linear functional, a map $\Phi_{\mJ}$ is a Jensen
% map for a polynomial class $\mJ$ if for any $L$ in this class and for any $f \in \mJ$, }
% \[
% L(f)\geq f\big(\Phi_{\mJ}(L)\big). 
% %\forall f\in\mJ.
% \]
The precise domain of functionals is specified in the relevant definitions.
A class admitting such a map is called a Jensen class, and its elements
are called Jensen-type polynomials. Importantly, the
associated Jensen map, in general, depends nonlinearly on the moments. For SOS-convex polynomials, the standard 
Jensen map is nothing but the linear natural moment map
$L\mapsto (L(x_1),\ldots,L(x_n))$.

Our first main example is the class of polynomials with essentially
nonpositive coefficients. Although these polynomials need not be convex,
their coefficient structure yields an explicit Jensen map built from
selected even moments. The key ingredient is a generalized moment
inequality for positive linear functionals, which bounds suitable mixed
even moments by products of pure moments. We also establish a composition principle. If an outer polynomial belongs
to a Jensen class and the image of a polynomial map has suitable geometric
properties---for example, if it is closed and convex---then the
composition inherits a Jensen-type property. This principle produces
further nonconvex examples, including componentwise univariate
compositions and compositions with structured quadratic maps. The latter
are motivated by polynomial loss functions arising in source localization,
an important engineering application.

We next consider a mixed system of Jensen-type polynomial inequalities with one universally quantified inequality: \footnote{For simplicity of notation, we state our results for the case in which only one inequality involves a universally quantified variable. The case of multiple universally quantified inequalities can be treated by extending the same arguments or by using Lagrange interpolation, at the cost of heavier notation.} 
$f_j(\bx)\leq 0,\ j=0,\ldots,J,\ 
f(\bx,\by)\leq 0\ \text{for all }\by\in\Y$, 
where $\Y$ is a compact set, $f_j$, $j=0, 1,\ldots,J$, and for each fixed $\by \in \Y$, $f(\cdot,\by)$ are all {\it Jensen-type polynomials}. Under some suitable mild assumptions,   
we establish SOS-based theorems of the alternative for the strict feasibility of this mixed 
system. This then leads to an SOS certificate for non-negativity of $f_0$ over a feasible set $\mathcal X$ described by finitely many polynomial inequalities
and an uncountable family of polynomial inequalities: 
\[
\mathcal X
:=
\left\{
\bx\in\mathbb{R}^n
\ \middle|\
f_j(\bx)\leq 0,\ j=1,\ldots,J,\quad
f(\bx,\by)\leq 0,\ \forall\,\by\in\mathcal Y
\right\}.
\]

Finally, we provide two applications together with illustrating numerical examples. In the first application, we examine the source localization problem and illustrate the exact SOS reformulation and solution recovery procedure via the Jensen map for the point-source localization problem in dimension $n=3$, under an appropriate dimension condition. In the second application, we apply these certificates to semi-infinite 
optimization involving Jensen-type polynomials. We construct two convergent semidefinite programming
hierarchies that provide upper and lower bounds on the optimal value. We
also establish finite-convergence and solution-extraction results,
including results based on a flatness rank condition and the associated
explicit Jensen maps. Numerical examples illustrate the bounds,
convergence behavior, and optimizer extraction. 

Compared with the existing moment-SOS hierarchies for universally quantified constraints in \cite{HKN2024}, our lower-bound hierarchy  has a simpler and more compact structure. The
multipliers for the finitely many constraints are nonnegative scalars,
and the SOS multiplier for the universally quantified constraint
depends only on the quantified variable $\by$. In \cite{HKN2024}, the corresponding multipliers are SOS polynomials in $\bx$ and $(\bx,\by)$, respectively. Consequently, for the same relaxation order,
our lower-bound hierarchy admits an SDP formulation of smaller size, which can be less costly to solve in general.  The complementary upper-bound hierarchy, which was not considered in \cite{HKN2024}, provides an optimality gap, including when the rank condition for finite convergence is not satisfied.  

The main contributions of the paper are as follows:
\begin{enumerate}
\item We introduce Jensen classes of polynomials and their associated
nonlinear Jensen maps, and establish their basic properties.

\item We show that polynomials with essentially nonpositive coefficients
form a Jensen class with an explicit Jensen map. We also develop
composition principles that yield further nonconvex examples, including
componentwise univariate compositions and compositions with structured
quadratic maps.

\item We prove theorems of the alternative for systems of strict polynomial
inequalities involving universally quantified variables, and derive the
corresponding SOS certificates.

\item We illustrate our results through two applications, with numerical examples. First, we show that source-localization problems admit an exact SOS reformulation under a simple dimension condition, and that an optimal solution can be recovered by solving an SDP and applying the explicit Jensen map. Second,  we construct two convergent SDP hierarchies for
semi-infinite polynomial optimization. We establish finite-convergence and
solution-extraction results.
\end{enumerate}

The remainder of the paper is organized as follows.
Section~\ref{sec::pre} introduces notation and recalls basic facts about
positive linear functionals, sums of squares, and moment representations.
Section~\ref{sec::Jensen} defines Jensen classes and establishes the general
composition principle. Section~\ref{sec:classes} develops nonconvex
examples, including polynomials with essentially nonpositive coefficients
and polynomials with suitable composition structures.
Section~\ref{sec::alt} proves theorems of the alternative for
systems involving universal quantifiers.
Section~\ref{sec::certificates} derives the associated nonnegativity
certificates. Section~\ref{sec::sipop} presents two applications in source localization and semi-infinite optimization illustrating the results together with numerical examples.
Some conclusions are summarized in Section \ref{sec::conclusion} and technical proofs are collected in the appendices.}

%\vskip 7pt
%The rest of this paper is organized as follows. 
%In Section \ref{sec::pre}, we review necessary notation and preliminaries. 
%Section \ref{sec::alt} presents a theorem of the alternative for homogeneous polynomials with 
%the sign structures and UQs. 
%In Section \ref{sec::certificates}, we establish nonnegativity certificates for polynomials 
%with sign structures and UQs. 
%Section \ref{sec::sipop} applies these nonnegativity certificates to problem \eqref{eq::sipop}.

\section{Notation and Preliminaries}\label{sec::pre} 
We collect some notation and basic concepts
which will be used in this paper. 
The symbol $\N$ (resp., $\RR$) denotes
the set of nonnegative integers (resp., real numbers). 
For $n_1,n_2 \in\N$ with $n_1 \le n_2$, we  
let $$[n_1:n_2]\coloneqq\{n_1, \ldots,n_2\}.$$ 
For $\ba=(\alpha_i)\in\N^n$, let $|\ba|=\sum_{i=1}^n\alpha_i$. 
For $k\in\N$, let $\N^n_{=k}\coloneqq\{\ba\in\N^n\mid |\ba|=k\}$
and $\N^n_k\coloneqq\{\ba\in\N^n\mid |\ba|\le k\}$.
%We use $|\cdot|$ to denote the cardinality of a set. 
For $\bu=(u_1,\ldots,u_n)\in \RR^n$, $\Vert \bu\Vert$ denotes the standard Euclidean norm of $\bu$.
For $t\in\re$, we use $\lceil t\rceil$ to denote the smallest integer that is not smaller than $t$.
For $i\in[1:n]$, denote by $\be_i$ the vector with a $1$ in the $i$-th coordinate and $0$'s elsewhere.
The symbol ${\bf 0}$ denotes the zero vector, whose dimension is clear from context.
For a set of vectors $\{\bv_1,\ldots,\bv_s\}\subset\RR^n$, the notation 
${\rm span}\{\bv_1,\ldots,\bv_s\}$ denotes the set of all linear combinations of them in $\RR^n$.
Denote by $\bS^n$ (resp. $\bS_+^n$) 
the $n\times n$ symmetric real matrix (resp. $n\times n$ positive semidefinite matrix) space. 
The symbol $\bv^T$ (resp., $C^T$) denotes the transpose of the vector $\bv$ (resp., the matrix $C$).
We denote by $\bx$ (resp., $\by$)
the $n$-tuple (resp., $\ny$-tuple) of variables $(x_1,\ldots,x_n)$ (resp.,
$(y_1,\ldots,y_\ny)$).
Let
$\RR[\bx]$ be the ring of multivariate polynomials in $\bx$,
and $\RR[\bx]_k$ denote the subset of polynomials of degree at most $k$. A polynomial 
$\phi\in\RR[\bx]$ can be written as $\phi=\sum_{\ba\in\N^n}\phi_{\ba}\bx^{\ba}$ with 
$\phi_{\ba}\in\RR$ and 
$\bx^{\ba}\coloneqq x_1^{\alpha_1}\cdots x_n^{\alpha_n}$. 
Let $\deg(\phi)$ denote the degree of the polynomial $\phi$.
{Throughout the paper, when evaluating monomials
$\bx^{\ba}=\prod_{i=1}^n x_i^{\alpha_i}$, we adopt the standard convention
$x_i^0=1$, including when $x_i=0$.}. The support of $\phi$ is defined by 
$\supp(\phi)\coloneqq\{\ba\in\N^n\mid \phi_{\ba}\ne0\}$.
For $A=\RR[\bx]$ or $A= \RR[\bx]_k$, we use  $A^*$ to denote the dual
space of linear functionals from $A$ to $\RR$.

For a compact set $\Y\subset\RR^\ny$, let $\mathcal{M}(\Y)$ denote
the space of finite signed Borel measures on $\Y$, and let
$\mathcal{M}_+(\Y)$ denote its cone of nonnegative measures. Let
$C(\Y)$ be the space of continuous functions from $\Y$ to $\RR$,
equipped with the supremum norm
$\|\gamma\|_\infty\coloneqq\sup_{\by\in\Y}|\gamma(\by)|$. 
By the Riesz representation theorem, $C(\Y)^*$ is identified with
$\mathcal{M}(\Y)$. Under this identification, the duality pairing is
\[
\langle\mu,\gamma\rangle
\coloneqq\int_{\Y}\gamma(\by)\,\ud\mu(\by),
\qquad
\mu\in\mathcal{M}(\Y),\quad \gamma\in C(\Y),
\]
and the dual norm of $\mu$ is its total variation
$|\mu|(\Y)$.

%\subsection{Sums of squares and moments}
We recall some background about SOS polynomials and the duality theory of moment matrices.
%For any $f(\bx)\in\RR[\bx]_k$, let  $\bf{f}$ denote its column vector of
%coefficients  in the canonical monomial basis of $\RR[Y]_k$. 
A polynomial $\phi(\bx)\in\RR[\bx]$ is said to be 
\emph{sum of squares} (SOS) if it can be written as $\phi(\bx)=\sum_{i=1}^t
\phi_i(\bx)^2$ for some $\phi_1(\bx),\ldots,\phi_t(\bx)\in\RR[\bx]$. 
%Let $\Sigma[\bx]$ denote the set of SOS polynomials in $\RR[\bx]$.
Denote by $\Sigma[\bx]$ the set of SOS polynomials in 
$\RR[\bx]$. 
Denote by $\Sigma[\bx]_m$ the set of SOS polynomials of degree at most $m$. Every SOS polynomial is nonnegative, but the converse does not hold
in general; see, for example, \cite{Reznick96someconcrete}. Hilbert's
theorem identifies several cases in which nonnegative polynomials are
SOS. Further recent results concerning SOS characterization for nonnegative polynomials with specific structure, for example, can be found at \cite{ADG2023}. 
\begin{theorem}[Hilbert's Theorem \cite{Lasserre09}]\label{th::hilbert}
	Every nonnegative polynomial $\phi\in\RR[\bx]$ on $\RR^n$ is SOS
	in the following cases: {\rm (i) $n=1$; \footnote{The fact that every nonnegative univariate polynomial is SOS extends to nonnegative separable polynomials. Here, a polynomial $\phi \in \RR[\bx]$ is separable if $\phi(\bx)=\sum_{i=1}^n\phi_i(x_i)$ for some univariate polynomial $\phi_i$.} (ii) $\deg(\phi) \le 2$; (iii) $n=2$}
	and $\deg(\phi) \le 4$.
\end{theorem}

%{\color{blue} Some background, such as SOS implies nonnegativity, but converse is not true. Recall the three cases SOS are equivalent to nonnegativity. Also recall SOS-convexity as well.}

%{\color{blue} Can we add some preliminaries for measure space, $C(Y)$ and their norms and $M_+(Y)$?} 

Let $\bh\coloneqq\{h_1,\ldots,h_r\}\subset\RR[\bx]$ be the set of polynomials, and consider
\[\K_{\bh}:=\{\bx\in\RR^n \mid h_1(\bx)\ge 0, \ldots, h_r(\bx)\ge 0\}.\] 
We denote by
\[
  \mQ(\bh)\coloneqq\left\{\sigma_0+\sum_{j=1}^r\sigma_jh_j\ \middle|\ 
\sigma_j \in \Sigma[\bx], \ \ j\in[0:r] \right\}
\]
the {\itshape quadratic module} generated by $\bh$ and denote by
\[
\Qk(\bh)\coloneqq\left\{\sigma_0+\sum_{j=1}^r\sigma_jh_j\ \middle|\ \sigma_0,\sigma_j \in \Sigma[\bx], \deg(\sigma_0),\deg(\sigma_jh_j)\le 2k,\ \ j\in[1:r]\right\}
\]
the $k$-th {\itshape truncated quadratic module}.
It is clear that if a polynomial $p\in \mQ(\bh)$, then
$p(\bx)\ge 0$ for any $\bx\in \K_{\bh}$ though the converse is not necessarily
true.
% Note that checking whether $f\in\mathcal{Q}_k(\bg)$ for a fixed $k\in\N$ is an
% SDP feasibility problem \cite{LasserreGlobal2001}.

% For $k\in\N$, denote $s(k)\coloneqq{n+k\choose n}$.
Given a (pseudo-moment) sequence of real numbers $\bz\coloneqq(z_{\ba})_{\ba\in \mathbb{N}^n}$,
% We call $\by$ a {\itshape truncated moment sequence} up to order $2k$ if there exists
% a Borel measure $\mu$ on $\RR^n$ such that
% \[y_{\ba}=\int \bx^{\ba}\ud\mu(\bx),\quad \forall \ba\in\N^n_{2k}.\]
% In this case, we say that $\by$ has a {\itshape representing measure} $\mu$.
the $k$-th order {\itshape moment matrix} is the matrix $\bM_k(\bz)$
indexed by $\mathbb{N}^n_{k}$ with the $(\ba,\bb)$-th entry being $z_{\ba+\bb}$. Given a polynomial $\phi(\bx)=\sum_{\ba\in\supp(\phi)}\phi_{\ba}\bx^{\ba}$, the $k$-th order {\itshape localizing
matrix} $\bM_{k}(\phi\bz)$ indexed by $\mathbb{N}^n_{k}$ is defined by $[\bM_{k}(\phi\bz)]_{\bb,\bgm}=\sum_{\ba\in\supp(\phi)}\phi_{\ba}z_{\ba+\bb+\bgm}$.
% Denote by ${\mathscr{M}_{2k}}$ the space of all sequences
% $\by=(y_{\ba})_{\ba \in \mathbb{N}^n_{2k}}\in\RR^{s(2k)}$ with
% order at most $2k$.
The Riesz functional $L_{\bz}$ on $\RR[\bx]$ is defined by $L_{\bz}\left(\phi\right)
\coloneqq\sum_{\ba\in\supp(\phi)}\phi_{\ba}z_{\ba}$.
% From the definition of the localizing moment
% matrix $\bM_{k-d_f}(f\by)$, it is easy to check that
% \begin{equation}\label{eq::M}
% \mathbf{h}^{\intercal}\bM_{k-d_f}(f\by)\mathbf{h}=\mathscr{L}_{\by}(f(\bx)h(\bx)^2),\quad \forall
% h(\bx)\in\RR[\bx]_{k-d_f}.
% \end{equation}
% Let $d_j\coloneqq\lceil\deg(g_j)/2\rceil$ for each $j=1,\ldots,s$.
% For any $\bv\in \K$, let $\zeta_{2k,\bv}\coloneqq[v^{\ba}]_{\ba\in\N^n_{2k}}$
% be the {\itshape Zeta vector} of $\bv$ up to degree $2k$, i.e.,
% \[
% 	\zeta_{2k,\bv}=
% 	[1\quad v_1\quad \cdots\quad v_n\quad
% 		v_1^2\quad v_1v_2\quad \cdots\quad v_n^{2k}].
% \]
% Then, $\bM_k(\zeta_{2k,\bv})\succeq 0$ and $\bM_{k-d_j}(g_j\zeta_{2k,\bv})\succeq 0$
% for $j=1,\ldots,s$. In fact, let $g_0=1$, then for each $j=0,1,\ldots,s$,
% \[
% 	\begin{aligned}
% 		\mathbf{h}^{\intercal}\bM_{k-d_j}(g_j\zeta_{2k,\bv})\mathbf{h}
% 		&=\mathscr{L}_{\zeta_{2k,\bv}}(g_j(\bx)h(\bx)^2)=g_j(\bv)h(\bv)^2\ge 0,\ \forall
% 		h(\bx)\in\RR[\bx]_{k-d_j}.
% 	\end{aligned}
% \]
\begin{definition}[Archimedean condition]\label{def::AC}
We say that $\mQ(\bh)$ is {\itshape Archimedean} if 
there exists $M>0$ such that $M-x_1^2-\cdots-x_n^2\in \mQ(\bh)$.
\end{definition}
Note that the Archimedean condition implies that $\K_{\bh}$ is compact and the converse is not necessarily true. However, for any compact set $\K_{\bh}$, we could always
force the associated quadratic module to be Archimedean
by adding a redundant constraint $M-x_1^2-\cdots-x_n^2\ge 0$ in the description of
$\K_{\bh}$ for sufficiently large $M$.

\begin{theorem}\label{th::PP}{\upshape\cite[{Putinar's Positivstellensatz\/}]{Putinar1993}}
Suppose that $\mQ(\bh)$ is Archimedean.
If a polynomial $\phi\in\RR[\bx]$ is positive on $\K_{\bh}$, then $\phi\in\Qk(\bh)$ for some $k\in\N$.
\end{theorem}

We also recall SOS-convexity, introduced by Helton and Nie
\cite{SRCSHN} (see \cite{ahmadi2013complete} for a further study). A polynomial $\phi\in\RR[\bx]$ is SOS-convex if there
exist an integer $r$ and a polynomial matrix
$\Phi\in\RR[\bx]^{r\times n}$ such that $\nabla^2\phi(\bx)=\Phi(\bx)^\mathsf{T}\Phi(\bx)$. 
SOS-convexity can be checked by solving an SDP problem 
\cite{SRCSHN}; in contrast, checking convexity of a general polynomial
is computationally difficult \cite{Ahmadi2013}. It is known that any convex quadratic function or convex separable polynomial is an SOS-convex polynomial, or more generally, sum of a convex separable polynomial and convex quadratic function.  Importantly, SOS-convex
polynomials also satisfy the following Jensen inequality. For a positive integer $m$, 
let $(\Sigma[\bx]_{m})^*$ be the positive dual cone of $\Sigma[\bx]_{m}$, that is, 
\[
(\Sigma[\bx]_{m})^*=\{L \mbox{ is a continuous linear functional on } \RR[\bx]_{m}\ \mid\ 
L(f) \ge 0 \mbox{ for all } f \in \Sigma[\bx]_{m}\}.
\]

\begin{proposition}{\upshape\cite[Theorem 2.6]{LasserreConvexity2009}}\label{prop::Jensen}
	Let $\phi\in\RR[\bx]_{m}$ be an SOS-convex polynomial. For any
	$L\in(\Sigma[\bx]_{m})^*$ satisfying $L(1)=1$, it holds that  
    $L(\phi(\bx))\ge
	\phi(L(x_1),\ldots,L(x_n)).$
\end{proposition}
%The following result plays a significant role
%in polynomial optimization with SOS-convex polynomials.
%\begin{lemma}{\upshape\cite[Lemma 8]{SRCSHN}}\label{lem::sosconvex}
%	Let $\phi\in\RR[\bx]$ be SOS-convex. If $\phi(\bu)=0$ and $\nabla \phi(\bu)=0$
%	for some $\bu\in\RR^n,$ then $\phi$ is an SOS polynonmial.
%\end{lemma}

\section{Jensen-type polynomials and Jensen class}\label{sec::Jensen}
In this section, we first introduce the definition and basic properties of Jensen-type polynomials.
Then we develop general composition principles that yield further nonconvex examples of Jensen
classes.
\subsection{Definition and basic properties}
{ Let $m$ be a positive integer. Recall that $\Sigma[\bx]_m$ is 
the cone of SOS polynomials in $\RR[\bx]_m$ and $(\Sigma[\bx]_m)^*$ is the positive dual cone of $\Sigma[\bx]_m$.
%, that is, 
%$$(\Sigma[\bx]_m)^*=\{L \mbox{ is a continuous linear functional on } \RR[\bx]_m:  L(f) \ge 0 \mbox{ for all } f \in \Sigma[\bx]_m\}.$$
Let $\mathcal{B}[\bx]_m$ be defined by
$$
\mathcal{B}[\bx]_m:=\{L\in (\Sigma[\bx]_m)^* \mid L(1)=1\}.
$$ We now introduce the notion of Jensen class and Jensen-type polynomial as follows. 
%\begin{definition}
% Let $\mJ\subset\RR[\bx]_m$ with a map $\Phi_{\mJ}: \overline{(\Sigma[\bx]_m)^*} \to \RR^n$. We say $f$   We call $\mJ\subset\RR[\bx]_m$ with a map $\Phi_{\mJ}: \overline{(\Sigma[\bx]_m)^*} \to \RR^n$
%    satisfying the condition
%    \[
%    L(f)\ge f(\Phi_{\mJ}(L)),\quad \forall\ f\in\mJ,\ L\in\overline{(\Sigma[\bx]_m)^*},
%    \tag{C}
%    \]
%    a set of Jensen-type polynomials.
%\end{definition}} 

{\begin{definition}[Jensen class \& Jensen-type polynomials]
Let $\mJ\subseteq\mathbb{R}[\bx]_m$ be a set of real polynomials with degree at most $m$. We say $\mJ$  is a \emph{Jensen class} if there exists a mapping
$
\Phi_{\mJ}:\mathcal{B}[\bx]_m\to\mathbb{R}^n
$
such that for all $f\in\mJ$ and $L\in\mathcal{B}[\bx]_m$, 
\begin{equation}\label{J}
L(f)\geq f\big(\Phi_{\mJ}(L)\big).
\qquad
\tag{C}
\end{equation}
For a Jensen class $\mJ$, a mapping $\Phi_{\mJ}$ satisfying \eqref{J} is called a \emph{Jensen map} associated with $\mJ$. Each $f\in\mJ$ is called a \emph{Jensen-type polynomial} with respect to $(\mJ,\Phi_{\mJ})$. 

Note that the condition \eqref{J} automatically holds with equality when $f$ is a constant polynomial. Also,  if $f$ satisfies \eqref{J} with some Jensen map $\Phi_{\mJ}$, then, for any $\alpha \in \mathbb{R}$,  $f+\alpha$ also satisfies \eqref{J} with the same Jensen map. Throughout this paper, {\emph  without loss of generality, we always assume that $\mJ$ contains all the constant polynomials, and is closed under addition with any constant polynomial.}
\end{definition}}

%{\color{blue} In the definition of Jensen class, the requirement $0 \in \mJ$ and $\mJ + \{p_{\alpha} \mid p_{\alpha} \equiv \alpha, \, \alpha \in \mathbb{R}\} =\mJ$ ensure that $\mJ$ contains every constant polynomials (and so, is nontrivial), and is closed under addition of constants. In particular, $f \in \mJ$ implies that $f +\alpha \in \mJ$ for any $\alpha \in \mathbb{R}$}. 
Denote by $\mJ_{\nm}^{\sosc}$ the set of SOS-convex polynomials in
$\mathbb{R}[\bx]_m$. 
By Proposition \ref{prop::Jensen}, the set
$\mJ_{\nm}^{\sosc}$ together with the map
\[
\Phi_{\mJ_{\nm}^{\sosc}}:
\mathcal{B}[\bx]_m \to \mathbb{R}^n,
\qquad
L \mapsto \big(L(x_1),\ldots,L(x_n)\big),
\]
satisfies condition~\eqref{J}. Here,
$\Phi_{\mJ_{\nm}^{\sosc}}$ is the natural moment map: its
$i$th component is the image of the coordinate function $x_i$
under the linear functional $L$. Consequently,
$\mJ_{\nm}^{\sosc}$ is a Jensen class with respect to the natural
Jensen map $\Phi_{\mJ_{\nm}^{\sosc}}$.

\begin{remark}[Comments on the Jensen map] 
For the SOS-convex polynomials, the associated natural Jensen map is a linear map. On the other hand, we note that, in the definition (C) of Jensen class and Jensen-type polynomial, the associated Jensen map $\Phi_{\mJ}$ need not be a linear map. As we will see later in Section~\ref{sec:classes}, this flexibility allows us to cover interesting nonconvex classes of polynomials which admit tractable SOS representations. 
\end{remark}

A simple consequence of the definition of a Jensen-type polynomial is that its nonnegativity admits an SOS
%sums-of-squares 
representation.  
\begin{lemma}\label{cor::sos}
    Let $\mJ\subset\RR[\bx]_m$ be a Jensen class with a Jensen map $\Phi_{\mJ}$.
    For any $f\in\mJ$, if $f^{\star}:=\inf_{\bx\in\RR^n}f(\bx)>-\infty$,
    then $f(\bx)-f^{\star}\in\Sigma[\bx]_m$. In particular, for any $f \in \mJ$, $f$ is nonnegative if and only if $f \in \Sigma[\bx]_m$.
\end{lemma}
\begin{proof} It suffices to establish the first statement as the second statement follows immediately. To do this, let $f^{\star}:=\inf_{\bx\in\RR^n}f(\bx)>-\infty$. Suppose to the contrary, that $f-f^{\star}\not\in\Sigma[\bx]_m$. 
    Since $\Sigma[\bx]_m$ is 
    closed and convex, by the separation theorem, there exists a linear functional
    $L\in(\Sigma[\bx]_m)^*$ such that $L(f-f^{\star})<0$,
    i.e., $L(f)<f^{\star} \, L(1)$. 
    Since $1\in \Sigma[\bx]_m$, we have $L(1)\ge 0$. We now divide the discussion into two cases.

    Case 1. $L(1)>0$. Let $\wt{L}=L/L(1)$. 
    Then, $\wt{L}\in\mathcal{B}[{\bx}]_m$.
    As $f(\bx)\ge f^{\star}$ for all $\bx \in \mathbb{R}^n$, for a Jensen map $\Phi_{\mJ}$ associated with $\mJ$, it holds that 
    \[
    f^{\star}\le f(\Phi_{\mJ}(\wt{L}))\le \wt{L}(f)=L(f)/L(1)<f^{\star},
    \]
    which is not possible, and so, a contradiction.

    Case 2.  $L(1)=0$. Fix an $L_0\in\mathcal{B}[{\bx}]_m$. 
%{\color{red}    If $L_0(f-f^{\star})<0$, then we can get the same contradiction
%    \[
%    f^{\star}\le f(\Phi_{\mJ}(L_0))\le L_0(f)<L_0(f^\star)=f^{\star}.
%    \]
%    Hence, we assume that $L_0(f-f^{\star})\ge 0$. (It seems that we don't need this part)} 
    As $L(f-f^{\star})<0$, there exists 
$\varepsilon>0$ such that $(L+\varepsilon L_0)(f-f^{\star})<0$. 
    Consider the linear functional
    $L_{\varepsilon}:=L+\varepsilon L_0$. We have $L_{\varepsilon}\in(\Sigma[\bx]_m)^*$ with $L_{\varepsilon}(1)=L(1)+\epsilon L_0(1)= \epsilon >0$.
    Then repeating the above arguments in the first case
    yields again a contradiction.
\end{proof}

\subsection{A general composition principle}
Next, we provide results which allows one to recognize Jensen class and Jensen-type polynomial with composition structure. To do this, we define the push-forward map induced by a polynomial mapping as follows. Let $m_1,m_2$ be nonnegative integers, and let $m=m_1m_2$. Let $h=(h_1,\ldots,h_p)$ where each $h_\ell\in\RR[\bx]_{m_2}$, and denote by $\bxi=(\xi_1,\ldots,\xi_p)$ a $p$-tuple of variables. Then, for each ${L} \in \mathcal{B}[\bx]_m$, we define the push-forward mapping $h_{\sharp}L: \RR[\bxi]_{m_1} \rightarrow \mathbb{R}$ induced by the polynomial mapping $h$ as  
\begin{equation}\label{eq:push-forward}
(h_{\sharp}L)(q)=L(q \circ h) \mbox{ for all } q \in \RR[\bxi]_{m_1}.    
\end{equation}

We first see that, under a compatible condition, composition of a Jensen class with a polynomial mapping induces another Jensen class. 
\begin{lemma}[Compositions under compatible conditions]  \label{lem::compatible_conditions}
Let $\mJ\subset\RR[\bxi]_{m_1}$ be a Jensen class of real polynomials on $\RR^p$ with Jensen map $\Phi_\mJ:\mathcal{B}[\bxi]_{m_1} \rightarrow \RR^p$. %Let $g \in \mJ$ and  
%let $f=g \circ h$ where  
Let $h=(h_1,\ldots,h_p)$ with $h_\ell\in\RR[\bx]_{m_2}$, $\ell\in[1:p]$.  
%$:\RR^n \rightarrow \RR$, $\ell=1,\ldots,p$,  and $h_\ell$ are all real polynomials with degree $m_2$. 
Let $m=m_1m_2$. Suppose that  
$\Phi_{\mJ}(h_\sharp L) \in h(\mathbb{R}^n) \mbox{ for all } {L} \in \mathcal{B}[\bx]_m$, that is, 
\[ h^{-1}(\Phi_{\mJ}(h_\sharp L)) \neq \emptyset \mbox{ for all } {L} \in \mathcal{B}[\bx]_m.
\]Then,  $\widehat{\mJ}=\{f \in \RR[\bx]_m \mid f=g \circ h, \ g \in \mJ\}$ is a Jensen class with Jensen map $\Phi_{\widehat{\mJ}}$, where $\Phi_{\widehat{\mJ}}:\mathcal{B}[\bx]_m \rightarrow \mathbb{R}^n$ is defined in such a way  that  \[
\Phi_{\widehat{\mJ}}(L) \in h^{-1}(\Phi_{\mJ}(h_\sharp L)).
\]
\end{lemma}
\begin{proof}
As $g \in \mJ$ and $\mJ$ is a Jensen class of real polynomials on $\RR^p$ with Jensen map $\Phi_\mJ:\mathcal{B}[\bxi]_{m_1} \rightarrow \RR^p$, for all ${L}_0 \in \mathcal{B}[\bxi]_{m_1}$, one has  
$L_0(g) \ge g(\Phi_\mJ(L_0))$. 
For each ${L} \in \mathcal{B}[\bx]_m$, recall that $h_{\sharp}L: \RR[\bxi]_{m_1} \rightarrow \mathbb{R}$ is given as 
\[
(h_{\sharp}L)(q)=L(q \circ h) \mbox{ for all } q \in \RR[\bxi]_{m_1}.
\]
Then, direct verification shows that, for all ${L} \in \mathcal{B}[\bx]_m$, $h_{\sharp}L(1)=1$ and $h_{\sharp}L \in (\Sigma[\bxi]_{m_1})^*$, and so,  
$h_{\sharp}L \in \mathcal{B}[\bxi]_{m_1}$.
Thus, for all $L \in \mathcal{B}[\bx]_m$ and $f \in \widehat{\mJ}$, there exists $g \in \mJ$ such that $f=g\circ h$ and so, 
\[
L(f)=L( g \circ h)=h_{\sharp}L(g) \ge g(\Phi_\mJ(h_{\sharp}L)).
\]
Now, with the assumption that $\Phi_\mJ(h_{\sharp}L) \in h(\mathbb{R}^n)$, there exists $\bu=(u_1,\ldots,u_n) \in \RR^n$ such that 
\[
h(u_1,\ldots,u_n)=\Phi_\mJ(h_{\sharp}L).
\]
For each $L\in\mathcal{B}[\bx]_m$, fix such a point $\bu$ and define 
$\Phi_{\widehat{\mJ}}(L)= \bu$. 
Then, we have 
\[
L(f)=L( g \circ h)=h_{\sharp}L(g) \ge g(\Phi_\mJ(h_{\sharp}L)) = g(h(\bu)) =f(\bu)= f(\Phi_{\widehat{\mJ}}(L)).
\]
Thus, the conclusion follows. 
\end{proof}

%{\color{red}}
%\[
%\sqrt[m]{L(y_1^m \circ h)}=\sqrt[m]{L(h_1^m)}
%\]

%}
The next lemma provides a convenient check for the compatibility condition under an image convexity assumption. We note that image convexity for quadratic maps is an important topic in geometry and closely relates to the concept of joint numerical range. 
\begin{lemma}\label{lemma:3.3}
Let $m_1,m_2$ be positive integers, and let $m=m_1m_2$. 
Let $\Phi_\mJ:\mathcal{B}[\bxi]_{m_1} \rightarrow \mathbb{R}^p$ be the natural moment map, defined by, 
\[
\Phi_\mJ(L_0)= (L_0(\xi_1),\ldots,L_0(\xi_p)), \ \forall \, L_0 \in \mathcal{B}[\bxi]_{m_1}.  
\]
Here, for each $\ell\in[1:p]$, $L_0(\xi_\ell)$ denotes the image of
the coordinate function 
$\bxi\mapsto \xi_\ell$ 
under the linear functional $L_0$. Let $h_\ell \in \RR[\bx]_{m_2}$, $\ell\in[1:p]$, and $h=(h_1,\ldots,h_p)$. Suppose that  $h(\RR^n)$ is a closed and convex set, and one of the following conditions
holds:
\begin{itemize}
\item[(a)] 
each $h_\ell$ is a separable polynomial, $\ell\in[1:p]$; 
\item[(b)] each $h_\ell$ is a quadratic function, $\ell\in[1:p]$;
\item[(c)] each $h_\ell$ is a bivariate quartic polynomial on $\bx=(x_1,x_2) \in \mathbb{R}^2$, $\ell\in[1:p]$.
\end{itemize}
Then, $\Phi_{\mJ}(h_\sharp L) \in h(\RR^n) \mbox{ for all } {L} \in \mathcal{B}[\bx]_m$. 
\end{lemma}
\begin{proof}
Note that 
\[
\begin{aligned}
    \Phi_{\mJ}(h_\sharp L) =(h_\sharp L(\xi_1),\ldots,h_\sharp L(\xi_p))
    =(L(\xi_1 \circ h),\ldots,L(\xi_p \circ h))
    =(L(h_1),\ldots,L(h_p)).
\end{aligned}
\]
To the contrary, suppose that $(L(h_1),\ldots,L(h_p)) \not\in h(\RR^n)$.
%see the conclusion, it suffices to show that 
%\[
%(L(h_1),\ldots,L(h_p)) \in K.
%\]
%Otherwise, by 
By the strict convex separation theorem and the assumption that $h(\RR^n)$ is closed and convex, 
there exist $(\lambda_1,\ldots,\lambda_p) \neq (0,\ldots,0)$ and $\gamma \in \RR$ such that 
\[
L\left(\sum_{\ell=1}^p \lambda_\ell h_\ell\right) =\sum_{\ell=1}^p \lambda_\ell L( h_\ell) < \gamma \le \sum_{\ell=1}^p \lambda_\ell h_{\ell}(\bx) \mbox{ for all }  \bx \in \RR^n. 
\]
Let $w(\bx)=\sum_{\ell=1}^p \lambda_\ell h_\ell(\bx)-\gamma$. Then, as $L(1)=1$, we have 
$L(w) < 0$ but $w(\bx)\ge 0$ for all $\bx \in \RR^n$. Next, we show that $w \in \Sigma[\bx]_m$. To see this, if each $h_\ell$ is separable, then $w$ is also   separable. As a nonnegative function, it must be SOS. Similarly, if each $h_\ell$ is quadratic (resp., bivariate quartic polynomial on $\bx=(x_1,x_2)$), then $w$ is a nonnegative quadratic function (resp., bivariate quartic polynomial on $\bx=(x_1,x_2)$), and so, is also SOS (Theorem \ref{th::hilbert}). As $L \in \mathcal{B}[\bx]_m$ and $w \in \Sigma[\bx]_m$, $L(w) \ge 0$ which contradicts the fact that $L(w)<0$. Thus, the conclusion follows.  
\end{proof}

We note that, as the following two examples illustrate,  the compatibility condition can fail in general if the closedness or convexity assumption of the image set $h(\mathbb{R}^n)$  is dropped. 
\begin{example}
Consider $h=(h_1,h_2)$ where $h_1(x)=x$ and $h_2(x)=x^2$ with $x \in \RR$. 
Then, $h(\RR)=\{(t,t^2): t \in \RR\}$ which is not convex. Let $m_2=\max\{{\rm deg}(h_1),{\rm deg}(h_2)\}=2$. For any positive integer $m_1$, set $m=m_1\, m_2=2m_1$. 
Consider the linear functional {$L \in \mathcal{B}[x]_m$} given by 
$L(p)=\frac{1}{2}(p(1)+p(-1))$. 
Direct verification shows that $L(h_1)=0$ and $L(h_2)=1$. This implies that the compatibility condition fails because 
$\Phi_{\mJ}(h_\sharp L)= (L(h_1),L(h_2))=(0,1) \notin h(\RR)$. 
%{\color{red} (Here, $\Phi_{\mJ}$ denotes the natural
%moment map on $\mathcal{B}[\xi_1,\xi_2]_{m_1}$. It seems that we set $m_2=m=2$, so $m_1=1$?)}
\end{example}

\begin{example}
Consider $h=(h_1,h_2)$ where $h_1(\bx)=x_1$ and $h_2(\bx)=x_1x_2$ with $\bx=(x_1,x_2) \in \RR^2$. Then, $h(\RR^2)=\{(0,0)\} \cup \big((\mathbb{R} \backslash\{0\}) \times \mathbb{R}\big)$ which is not closed. Let $m_2=\max\{{\rm deg}(h_1),{\rm deg}(h_2)\}=2$. For any positive integer $m_1$, set $m=m_1\, m_2=2m_1$. 
Consider the linear functional {$L \in \mathcal{B}[x]_m$} defined by 
$L(p)=\frac{1}{2}(p(1,1)+p(-1,-1))$. 
Direct verification shows that $L(h_1)=0$ and $L(h_2)=1$. This implies that the compatibility condition fails because 
$\Phi_{\mJ}(h_\sharp L)= (L(h_1),L(h_2))=(0,1) \notin h(\RR^2)$. %{\color{red}(same question as above)}
\end{example}

\section{Classes of nonconvex Jensen-type polynomials}\label{sec:classes}
In this section, we provide various classes of nonconvex Jensen-type  polynomials. 
\subsection{Polynomials with ENPCs}\label{sec::epnc}
We first see that polynomials with essentially non-positive coefficients (ENPC) form a Jensen class with a specific nonlinear Jensen map. A polynomial with ENPC (or simply ENPC polynomial) is a polynomial with even degree $m$, whose coefficients are all nonpositive except those associated
with the monomials $x_i^m$ and the constant term. In the homogeneous case, polynomial of this form arises naturally in studying the maximum eigenvalue of the tensors with nonnegative entries, which has applications in high-order Markov chains \cite{HLQ2016}. Moreover, in the homogeneous quadratic setting, this requirement means the Hessian matrix of the quadratic function has non-positive off-diagonal elements (often referred to as $Z$-matrix), which has been studied in nonconvex quadratic optimization \cite{li,kim2003exact}.   

%\subsection{Subclass 1. Polynomials with ENPCs}
Let $m$ be an even integer. Consider the set $\mJ_{\nm}^{\enpc}$ of polynomials with ENPCs in $\RR[\bx]_m$, i.e.,  
\[
\mJ_{\nm}^{\enpc}:=\left\{\sum_{i=1}^n h_{m,i}x_i^m+\sum_{\ba\in\Omega\setminus\{\bf 0\}} h_{\ba}\bx^{\ba}+h_{\bf 0}
\ \Bigg|h_{m,i}, \ h_{\ba},\ h_{\bf 0}\in\RR \mbox{ and }  h_{\ba}\le 0, \forall \ba\in\Omega\setminus \{\bf 0\}
\right\},
\]
where $\Omega:= \N^n_m \backslash \{m\be_1, \ldots, m\be_n\}$, 
and the map
\begin{equation}\label{eq::enpc_map}
    \Phi_{\mJ_{\nm}^{\enpc}}:  \mathcal{B}[\bx]_m \to \RR^n, \quad L\mapsto 
\left(\sqrt[m]{L(x_1^m)},\ldots, \sqrt[m]{L(x_n^m)}\right). 
\end{equation}
Here, $L(x_i^m)$ denotes the value of the linear functional $L$ at the
polynomial $\bx\mapsto x_i^m$. Since $m$ is even and
$L \in \mathcal{B}[\bx]_m$, one has $L(x_i^m)\geq0$, and hence the $m$th-roots in
\eqref{eq::enpc_map} are well-defined.

The main result relies on the following Hölder-type inequality for even
moments, whose proof is given in the appendix.

%Consider the convex set
%\begin{equation}\label{eq::W}
%    \mathcal{W}:=\left\{\bw\in\RR^{|\N^n_{m}|} \mid w_0=1,\ \bM_{m/2}(\bw)\succeq 0\right\}.
%\end{equation}
%We write $\be=(\be_1,\ldots,\be_n)$ for the standard basis of $\RR^n$.  
%For any point $\bw=(w_{\ba})_{\ba\in\N^n_m}\in\RR^{|\N^n_m|}$, 
%let 
%\[
%\bx_{\bw}:=\left(\sqrt[m]{w_{m\be_1}},\ldots,\sqrt[m]{w_{m\be_n}}\right)\in\RR^n.
%\]

\begin{lemma}\label{lemma:claim}
For any $s\in\mathbb{N}$ with $s>0$, any
$L\in\mathcal{B}[\bx]_{2s}$, and any
$\bv\in\mathbb{N}^n$ satisfying $|\bv|\le s$, one has
\[
L(\bx^{2\bv})
\le
\prod_{i=1}^n L(x_i^{2s})^{v_i/s}.
\tag{C'}
\]
\end{lemma}
\begin{proof}
See Appendix \ref{sec::proof0}.
\end{proof}
\begin{proposition}[Subclass 1: Polynomials with ENPCs]\label{lem::extremal}
Let $m$ be an even integer. Consider the set $\mJ_{\nm}^{\enpc}$ of polynomials with ENPCs in $\RR[\bx]_m$ and the map $\Phi_{\mJ_{\nm}^{\enpc}}$. 
\begin{itemize} 
\item[{\rm (i)}]  For any $L\in \mathcal{B}[\bx]_m$, it holds that
\begin{equation}\label{eq::ineq_mom}
\left(\Phi_{\mJ_{\nm}^{\enpc}}(L)\right)^{\ba}=\prod_{i=1}^n\left(\sqrt[m]{L(x_i^m)}\right)^{\alpha_i}\ge |L(\bx^{\ba})|,\quad\forall \ba\in\N^n_m,
\end{equation}
and the equality holds for $\ba=m\be_i$, $i\in[1:n]$ or $\ba ={\bf 0}$.
%    The set of extremal points of $\mathcal{W}_k$
%    the convex set 
%\begin{equation}\label{eq::W}
%    \left\{\bw\in\RR^{|\N^n_{2k}|} \mid w_0=1,\ \bM_k(\bw)\succeq 0\right\}
%\end{equation}
%is $\{(\bu^{\ba})_{\ba\in\N^n_{2k}} \mid \bu\in\RR^n\}$.

\item[{\rm (ii)}]  The set of polynomials $\mJ_{\nm}^{\enpc}$ is a Jensen class with the Jensen map $\Phi_{\mJ_{\nm}^{\enpc}}$.
\end{itemize}
\end{proposition}}
\begin{proof} 

[Proof of {\rm (i)}] Let $m=2r$ be even. As $L \in \mathcal{B}[{\bx}]_m$, we first observe that, for any $\phi,q \in \mathbb{R}[{\bx}]_r$, the following 
polynomial Cauchy--Schwarz inequality
\[
L(\phi q)^2\le L(\phi^2)L(q^2).
\tag{CS}
\]
Indeed, for every $t\in\mathbb{R}$,
$0\le L\big((\phi+tq)^2\big)
=
L(\phi^2)+2tL(\phi q)+t^2L(q^2)$, 
and hence the discriminant of the quadratic polynomial on the right-hand
side is nonpositive. This shows that 
(CS) holds. 

Let
$\ba\in\mathbb{N}_m^n$ and choose
$\bb,\bgamma\in\mathbb{N}^n$ such that
\[
\ba=\bb+\bgamma,
\qquad
|\bb|\le r,
\qquad
|\bgamma|\le r.
\]
Such a decomposition always exists because $|\ba|\le2r$.

By \textup{(CS)},
\[
\left|L(\bx^\ba)\right|^2
=
\left|L(\bx^\bb\bx^\bgamma)\right|^2
\le
L(\bx^{2\bb})L(\bx^{2\bgamma}).
\]
Applying Lemma \ref{lemma:claim} with $s=r$ to both factors gives
\[
\begin{aligned}
\left|L(\bx^\ba)\right|^2
&\le
\left(
\prod_{i=1}^n L(x_i^{2r})^{\beta_i/r}
\right)
\left(
\prod_{i=1}^n L(x_i^{2r})^{\gamma_i/r}
\right)=
\prod_{i=1}^n L(x_i^{2r})^{\alpha_i/r}.
\end{aligned}
\]
Taking square roots and using $m=2r$, we obtain
\[
\left|L(\bx^\ba)\right|
\le
\prod_{i=1}^n L(x_i^m)^{\alpha_i/m}
=
\left(
\Phi_{\mJ_{\nm}^{\mathrm{enpc}}}(L)
\right)^\ba.
\]
Clearly, if $\ba={\bf 0}$, then equality holds in the above inequality. In addition, if $\ba=m\be_i$, then equality also holds because
\[
\left(
\Phi_{\mJ_{\nm}^{\mathrm{enpc}}}(L)
\right)^{m\be_i}
=
L(x_i^m)=L(\bx^{m\be_i}).
\]

[Proof of {\rm (ii)}] Finally, let $f  \in\mJ_{\nm}^{\mathrm{enpc}}$. Then, 
$f(\bx)
=
\sum_{i=1}^n f_{m,i}x_i^m
+
\sum_{\ba\in\Omega\backslash \{{\bf 0}\}}
%{\substack{\ba\in\mathbb{N}_m^n \backslash \{{\bf 0}\}\\
%\ba\notin\{m\be_1,\ldots,m\be_n\}}}
f_\ba\bx^\ba +f_{\bf 0}$, 
where $f_\ba\le0$ for all $\ba\in\Omega\backslash \{{\bf 0}\}$. 
It then follows from part {\rm (i)}
%, for all $\alpha\in\mathbb{N}_m^n$ 
%\[
%L(\bx^\ba)
%\le
%\left(
%\Phi_{\mJ_{\nm}^{\mathrm{enpc}}}(L)
%\right)^\ba.
%\]
%Therefore, 
that for all $\ba\in\Omega\backslash \{{\bf 0}\}$,
%$\ba\in\mathbb{N}_m^n \backslash \{{\bf 0}\}$ with $\alpha \neq m \be_i$, $i=1,\ldots,n$, 
\[
L(f_\ba  \, \bx^\ba) = f_\ba \, L(\bx^\ba)
\ge
f_\ba
\left(
\Phi_{\mJ_{\nm}^{\mathrm{enpc}}}(L)
\right)^\ba.
\] 
Moreover, 
%for the monomials $x_i^m$ or the constant term, equality holds in the above expression, that is, 
for $\ba=m \be_i$, $i\in[n]$, or $\ba=\bf 0$, 
\[
L(f_\ba  \, \bx^\ba) = 
f_\ba
\left(
\Phi_{\mJ_{\nm}^{\mathrm{enpc}}}(L)
\right)^\ba.
\]
Summing over all monomials gives
$L(f)
\ge
f\left(
\Phi_{\mJ_{\nm}^{\mathrm{enpc}}}(L)
\right)$. Thus, the conclusion follows. 
%So, $\mJ_{\nm}^{\mathrm{enpc}}$ is a Jensen class with Jensen map
%$\Phi_{\mJ_{\nm}^{\mathrm{enpc}}}$.
\end{proof}

\subsection{Componentwise univariate compositions}\label{sec::comp} 
Next, using the composition principles developed in the previous section, we identify several classes of nonconvex Jensen-type polynomials with a composition structure.  %{\color{blue} Make the notation consistent {\rm sos-c} or {\rm c-sosc}}
\begin{proposition}[Subclass 2: Composition of SOS-convexity and componentwise  univariate polynomials]\label{prop::subclass2}
Let $h=(h_1,\ldots,h_n)$ where each $h_i$ is a univariate polynomial. Let $m_g, m_h$ be positive integers such that {$m_h =\max_{1\le i\le n}\deg(h_i)$}, and set $m:=m_gm_h$. Define
\[
\mJ^{\csosc}_{h,m}:=\left\{g(h_1(x_1),\ldots,h_n(x_n))\in\RR[\bx]_m\mid
g \in \RR[\bxi]_{m_g} \text{ is SOS-convex on } \RR^n\right\}. 
\]
Define the map $\Phi_{\mJ^{\csosc}_{h,m}}: \mathcal{B}[\bx]_m \to \RR^n$ 
in the following way: for any $L\in\mathcal{B}[\bx]_m$, 
%then let $\varphi_3(L)$ be an arbitrarily fixed point in $\RR^n$; otherwise, 
let $\Phi_{\mJ^{\csosc}_{h,m}}(L)=\bu$ where $\bu=(u_1,\ldots,u_n)\in\RR^n$ is chosen so that it satisfies 
$h_i(u_i)=L(h_i(x_i))$ for all $i\in[1:n]$. Then,  \begin{itemize}
\item[{\rm (i)}] the set $\mJ^{\csosc}_{h,m}$ is a Jensen class with the associated Jensen map $\Phi_{\mJ^{\csosc}_{h,m}}$.
\item[{\rm (ii)}] For any $ f \in \mJ^{\csosc}_{h,m}$, if $f^{\star}:=\inf_{\bx\in\RR^n}f(\bx)>-\infty$,
    then $f(\bx)-f^{\star}\in\Sigma[\bx]_m$. In particular,  for any $f \in \mJ^{\csosc}_{h,m}$, $f$ is nonnegative if and only if $f \in \Sigma[\bx]_m$. 
\end{itemize}
\end{proposition}
\begin{proof}
Firstly, we see that 
\[
h(\mathbb{R}^n)=\{(h_1(x_1),\ldots,h_n(x_n)) \mid (x_1,\ldots,x_n) \in \RR^n\}= \prod_{i=1}^n h_i(\mathbb{R}).
\]
Note that the class of all SOS-convex polynomials forms a Jensen class associated with the natural moment map as the Jensen map. By Lemmas \ref{lem::compatible_conditions} and \ref{lemma:3.3},
to see the conclusion, it suffices to show that, for each $i=1,\ldots,n$, $h_i(\mathbb{R})$ is closed and convex. If $h_i$ takes constant value, then trivially, $h_i(\mathbb{R})$ is closed and convex. Thus, we assume that each $h_i$ is not a constant polynomial. We now divide the discussion in three cases. 

Case 1. $\deg(h_i)$ is odd. In this case,  $h_i(\RR)=\RR$ and so, the conclusion follows trivially; 

Case 2. $\deg(h_i)$ is even and its leading coefficient is positive. Then, $h_i(\RR)=[l_i, +\infty)$ for some $l_i\in\RR$ which is closed and convex. 
%Moreover, as $\psi_i(t)-l_i$ is nonnegative for all $t \in \R$, this polynomial is SOS. So,  
%    \[
%    L(\psi_i(x_i)-l_i)=L(\psi_i(x_i))-l_i\ge 0,
%    \]
%    which implies that $L(\psi_i(x_i))\in\psi_i(\RR)$.

    Case 3. $\deg(h_i)$ is even and its leading coefficient is negative. Then, $h_i(\RR)=(-\infty,u_i]$ for some $u_i\in\RR$ which is also closed and convex. 
    %Again, a similar argument applying to $t \mapsto u_i-\psi_i(t)$ gives the desired conclusion $L(\psi_i(x_i))\in\psi_i(\RR)$. 

    Thus, in all three cases, $h_i(\RR)$ is closed and convex. Hence
$h(\RR^n)$ is closed and convex, and the result follows from
Lemmas~\ref{lem::compatible_conditions} and~\ref{lemma:3.3}.
\end{proof}

\begin{remark}[Alternative approach for the SOS representation ] The SOS representation in part {\rm (ii)} of Proposition \ref{prop::subclass2} can be derived independently. Indeed, assume without loss of generality that each $h_i$ is not a constant polynomial, and let 
\[
\begin{aligned}
   I_1&=\{i: {\rm deg}(h_i) \mbox{ is odd}\},\\
   I_2&=\{i: {\rm deg}(h_i) \mbox{ is even and leading coefficient is positive}\},\\
   I_3&=\{i: {\rm deg}(h_i) \mbox{ is even and leading coefficient is negative.}\}. 
\end{aligned}  
\]
It follows that
\[
h(\RR^n)=\prod_ {i=1}^n h_i(\RR)= \prod_ {i \in I_1} \RR \times \prod_ {i \in I_2} [\ell_i,+\infty) \times \prod_ {i \in I_3} (-\infty,u_i], 
\]
where $\ell_i=\inf\{h_i(x_i): x_i \in \mathbb{R}\}$, $i \in I_2$ and $u_i=\sup\{h_i(x_i): x_i \in \mathbb{R}\}$, $i \in I_3$. Direct verification shows that the interior of $h(\RR^n)$ is nonempty. 
Recall that $f=g\circ h$ where $g$ is SOS-convex polynomial and $h(\bx)=(h_1(x_1),\ldots,h_n(x_n))$ with each $h_i$ being a univariate polynomial. So, if $f \ge 0$,  then 
\[
g(\xi_1,\ldots,\xi_n) \ge 0 \mbox{ on } h(\RR^n)=\prod_ {i \in I_1} \RR \times \prod_ {i \in I_2} [\ell_i,+\infty) \times \prod_ {i \in I_3} (-\infty,u_i],
\]
Note that the interior of $h(\mathbb{R}^n)$ is nonempty. So, by Lagrangian duality, there exist $\lambda_i \ge 0$, $i \in I_2 \cup I_3$ such that
\[
\sigma(\xi_1,\ldots,\xi_n)=g(\xi_1,\ldots,\xi_n) + \sum_{i \in I_2} \lambda_i (\ell_i-\xi_i) 
+ \sum_{i \in I_3} \lambda_i (\xi_i-u_i) \ge 0 \mbox{ for all }  (\xi_1,\ldots,\xi_n) \in \RR^n. 
\]
Since \(\sigma\) is SOS-convex and nonnegative, it is an SOS polynomial. Thus, by letting $\xi_i= h_i(x_i)$, $i=1,\ldots,n$, we have the following representation
\[
f(\bx)=\sigma(h_1(x_1),\ldots,h_n(x_n)) + \sum_{i \in I_2} \lambda_i (h_i(x_i)-\ell_i) + \sum_{i \in I_3} \lambda_i (u_i-h_i(x_i))
\]
Note that $\ell_i=\inf h_i$, $i \in I_2$ and $u_i=\sup h_i$, $i \in I_3$, and each $h_i$ is a univariate polynomial, and so, 
\[
h_i(x_i)-\ell_i =\bar \sigma_i(x_i), \ i \in I_2, \mbox{ and } u_i-h_i(x_i) =\bar \sigma_i(x_i), \ i \in I_3,
\]
for some SOS polynomials $\bar \sigma_i$. Thus, we see that $f$ is SOS.
\end{remark}

Below, we give a simple illustrative non-convex example which shows that this composition class goes beyond SOS-convexity, and is distinct
from the class of polynomials with ENPCs. 
\begin{example}\label{example:4.1}
Consider 
$f(x_1,x_2) = x_1^8 + x_2^8 - 2x_1^4 x_2^4 + x_1^4$. 
Note that $f=g\circ h$ where $g(\xi_1, \xi_2) = (\xi_1 - \xi_2)^2 + \xi_1 + 1$ which is SOS-convex, and 
\[
h(x_1,x_2)=(h_1(x_1),h_2(x_2)) 
\]
with $h_1(t)=h_2(t)=t^4-1$. So, $f \in \mJ^{\csosc}_{h,m}$ with $m=8$.

On the other hand, $f$ is not a polynomial with essentially non-positive coefficients. Moreover, a direct computation gives
\[
\nabla^2 f(1,\tfrac12) =
\begin{pmatrix}
66.5 & -4 \\
-4 & -5.125
\end{pmatrix},
\]
whose determinant is \(66.5 \cdot (-5.125) - 16 < 0\). So \(f\) is not convex and therefore not SOS-convex.  \qed
\end{example} 

{\begin{proposition}[Subclass 3: Composition of polynomials with ENPCs
and componentwise univariate polynomials]
Let $h=(h_1,\ldots,h_n)$, where each $h_i$ is a univariate polynomial
satisfying either
\begin{itemize}
\item[{\rm (a)}] $\deg(h_i)$ is odd; or
\item[{\rm (b)}] $\deg(h_i)$ is even and its leading coefficient is positive.
\end{itemize}
Let $m_h= \max_{1\le i\le n}\deg(h_i)$.   
Let $m_g$ be an even positive integer, and set $m:=m_gm_h$. Define
\[
\mJ^{\cenpc}_{h,m}
:=
\left\{
g\bigl(h_1(x_1),\ldots,h_n(x_n)\bigr)\in\mathbb{R}[\bx]_m
\ \middle|\
g(\xi_1,\ldots,\xi_n)\in
\mJ_{\subnm{n}{m_g}}^{\enpc}
\subset\mathbb{R}[\bxi]_{m_g}
\right\}.
\]
Define the map
$\Phi_{\mJ^{\cenpc}_{h,m}}:
\mathcal{B}[\bx]_m\to\mathbb{R}^n$ 
as follows: for any $L\in\mathcal{B}[\bx]_m$, let
$\Phi_{\mJ^{\cenpc}_{h,m}}(L)=\bu$, where
$\bu=(u_1,\ldots,u_n)\in\mathbb{R}^n$ is chosen so that it satisfies
\[
h_i(u_i)=\sqrt[m_g]{L\bigl(h_i^{m_g}(x_i)\bigr)}
\qquad\text{for each }i\in[1:n].
\]
Then, the following statements hold. 
\begin{itemize}
\item[{\rm (i)}] The set $\mJ^{\cenpc}_{h,m}$ is a Jensen class with
associated Jensen map $\Phi_{\mJ^{\cenpc}_{h,m}}$.

\item[{\rm (ii)}] For any $f\in\mJ^{\cenpc}_{h,m}$, if
$f^\star:=\inf_{\bx\in\mathbb{R}^n}f(\bx)>-\infty$, 
then $f-f^\star\in\Sigma[\bx]_m$. In particular, for any
$f\in\mJ^{\cenpc}_{h,m}$, $f$ is nonnegative if and only if
$f\in\Sigma[\bx]_m$.
\end{itemize}
\end{proposition}

\begin{proof}
Let $h_\sharp L$ be defined as in \eqref{eq:push-forward}, with
$m_1=m_g$ and $m_2=m_h$. By Lemma~\ref{lem::compatible_conditions},
it suffices to show that
\[
\Phi_{\mJ^{\enpc}_{n,m_g}}(h_\sharp L)\in h(\mathbb{R}^n)
\qquad
\text{for every }L\in\mathcal{B}[\bx]_m.
\]
By the definitions of $\Phi_{\mJ^{\enpc}_{n,m_g}}$ and $h_\sharp L$,
\[
\begin{aligned}
\Phi_{\mJ^{\enpc}_{n,m_g}}(h_\sharp L)
&=
\left(
\sqrt[m_g]{(h_\sharp L)(\xi_1^{m_g})},\ldots,
\sqrt[m_g]{(h_\sharp L)(\xi_n^{m_g})}
\right)\\
&=
\left(
\sqrt[m_g]{L(\xi_1^{m_g}\circ h)},\ldots,
\sqrt[m_g]{L(\xi_n^{m_g}\circ h)}
\right)\\
&=
\left(
\sqrt[m_g]{L(h_1^{m_g}(x_1))},\ldots,
\sqrt[m_g]{L(h_n^{m_g}(x_n))}
\right).
\end{aligned}
\]
For each $i=1,\ldots,n$, we consider two cases.

Case (a). If $\deg(h_i)$ is odd, then $h_i(\mathbb{R})=\mathbb{R}$, so trivially 
$\sqrt[m_g]{L(h_i^{m_g}(x_i))}\in h_i(\mathbb{R})$.

Case (b). If $\deg(h_i)$ is even and its leading coefficient is positive, then
$h_i(\mathbb{R})=[\ell_i,+\infty)$ for some $\ell_i\in\mathbb{R}$.
Consequently,
$h_i^{m_g}(x_i)\geq \bigl(\max\{0,\ell_i\}\bigr)^{m_g}$ 
for all $x_i\in\mathbb{R}$. 
Thus
\[
h_i^{m_g}(x_i)-\bigl(\max\{0,\ell_i\}\bigr)^{m_g}
\in\Sigma[\bx]_m.
\]
Since $L\in\mathcal{B}[\bx]_m$, it follows that
$L(h_i^{m_g}(x_i))
\geq \bigl(\max\{0,\ell_i\}\bigr)^{m_g}$. 
Therefore,
$\sqrt[m_g]{L(h_i^{m_g}(x_i))}
\geq\max\{0,\ell_i\}$, 
and hence $\sqrt[m_g]{L(h_i^{m_g}(x_i))}\in h_i(\mathbb{R})$.

In either case,
$\Phi_{\mJ^{\enpc}_{n,m_g}}(h_\sharp L)\in h(\mathbb{R}^n)$
for every $L\in\mathcal{B}[\bx]_m$. 
The result now follows from Lemma~\ref{lem::compatible_conditions}.
\end{proof}}

{We note that, if the assumptions {\rm (a)} and {\rm (b)} are not satisfied, then %$\Phi_{\mJ^{\cenpc}_{h,m}}$ need not to be well-defined, and so, 
the conclusion of the previous proposition may fail. Indeed, 
let $n=2$, $m_h=2$, $m_g=4$, $m=m_gm_h =8$, and define
\[
h:\mathbb{R}^2\to\mathbb{R}^2,
\qquad
h(x_1,x_2):=(-x_1^2,-x_2^2).
\]
Thus, both components of $h$ are quadratic univariate functions with negative
leading coefficients. Define
\[
g(\xi_1,\xi_2)
:=
-\xi_1^2\xi_2-\xi_1\xi_2^2-3\xi_1\xi_2+1.
\]
Then $g\in\mJ_{\subnm{2}{4}}^{\enpc}$ and
%$g\in\mathbb{R}[\bxi]_4$, and all its nonconstant coefficients
%are negative. The 
the composition $f:=g\circ h$ is
\[
\begin{aligned}
f(x_1,x_2)
&=
g(-x_1^2,-x_2^2)=
x_1^4x_2^2+x_1^2x_2^4-3x_1^2x_2^2+1, 
\end{aligned}
\]
which is the Motzkin polynomial. It is well known that Motzkin polynomial is a nonnegative polynomial but is not SOS (see, for example \cite{Lasserre09,Nie2023}). 
%By the arithmetic--geometric mean inequality,
%\[
%x_1^4x_2^2+x_1^2x_2^4+1
%\ge
%3\sqrt[3]{x_1^6x_2^6}
%=
%3x_1^2x_2^2.
%\]
%Consequently,
%\[
%f(x_1,x_2)\ge 0
%\qquad
%\text{for all }%(x_1,x_2)\in\mathbb{R}^2.
%\]
%Nevertheless, the Motzkin %polynomial is not a sum of squares of
%polynomials; that is,
%\[
%f\notin\Sigma[x_1,x_2].
%\]
Thus,
$f=g\circ h\ge0$ 
 on $\mathbb{R}^2$,
but
%\[
$f\notin\Sigma[x_1,x_2]$.
%\]
This shows that the SOS conclusion in part~{\rm (ii)} may fail if the
assumptions on the degrees and leading coefficients of the components of $h$ are removed.}

\subsection{Structured quadratic-map compositions}
Next, we consider the function which can be written as composition of SOS-convex polynomial and quadratic maps. An illustrative example is the loss function in source localization problem, which takes the form 
\[
f(\bx)= \sum_{\ell=1}^p (\|\bx-{\bv}_\ell\|_U^2 -r_\ell^2)^{2d},
\]
where ${\bv}_\ell \in \RR^n$ is the given location of the sensor $\ell$, $r_\ell > 0$ is the measured distance to each sensor $\ell$ (often corrupted by noise), $p \in \mathbb{N}$ and $\|\cdot\|_U$ is the norm induced by a positive definite symmetric matrix $U$ and is given by $\|\bv\|_U=\sqrt{\bv^TU\bv}$ for all $\bv$. In this case, $f=g\circ h$ where 
$g(\bxi)=\sum_{\ell=1}^p \xi_\ell^{2d}$ and $h(\bx)=(h_1(\bx),\ldots,h_p(\bx))$ 
where $h_\ell(\bx)=\|\bx-{\bv}_\ell\|_U^2 -r_\ell^2,$ $\ell=1,\ldots,p$. 
\begin{proposition}[Subclass 4: Composition of SOS-convexity and structured quadratic maps]\label{prop::subsclass4}
Let $h=(h_1,\ldots,h_p)$ where  $h_\ell$, $\ell=1,\ldots,p$, are quadratic functions. Let $m_g$ be a positive integer, $m=2m_g$, and define
\[
\mJ^{\csosc}_{h,m}:=\{g(h_1(\bx),\ldots,h_p(\bx))\in\RR[\bx]_m\mid
g(\xi_1,\ldots,\xi_p) \in \mathbb{R}[\bxi]_{m_g} \text{ is SOS-convex on } \RR^p\},
\]    Define the map $\Phi_{\mJ^{\csosc}_{h,m}}: \mathcal{B}[\bx]_m \to \RR^n$ 
as, for any $L\in\mathcal{B}[\bx]_m$, 
%then let $\varphi_3(L)$ be an arbitrarily fixed point in $\RR^n$; otherwise, 
 $\Phi_{\mJ^{\csosc}_{h,m}}(L)=\bu$ where $\bu \in\RR^n$ satisfies 
$h_\ell(\bu)=L(h_\ell(\bx))$ for all $\ell \in[1:p]$. 
Suppose that one of the following conditions hold:
\begin{itemize}
\item[{\rm (a)}] $p=3$, {$n \ge 3$}, $h_\ell(\bx)=\bx^TA_\ell\bx$ with $A_\ell \in \mathbb{S}^{n}$ and there exist $\lambda_\ell \in \mathbb{R}$, $\ell=1,2,3$ such that $\displaystyle \sum_{\ell=1}^3 \lambda_\ell A_\ell \succ 0$. 
\item[{\rm (b)}] $h_\ell(\bx)=\bx^TA_\ell\bx+\alpha_\ell$ with $A_\ell \in \mathbb{S}^{n}$,  and $\{A_1,\ldots,A_p\}$ are simultaneously diagonalizable.  
\item[{\rm (c)}] $h_\ell(\bx)=\bx^TA\bx+\bsa_\ell^T\bx+\alpha_\ell$, with $A \in \mathbb{S}^n$ is positive definite, $\bsa_\ell \in \mathbb{R}^n$ and $\alpha_\ell \in \mathbb{R}$,  and ${\rm dim} \, {\rm span}\{\bsa_2-\bsa_1,\ldots,\bsa_p-\bsa_1\} \le n-1$. \end{itemize}
Then, the following statements hold.  \begin{itemize}
\item[{\rm (i)}] The set $\mJ^{\csosc}_{h,m}$ is a Jensen class with the associated Jensen map $\Phi_{\mJ^{\csosc}_{h,m}}$.
\item[{\rm (ii)}] For any $f \in \mJ^{\csosc}_{h,m}$,  if $f^{\star}:=\inf_{\bx\in\RR^n}f(\bx)>-\infty$,
    then $f(\bx)-f^{\star}\in\Sigma[\bx]_m$.  In particular, for any $f \in \mJ^{\csosc}_{h,m}$, $f$ is nonnegative if and only if $f \in \Sigma[\bx]_m$. 
\end{itemize}

\end{proposition}
\begin{proof}
From Lemma \ref{lemma:3.3}, it suffices to show that $h(\RR^n)$ is closed and convex.

Case {\rm (a)}: It follows from Polyak's result \cite[Theorem 2.1]{polyak1998convexity}.

Case {\rm (b)}: As $\{A_1,\ldots,A_p\}$ is simultaneously diagonalizable, 
there exists an orthonormal matrix $W \in \RR^{n \times n}$ such that $\hat{A}_\ell:=W^TA_\ell W$ are diagonal matrices. Thus, 
\[
h(\RR^n)=\{(\bx^T\hat{A}_1\bx,\ldots,\bx^T\hat{A}_p\bx) \mid \bx \in \RR^n\}+(\alpha_1,\ldots,\alpha_p).
\]
Write
$\widehat A_\ell
=
\operatorname{Diag}(a_{\ell 1},\ldots,a_{\ell n})$, $\ell=1,\ldots,p$.  
Setting
$t_i=x_i^2\ge0$, we obtain
\[
\begin{aligned}
h(\mathbb{R}^n)
=
\left\{
\left(
\sum_{i=1}^n a_{1i}t_i,\ldots,
\sum_{i=1}^n a_{pi}t_i
\right)\ \middle|\
(t_1,\ldots,t_n)\in\mathbb{R}^n_+
\right\}+
(\alpha_1,\ldots,\alpha_p), 
\end{aligned}
\]
% Thus, $h(\mathbb{R}^n)$ is a finitely generated polyhedral cone
% translated by $(\alpha_1,\ldots,\alpha_p)$. In particular, it is closed
% and convex.
which is a convex polyhedral set, and so, is closed and convex.

Case {\rm (c)}: Let us show that $K=h(\RR^n)$ is closed and convex. Note that convexity is preserved under invertible linear maps, to see $K$ is convex, it suffices to show that 
\[
K'=\{(h_1(\bx),h_2(\bx)-h_1(\bx),\ldots, h_p(\bx)-h_1(\bx)): \bx \in \RR^n\}
\]
is convex. Let 
$\bxi=(\xi_1,\ldots,\xi_p) \in K'$ and $\bz=(z_1,\ldots,z_p) \in K'$. 
So, there exist $\bu$ and $\bv$ such that 
\[
\xi_1= \bu^TA\bu+ \bsa_1^T\bu+\alpha_1,  \ \xi_\ell=  (\bsa_\ell-\bsa_1)^T\bu+(\alpha_\ell-\alpha_1), \ \ell=2,\ldots,p\]
and \[ z_1= \bv^TA\bv+ \bsa_1^T\bv+\alpha_1, \ z_\ell= (\bsa_\ell-\bsa_1)^T\bv+(\alpha_\ell-\alpha_1), \ \ell=2,\ldots,p.
\]

Fixing a $\lambda\in (0, 1)$, we next show that $\lambda \bxi+(1-\lambda) \bz \in K'$. 
Let {$\bw=\lambda \bu+(1-\lambda) \bv$}. Noting that $A$ is positive definite,
by convexity, 
\[
\lambda \xi_1+(1-\lambda)z_1 \ge \bw^TA\bw+ \bsa_1^T\bw+\alpha_1 
\ \ \mbox{ and }\ \ \lambda \xi_\ell+(1-\lambda)z_\ell = (\bsa_\ell-\bsa_1)^T\bw+(\alpha_\ell-\alpha_1), \ \ell=2,\ldots,p. 
\]  
From the assumption, there exists $\bd \in \RR^n \backslash\{\bf 0\}$ such that 
$(\bsa_i-\bsa_1)^T\bd=0$. Thus, consider $\bw(t)=\bw + t \bd$ with $t \ge 0$. 
We see that, for all $t \ge 0$, 
\[
\lambda \xi_\ell+(1-\lambda)z_\ell = (\bsa_\ell-\bsa_1)^T\bw(t)+(\alpha_\ell-\alpha_1), \ \ell=2,\ldots,p
\]
and 
\[
\bw(t)^TA\bw(t)+ \bsa_1^T\bw(t)+\alpha_1 \rightarrow +\infty \ \ \mbox{ as }\ \   t \rightarrow +\infty.
\]
Thus, by the intermediate value theorem, there exists $t_0 \ge 0$ such that
\[
\lambda \xi_1+(1-\lambda)z_1= \bw(t_0)^TA\bw(t_0)+ \bsa_1^T\bw(t_0)+\alpha_1.
\]
Thus, $\lambda \bxi+(1-\lambda) \bz \in K'$, and so, $K'$ is convex. Thus, $K$ is also convex. 

To see $K$ is closed, let $(\xi_1^{(k)},\ldots,\xi_p^{(k)}) \in K$ with 
$(\xi_1^{(k)},\ldots,\xi_p^{(k)}) \rightarrow (\xi_1,\ldots,\xi_p)$. 
Then, there exist $\bx^{(k)} \in \RR^n$ such that  
\[
\xi_\ell^{(k)}=(\bx^{(k)})^TA(\bx^{(k)})+ \bsa_\ell^T(\bx^{(k)})+\alpha_\ell \rightarrow \xi_\ell.
\]
We now see that $\{\bx^{(k)}\}$ is bounded. (Otherwise, by passing to a subsequence if necessary, we assume $\|\bx^{(k)}\| \rightarrow \infty$.
By further passing to a subsequence, we can further assume that $\frac{\bx^{(k)}}{\|\bx^{(k)}\|} \rightarrow \hat{\bx} \neq 0$ for some $\hat{\bx}\in\RR^n$. 
Then, by dividing $\|\bx^{(k)}\|^2$ and passing to limit,  we obtain that $\hat{\bx}^TA\hat{\bx}=0$ which contradicts
to the fact that $A$ is positive definite.) 
Thus, $\{\bx^{(k)}\}$ is bounded. 
By passing to a subsequence, we may assume that $\bx^{(k)} \rightarrow \bar{\bx}$, 
and so, by passing to the limit, 
\[
\xi_\ell=\bar{\bx}^TA\bar{\bx}+ \bsa_\ell^T\bar{\bx}+\alpha_\ell, \ \ell=1,\ldots,p.
\]
Thus, $(\xi_1,\ldots,\xi_p)\in K$, and hence $K$ is closed. 
\end{proof} 

In passing, we note that the dimension condition in ${\rm (c)}$ can be efficiently checked. Also, the Jensen map indeed has a closed-form in this case. We will illustrate this in the applications section (Section 7) with the source localization problem. 
{ \subsection{Combinations of Jensen-type polynomials}

Let $m$ and $m'$ be positive even integers with $m\leq m'$.
For any $L\in\mathcal{B}[\bx]_{m'}$, its restriction
$L|_{\mathbb{R}[\bx]_m}$ belongs to $\mathcal{B}[\bx]_m$.

\begin{proposition}\label{prop::combination}
Let $\mJ_1\subseteq\mathbb{R}[\bx]_m$ and
$\mJ_2\subseteq\mathbb{R}[\bx]_{m'}$ be Jensen classes with associated
Jensen maps $\Phi_1$ and $\Phi_2$, respectively. Suppose that
\[
\Phi_1\bigl(L|_{\mathbb{R}[\bx]_m}\bigr)=\Phi_2(L)
\qquad
\text{for every }L\in\mathcal{B}[\bx]_{m'}.
\]
Then $\mJ_1+\mJ_2$ is a Jensen class in
$\mathbb{R}[\bx]_{m'}$, with associated Jensen map $\Phi_2$.
\end{proposition}

\begin{proof}
Take $f\in\mJ_1$, $g\in\mJ_2$, and
$L\in\mathcal{B}[\bx]_{m'}$. Since
$L|_{\mathbb{R}[\bx]_m}\in\mathcal{B}[\bx]_m$, the Jensen inequalities
for $\mJ_1$ and $\mJ_2$ give us that 
\[
\begin{aligned}
L(f+g)
=L|_{\mathbb{R}[\bx]_m}(f)+L(g) 
\geq f\bigl(\Phi_1(L|_{\mathbb{R}[\bx]_m})\bigr)
      +g\bigl(\Phi_2(L)\bigr)
&=f\bigl(\Phi_2(L)\bigr)+g\bigl(\Phi_2(L)\bigr)\\
&=(f+g)\bigl(\Phi_2(L)\bigr).
\end{aligned}
\]
Hence $\mJ_1+\mJ_2$ is a Jensen class with Jensen map $\Phi_2$.
\end{proof}

We now apply Proposition~\ref{prop::combination} to combine the ENPC
class with a class of SOS-convex compositions. Fix a positive even
integer $m$, and let $m_g$ be a positive even integer. Set
$m':=m \, m_g$ and define
\[
h'(\bx):=(x_1^m,\ldots,x_n^m).
\]
Consider the Jensen class
\[
\mJ^{\csosc}_{h',m'}
:=
\left\{
g(x_1^m,\ldots,x_n^m)\in\mathbb{R}[\bx]_{m'}
\ \middle|\
g\in\mathbb{R}[\bxi]_{m_g}
\text{ is SOS-convex on }\mathbb{R}^n
\right\}.
\]
For $L\in\mathcal{B}[\bx]_{m'}$, its associated Jensen map is
\[
\Phi_{\mJ^{\csosc}_{h',m'}}(L)
=
\left(
\sqrt[m]{L(x_1^m)},\ldots,\sqrt[m]{L(x_n^m)}
\right).
\]
Indeed, if $\bu=\Phi_{\mJ^{\csosc}_{h',m'}}(L)$, then
$u_i^m=h'_i(u_i)=L\bigl(h'_i(\bx)\bigr)=L(x_i^m)$, $i=1,\ldots,n$.  
This is the Jensen map given by Subclass~2. Moreover, it holds that 
$\Phi_{\mJ^{\enpc}_{\subnm{n}{m}}}(L|_{\RR[\bx]_m})=\Phi_{\mJ^{\csosc}_{h',m'}}(L)$
for all $L\in\mathcal{B}[\bx]_{m'}$.
Therefore, Proposition~\ref{prop::combination} implies that
$\mJ^{\enpc}_{\subnm{n}{m}}+\mJ^{\csosc}_{h',m'}
$ 
is a Jensen class with associated Jensen map
$\Phi_{\mJ^{\csosc}_{h',m'}}$.

\begin{corollary}
Let
$\phi(\bx)=f(\bx)+g(x_1^m,\ldots,x_n^m)$, 
where $f\in\mJ^{\enpc}_{\subnm{n}{m}}$ and
$g\in\mathbb{R}[\bxi]_{m_g}$ is SOS-convex. If
$\phi^\star:=\inf_{\bx\in\mathbb{R}^n}\phi(\bx)>-\infty$,
then
\[
\phi-\phi^\star\in\Sigma[\bx]_{m'}.
\]
\end{corollary}}

\begin{example}\label{example:4.2}
Let \(m=2\). Consider the ENPC polynomial
$f(x_1,x_2,x_3) = -x_1x_2 - x_1x_3 - x_2x_3\in\mJ^{\enpc}_{\subnm{3}{2}}$.
Let 
$g(\xi_1,\xi_2,\xi_3) = (\xi_1+\xi_2+\xi_3)^2 = \xi_1^2+\xi_2^2+\xi_3^2 + 2\xi_1\xi_2 + 2\xi_1\xi_3 + 2\xi_2\xi_3$, 
which is SOS-convex. Define the polynomial
\[
\phi(x_1,x_2,x_3) := f(x_1,x_2,x_3) + g(x_1^2,x_2^2,x_3^2)=(x_1^2+x_2^2+x_3^2)^2 - (x_1x_2+x_1x_3+x_2x_3).
\]
%The following properties hold for $h$:
%\begin{enumerate}[\upshape (i)]
%\item \(h\) does not have ENPCs. 
%\item \(h\) is \emph{not} a simple SOS-convex substitution \(g(\psi_1(x_1),\psi_2(x_2),\psi_3(x_3))\). If it were, the presence of the term \(-x_1x_2\) would force the linear parts of \(\psi_1,\psi_2\) to be nonzero, but then the Hessian of \(g\) would be indefinite at the origin (since \(g\) would have a cross term but no compensating quadratic terms in \(y_1^2,y_2^2,y_3^2\)), contradicting SOS-convexity.
%\item 
Since 
\[
\begin{aligned}
    \phi(\bx)&=(x_1^2+x_2^2+x_3^2)^2 - (x_1x_2+x_1x_3+x_2x_3)\\
    &\ge (x_1^2+x_2^2+x_3^2)^2-(x_1^2+x_2^2+x_3^2)\\
    &=\left(x_1^2+x_2^2+x_3^2-\frac{1}{2}\right)^2-\frac{1}{4}\ge -\frac{1}{4}.
\end{aligned}
\]
and the equality holds when \(x_1=x_2=x_3=\pm 1/\sqrt{6}\), it implies that
$\phi^{\star} := \inf_{\mathbb{R}^3} \phi = -\frac14.$ We have the SOS expression 
%\item The shifted polynomial is SOS:
\[
\phi(\bx) - \phi^{\star} = \phi(\bx) + \frac14 = \left(x_1^2+x_2^2+x_3^2 - \frac12\right)^2
+ \frac12\big((x_1-x_2)^2+(x_2-x_3)^2+(x_3-x_1)^2\big)\in\Sigma[\bx]_4.
\]
%\end{enumerate}
\end{example}

{

\section{Theorems of the alternative for Jensen-type polynomials with universal quantifiers}\label{sec::alt}

In this section, we establish a theorem of the alternative for
systems of Jensen-type polynomial inequalities involving universal
quantifiers. Let $\Y\subseteq\RR^\ny$ with $\ny \in \mathbb{N}$, 
$f_j\in\mathbb{R}[\bx]$ for $j\in[0:J]$, and
$f\in\mathbb{R}[\bx,\by]$.  Throughout the remainder of the paper, we impose the following assumption. \begin{assumption}\label{assump::0}
(i) $m\in\N$ is even and $\mJ\subset\RR[\bx]_m$ is a Jensen class of polynomials with Jensen map  $\Phi_{\mJ}$;\ \
(ii) $f_j \in \mJ$, for all $j \in [0:J]$,
$f(\cdot,\by)\in\mJ$ for all  $\by\in\Y$ and $\Y$ is compact. 
\end{assumption}

Recall that $\mathcal{M}_+(\Y)$ denotes the space of nonnegative
Borel measures supported on $\Y$. Let
\[
\Theta_m(\bx)=\sum_{i=1}^n x_i^m+1.
\]
The main result in this section is the following theorem of the alternative.

\begin{theorem}[Theorem of the alternative for Jensen-type polynomials with universal quantifiers]\label{th::alt}
Under Assumption \ref{assump::0}, exactly one of the following two statements holds:
    \begin{enumerate}
        \item[\upshape (i)] 
        $[f_0(\bx)<0, \ldots, f_J(\bx)<0, f(\bx,\by)<0 \, \forall \, \by \in \Y ]$ has a solution; %$\bu\in\RR^n$;  %such that $f_j(\bu)<0$ and 
        %$f(\bu,\bv)<0$ for all $j\in[1:J]$ and $\bv\in\Y$;
        \item[\upshape (ii)] there exist $\lambda_j\ge 0$, $j\in[0:J]$, and $\mu\in\mathcal{M}_+(\Y)$ with $\sum_{j=0}^{J}\lambda_j+\mu(\Y)=1$ such that 
        \[
        \sum_{j=0}^J\lambda_j f_j(\bx)+\int_{\Y} f(\bx,\by)\ud\mu(\by)\in\Sigma[\bx]_m.
        \] 
    \end{enumerate}    
Let $\nu$ be an arbitrary, fixed, finite nonnegative Borel measure with $\supp(\nu)=\Y$.
If $\Theta_m(\bx)\in\mJ$, then
%the statement 
${\rm (ii)}$ is further equivalent to 
\begin{enumerate}
        \item[\upshape (iii)] for any $\varepsilon>0$, there exist $\lambda_{\varepsilon,j}\ge 0$, $j=0,\ldots,J$, and
        $\sigma_{\varepsilon}\in\Sigma[\by]$
        with $\sum_{j=0}^J\lambda_{\varepsilon,j}+\int_{\Y}\sigma_{\varepsilon}(\by)\ud\nu(\by)=1$
        such that 
        \[
        \sum_{j=0}^J\lambda_{\varepsilon,j} f_j(\bx)+\int_{\Y} \sigma_{\varepsilon}(\by)f(\bx, \by)\ud\nu(\by)+\varepsilon\Theta_m(\bx)\in\Sigma[\bx]_m.
        \]
    \end{enumerate}
\end{theorem}

Before we proceed to the proof of this theorem, we remark that this extends the main result in \cite{HLQ2016}, the theorem of alternative for finitely many homogeneous ENPC polynomials, to Jensen-type polynomials with universal quantifiers. To see this, we first show that the following corollary holds. 
\begin{corollary}[Theorem of the alternative for finitely many Jensen-type polynomials] 
Let $m$ be even, and let $\mJ \subseteq \mathbb{R}[\bx]_m$ be a Jensen class of polynomials.  Suppose that $f_j \in \mJ$ for every $j\in[0:J]$. Then, exactly one of the following two statements holds:
    \begin{enumerate}
        \item[\upshape (i)] $[f_0(\bx)<0, \ldots, f_J(\bx)<0]$ has a solution;
        \item[\upshape (ii)] there exist $\lambda_j\ge 0$, $j\in[0:J]$, with $\sum_{j=0}^{J} \lambda_j=1$ such that $\sum_{j=0}^J\lambda_j f_j(\bx) \in\Sigma[\bx]_m$. 
    \end{enumerate}  
\end{corollary}
\begin{proof}
It suffices to prove [Not(i) $\Longrightarrow$ (ii)] as the reverse implication holds always.   
Let $\epsilon>0$ and consider the special case where $f(\bx,\by) \equiv -\epsilon$ in Theorem \ref{th::alt}. Then, Theorem \ref{th::alt} implies that 
there exist $\lambda_j\ge 0$, $j\in[0:J]$, and  $\mu\in\mathcal{M}_+(\Y)$ with $\sum_{j=0}^J \lambda_j+\mu(\Y)=1$ such that 
        \[
        \sum_{j=0}^J\lambda_j f_j(\bx) =\sigma(\bx)+\epsilon \, \mu(\Y).
        \] 
for some $\sigma \in \Sigma[\bx]_m$. Clearly, $\mu(\Y) \neq 1$. Otherwise $\lambda_j=0$ for all $j=0,1,\ldots,J$ which leads to $0=\sum_{j=0}^J\lambda_j f_j(\bx)  =\sigma(\bx)+\epsilon>0$. This shows that $\boldsymbol{\lambda}=(\lambda_0,\lambda_1,\ldots,\lambda_J) \neq {\bf 0}$ and  
        $\sum_{j=0}^J\lambda_j f_j(\bx) \in \Sigma[\bx]_m$.  
Thus, the conclusion follows by dividing $\sum_{j=0}^J\lambda_j \neq 0$.         
\end{proof}
\begin{remark}[Discussion on the link with the literature]
If one specializes each $f_j$ to be a homogeneous ENPC polynomial, then the above corollary collapses to the alternative theorem provided in \cite{HLQ2016}. Thus, the preceding theorem of the alternative extends the known result of alternative theorem from finitely many homogeneous ENPC polynomials to Jensen-type polynomials (which covers nonhomogeneous ENPC polynomials) with universal quantifiers.
\end{remark}

%{\color{blue} Can be extended to finitely many quantifiers--- explain briefly (for simplicity)}

We now turn to the proof of Theorem \ref{th::alt}. To do this, we need two intermediate results and some standard notions. A family of functions $\{p_t(\by)\mid t\in T\}$ on $\Y$
is called {\itshape equicontinuous} if for any $\varepsilon>0$, there 
exists a $\delta>0$ such that for all $t\in T$ and all $\bv^{(1)}$,
$\bv^{(2)}\in\Y$ with $\Vert \bv^{(1)}-\bv^{(2)}\Vert<\delta$, it holds
$\vert p_t(\bv^{(1)})-p_t(\bv^{(2)})\vert<\varepsilon$.
It is called {\itshape pointwise bounded} if for each $\bv\in\Y$, 
it holds $\sup_{t\in T}|p_t(\bv)|<\infty$. 
If $\Y$ is compact, and the family $\{p_t(\by)\mid t\in T\}$ is equicontinuous and pointwise bounded, 
then it can be proved that $\{p_t(\by)\mid t\in T\}$ is uniformly bounded on $\Y$. 
Moreover, we have the following result. Given a family of equicontinuous, pointwise bounded functions \(\{p_t\}_{t\in T}\) and nonnegative measures \(\mu\) and \(\nu\), there exists an SOS polynomial \(\sigma_\varepsilon\) such that $\int_{\Y}p_t(\by)\ud\mu(\by)$ can be approximated by $\int_{\Y}p_t(\by)\sigma_{\varepsilon}(\by)\ud\nu(\by)$ uniformly in $t \in T$. Related approximation results have been established for general noncompact domains under the so-called Carleman condition, but they do not provide the precise form needed here (see, e.g., \cite{HKN2024,Lasserre2015}). For completeness, we provide a proof in the appendix.} %{\color{blue} Can we say sth. e.g. link to Jiawang or Lasserre's previous work? any difference or just for completeness/self-containment purpose?}

\begin{proposition}\label{prop::sigma}
    Let $\Y$ be compact, and let $\nu$ be an arbitrary, fixed, finite nonnegative Borel measure on
    $\Y$ with $\supp(\nu)=\Y$. Given a family of equicontinuous and pointwise bounded functions 
    $\{p_t(\by)\mid t\in T\}$ on $\Y$ and nonnegative measure
    $\mu\in\mathcal{M}_+(\Y)$, for any $\varepsilon>0$, 
    there exists an SOS polynomial $\sigma_{\varepsilon}(\by)\in\Sigma[\by]$
    such that $\int_{\Y}\sigma_{\varepsilon}(\by)\ud\nu(\by)=\mu(\Y)$ and 
    \begin{equation}\label{eq::f}
        \left\vert\int_{\Y}p_t(\by)\ud\mu(\by)-\int_{\Y}p_t(\by)\sigma_{\varepsilon}(\by)\ud\nu(\by)\right\vert<\varepsilon,\ \ \forall\ t\in T.
    \end{equation}
\end{proposition}
\begin{proof}
    See Appendix \ref{sec::proof}.
\end{proof}
%{\color{red}
%\begin{remark}
%    The basic idea of the proof is simple. For a fixed $\mu$, we first construct a continuous nonnegative function $g_\varepsilon(\by)$ that serves as a density with respect to 
%    $\nu$ on $\Y$, such that $\int_{\Y} p_t(\by) g_\varepsilon(\by)\ud\nu(\by)$
%approximates $\int_{\Y} p_t(\by) \, \ud\mu(\by)$
%uniformly in $t\in T$. We then approximate the measure $g_\varepsilon(\by)\ud\nu(\by)$ by measures 
%whose densities are SOS polynomials $\sigma_{\varepsilon}(\by)$. 
%This second step is easy in our setting because $\Y$ is compact and the Stone--Weierstrass theorem 
%applies. This step is also considered in \cite{HKN2024}. However, since the integration domain is assumed to be more general there, the argument is more involved and requires the Carleman condition.
%\end{remark}}

For any $L\in(\RR[\bx]_m)^*$, the notation $L(f(\bx,\by))$
means that $L$ acts only on the $\bx$-variables, with $\by$
treated as parameters. More precisely, if
$f(\bx,\by)=\sum_{\ba\in\N^n_m} f_{\ba}(\by)\bx^{\ba}$,
then
\[
L(f(\bx,\by))
=\sum_{\ba\in\N^n_m} f_{\ba}(\by)L(\bx^{\ba})
\in\RR[\by].
\]

\begin{lemma}\label{lem::sos}
Suppose that Assumption \ref{assump::0} holds.
    For any $\lambda_j\ge 0$, $j\in[0:J]$, and $\mu\in\mathcal{M}_+(\Y)$, if the polynomial
        \[        \phi(\bx):=\sum_{j=0}^J \lambda_j f_j(\bx)+\int_{\Y} f(\bx,\by)\ud\mu(\by)
        \] 
        is nonnegative on $\RR^n$, then it is an SOS polynomial.
\end{lemma}
\begin{proof}
    Fix an arbitrary linear functional $L\in\mathcal{B}[\bx]_m$  and write $f(\bx,\by)=\sum_{\ba\in\N^n_m}f_{\ba}(\by)\bx^{\ba}$. 
    Recall from Assumption \ref{assump::0}  that $f_j \in \mJ$ and $f(\cdot,\by) \in \mJ$ for all $j=0,\ldots,J$ and $\by \in \Y$. By the Jensen-type inequality in the 
    condition (C), one has $L(f_j) \ge f_j(\Phi_{\mJ}(L))$ and 
    \[
    \begin{aligned}
       L\left(\int_{\Y} f(\bx,\by)\ud\mu(\by)\right)
    =L\left(\sum_{\ba\in\N^n_m}\int_{\Y}f_{\ba}(\by)\ud\mu(\by)\bx^{\ba}\right) =&\sum_{\ba\in\N^n_m}\int_{\Y}f_{\ba}(\by)\ud\mu(\by)L(\bx^{\ba})\\
    =&\int_{\Y}\sum_{\ba\in\N^n_m}f_{\ba}(\by)L(\bx^{\ba})\ud\mu(\by)\\
    =&\int_{\Y} L(f(\bx,\by))\ud\mu(\by)\ge \int_{\Y} f(\Phi_{\mJ}(L),\by)\ud\mu(\by), 
    \end{aligned}
    \]
    where the last inequality follows by the fact that $L(f(\cdot,\by)) \ge f(\Phi_{\mJ}(L),\by)$ for each $\by\in \Y$. 
    Note that 
      $L(\phi)= \sum_{j=0}^J  \lambda_j  L(f_j)+ L\left(\int_{\Y} f(\bx,\by)\ud\mu(\by)\right)$. 
It follows that 
    \[
    L(\phi) \ge  \sum_{j=0}^J \lambda_j f_j(\Phi_{\mJ}(L))+ \int_{\Y} f(\Phi_{\mJ}(L),\by)\ud\mu(\by) = \phi(\Phi_{\mJ}(L)).
    \]
    So, $\phi$ satisfies the Jensen inequality in the condition (C). Then, the conclusion follows using the same method of proof as in Lemma \ref{cor::sos}. 
\end{proof}

Now, we are ready to prove Theorem \ref{th::alt}.

\noindent{\itshape Proof of Theorem \ref{th::alt}}
\noindent We first prove that [Not(i) $\Longleftrightarrow$ (ii)]. 

\vskip 5pt
\noindent[Not(i) $\Longrightarrow$ (ii)] Suppose that (i) fails. That is,
there is no $\bu\in\RR^n$ such that $f_j(\bu)<0$ and $f(\bu,\by)<0$ for all $j\in[0:J]$ and $\by\in\Y$. We first claim that there is no $L\in\mathcal{B}[\bx]_m$ satisfying
\[
L(f_j)<0,\qquad j\in[0:J] \mbox{ and }
L(f(\cdot,\by))<0,\qquad \forall\,\by\in\Y.
\]
Indeed, if such an $L$ existed, then the  Jensen-type inequality in the condition (C) implies that 
\begin{equation}\label{eq::claim}
f_j(\Phi_{\mJ}(L))\le L(f_j)<0,\ \ \forall\ j\in[0:J],\ \
f(\Phi_{\mJ}(L), \by)\le L(f(\cdot,\by))<0,\ \ \forall\ \by\in\Y.
\end{equation}
That is, (i) holds true for the point $\Phi_{\mJ}(L)\in\RR^n$, a contradiction.

Recall that $C(\Y)$ denotes the space of continuous functions on $\Y$. Let
\[
C_{-}(\Y):=\{\gamma\in C(\Y) \mid \gamma(\by) \le 0, \forall \by\in\Y\}\quad\text{and}\quad
C_{--}(\Y):=\{\gamma\in C(\Y) \mid \gamma(\by) < 0, \forall \by\in\Y\}.
\]
Denote by $\RR^{J+1}_{-}$ (resp., $\RR^{J+1}_{--}$) the nonpositive (resp., negative) orthant
in $\RR^{J+1}$. 
%Denote by $C(\Y)$ the space of continuous functions $\gamma: \Y \to \RR$.
%Its dual is the space $\mathcal{M}(\Y)$ with the scalar 
%product $\mu\in\mathcal{M}(\Y)$ and $\gamma\in C(\Y)$ given by the integral 
%$\langle \mu, \gamma\rangle:=\int_{\Y} \gamma(\by)\ud\mu(\by)$.
Consider 
\[
K_1:=\left\{\left(L(f_0),\ldots,
L(f_J), L(f(\cdot,\by))\right) \ \big|\ L\in \mathcal{B}[\bx]_m\right\}
\quad\text{and}\quad 
K_2:= \RR^{J+1}_{-}\times C_{-}(\Y). 
%\mbox{ with }
%C_{-}(\Y):=\{\gamma\in C(\Y) \mid \gamma(\by) \le 0, \forall \by\in\Y\}.
\]
Clearly,  $K_1$ and $K_2$ are  convex subsets  of $\RR^{J+1}\times C(\Y)$ and 
$K_2$ is a convex cone with ${\rm int} K_2=\mathbb{R}^{J+1}_{--} \times C_{--}(\Y)$.
%where $C_{--}(\Y):=\{\gamma\in C(\Y) \mid \gamma(\by) < 0, \forall \by\in\Y\}$. 
Note that   
$K_1 \cap {\rm int} K_2=\emptyset$. By the convex separation theorem  \cite[Theorem 1.1.3]{zalinescu2002convex} and noting that $K_2$ is a cone, 
there exist $\td{\lambda}_j\ge 0$, $j\in[0:J]$, and $\td{\mu}\in\mathcal{M}_+(\Y)$, such that $\sum_{j=0}^J \td{\lambda}_j+\td{\mu}(\Y)>0$ and 
%\begin{equation}\label{eq::K1}
%\int_{\Y} \gamma(\by)\ud\td{\mu}(\by)\ge 0,\ \forall \gamma\in K_1,
%\end{equation}
%and 
\begin{equation}\label{eq::K2}
\sum_{j=0}^J \td{\lambda}_jL(f_j)+\int_{\Y} L(f(\cdot,\by))\ud\td{\mu}(\by)\ge 0,\ \
\forall L\in \mathcal{B}[\bx]_m.
\end{equation}
%By \eqref{eq::K1}, we have $\td{\mu}\in\mathcal{M}_+(\Y)$. 
For $\bu\in\RR^n$, considering the linear functional $L_{\bu}\in \mathcal{B}[\bx]_m$ 
with $L_{\bu}(\bx^{\ba})=\bu^{\ba}$ for all $\ba\in\N^n_m$, \eqref{eq::K2} implies that 
\[
\sum_{j=0}^J \td{\lambda}_jf_j(\bu)+\int_{\Y} f(\bu,\by)\ud\td{\mu}(\by)\ge 0, \ \forall \ \bu\in\RR^n.
\]
Let
\[
\lambda_j=\frac{\td{\lambda}_j}{\sum_{j=0}^{J}\td{\lambda}_j+\td{\mu}(\Y)}, \ \ 
j\in[0:J],\ \ \mu=\frac{\td{\mu}}{\sum_{j=0}^{J}\td{\lambda}_j+\td{\mu}(\Y)}.
\]
Then, $\lambda_j \ge 0$, $j \in[0:J]$, $\mu \in \mathcal{M}_+(\Y)$ and $\sum_{j=0}^{J}\lambda_j+\mu(\Y)=1$. Dividing $\sum_{j=0}^{J}\td{\lambda}_j+\td{\mu}(\Y)$ in \eqref{eq::K2} and using Lemma \ref{lem::sos}, one sees that  (ii) holds.

[(ii) $\Longrightarrow$ Not(i)] Let {(ii)} hold and suppose to the contrary that (i) also holds. Then,  there exists  $\bu\in\RR^n$ such that $f_j(\bu)<0$ and $f(\bu,\by)<0$ for all $j\in[0:J]$ and $\by\in\Y$. Then, (ii) implies that 
\[
\sum_{j=0}^J \lambda_j f_j(\bu)+\int_{\Y} f(\bu,\by)\ud\mu(\by) \ge 0.
\]
This is impossible as $\sum_{j=0}^J \lambda_j+\mu(\Y)=1$, $f_j(\bu)<0$ and $f(\bu,\by)<0$ for all $j\in[0:J]$ and $\by\in\Y$. Thus, the desired implication holds.   

Next, let $\nu$ be an arbitrary, fixed, finite nonnegative Borel measure with $\supp(\nu)=\Y$, and suppose that $\Theta_m(\bx)\in\mJ$. We now show [(ii) $\Leftrightarrow$ { (iii)}].  Noting that we have already shown that 
[Not (i) $\Leftrightarrow$ { (ii)}],  the desired conclusion follows if we show  [(ii) $\Longrightarrow$ (iii)] and [(iii) $\Longrightarrow$ Not(i)].

\vskip 5pt
\noindent[(ii) $\Longrightarrow$ (iii)]
Suppose that (ii) holds for $\lambda_j\ge 0$, $j\in[0:J]$, and $\mu\in\mathcal{M}_+(\Y)$. 
%Let $\td{\bx}=(x_1,\ldots,x_n,x_{n+1})$ and consider the homogenization
Let $\td{\bx}=(\bx,x_{n+1})=(x_1,\ldots,x_n,x_{n+1})$.
For each $j\in[0:J]$, as $\deg(f_j)\le m$, we can  write
$f_j(\bx)=\sum_{\ba\in\N^n_m}c_{j,\ba}\bx^\ba$. 
The homogenization of $f_j$  is a homogeneous
polynomial with variable $\td{\bx}$ and degree $m$, and is given by  
\[
\td{f}_j(\td{\bx})
=
\sum_{\ba\in\N^n_m}c_{j,\ba}\bx^\ba
x_{n+1}^{m-|\ba|}.
\]
Similarly,  for 
$f(\bx,\by)
=
\sum_{\ba\in\N^n_m}f_{\ba}(\by)\bx^{\ba}$, 
its homogenization is given by  
$\td{f}(\td{\bx},\by)
=
\sum_{\substack{\ba\in\N^n_m}}
f_{\ba}(\by) \bx^\ba
x_{n+1}^{m-|\ba|}$.   
In particular, for $x_{n+1} \neq 0$, one has     
\begin{equation}\label{eq::homo0}
\td{f}_j(\td{\bx})=x_{n+1}^m f_j\left(\frac{\bx}{x_{n+1}}\right),\ j\in[0:J],\ \td{f}(\td{\bx},\by)=x_{n+1}^m f\left(\frac{\bx}{x_{n+1}},\by\right).
    \end{equation}
Note that homogenization preserves SOS property. Then, (ii) implies that 
\[
\sum_{j=0}^J \lambda_j \td{f}_j(\td{\bx})+\int_{\Y}\td{f}(\td{\bx},\by)\ud\mu(\by)
\]
is an SOS polynomial in the variable $\td{\bx}$.
Write $\td{f}(\td{\bx},\by)=\sum_{\td{\ba}\in\N^{n+1}_{=m}}\td{f}_{\td{\ba}}(\by)\td{\bx}^{\td{\ba}}$.

Fix an $\varepsilon>0$. Consider the family 
$\{\td{f}_{\td{\ba}}(\by)\}_{\td{\ba}\in\N^{n+1}_{=m}}$ of finitely many polynomial functions
on $\Y$. It is clear that $\{\td{f}_{\td{\ba}}(\by)\}_{\td{\ba}\in\N^{n+1}_{=m}}$
is equicontinuous and pointwise bounded on $\Y$. Let
\begin{equation}\label{eq::N}
N:=\binom{n+m}{m}\max\{\vert\td{\bu}^{\td{\ba}}\vert\mid u_1^m+\cdots+u_n^m+u_{n+1}^m=1,\ \td{\ba}\in\N_{=m}^{n+1}\},
\end{equation}
which is well-defined. Recall that $\nu$ is a  finite nonnegative Borel measure with $\supp(\nu)=\Y$. 
By Proposition~\ref{prop::sigma}, 
there exists $\sigma_{\varepsilon}(\by)\in\Sigma[\by]$
such that $\int_{\Y}\sigma_{\varepsilon}(\by)\ud\nu(\by)=\mu(\Y)$ and 
\begin{equation}\label{eq::ineq2}
\left\vert\int_{\Y}\td{f}_{\td{\ba}}(\by)\ud\mu(\by)-\int_{\Y}\td{f}_{\td{\ba}}(\by)\sigma_{\varepsilon}(\by)\ud\nu(\by)\right\vert<\frac{\varepsilon}{N},\ \ \forall\ \td{\ba}\in \N_{=m}^{n+1}.
\end{equation}
For $j\in[0:J]$, let $\lambda_{\varepsilon,j}=\lambda_j$.
Then for all $\td{\bu}\in\RR^{n+1}$ with $\sum_{i=1}^{n+1}u_i^m=1$, 
as $\sum_{j=0}^J \lambda_j \td{f}_j(\td{\bu})+\int_{\Y}\td{f}(\td{\bu},\by)\ud\mu(\by)\ge 0$ by (ii), 
it holds that 
\begin{equation}\label{eq::nn1}
\begin{aligned}
\sum_{j=0}^J \lambda_{\varepsilon,j} \td{f}_j(\td{\bu})
+\int_{\Y} \sigma_{\varepsilon}(\by)\td{f}(\td{\bu}, \by)\ud\nu(\by)
&\ge \int_{\Y}\sigma_{\varepsilon}(\by)\td{f}(\td{\bu},\by)\ud\nu(\by)-\int_{\Y}\td{f}(\td{\bu},\by)\ud\mu(\by)\\
\ge &-\sum_{\td{\ba}\in\N_{=m}^{n+1}}\left\vert\int_{\Y}\td{f}_{\td{\ba}}(\by)\ud\mu(\by)
-\int_{\Y}\td{f}_{\td{\ba}}(\by)\sigma_{\varepsilon}(\by)\ud\nu(\by)\right\vert\vert \td{\bu}^{\td{\ba}}\vert\\
\ge&-\frac{\varepsilon}{N}\sum_{\td{\ba}\in\N_{=m}^{n+1}}\vert \td{\bu}^{\td{\ba}}\vert\ge -\varepsilon.
\end{aligned}
\end{equation}
Noting that $\td{f}_j(\td{\bx})$'s and $\td{f}(\td{\bx}, \by)$ are all homogeneous in $\td{\bx}$, it follows that 
for any $\td{\bx}\in\RR^{n+1}$, it holds that 
\begin{equation}\label{eq::nn2}
\sum_{j=0}^J \lambda_{\varepsilon,j} \td{f}_j(\td{\bx})
+\int_{\Y}\sigma_{\varepsilon}(\by)\td{f}(\td{\bx}, \by)\ud\nu(\by)+\varepsilon\sum_{i=1}^{n+1}x_i^m\ge 0.
\end{equation}
Therefore, for any $\bx\in\RR^n$, we have
\begin{equation}\label{eq::nn3}
\begin{aligned}
&\sum_{j=0}^J \lambda_{\varepsilon,j} f_j(\bx)
+\int_{\Y}\sigma_{\varepsilon}(\by)f(\bx, \by)\ud\nu(\by)+\varepsilon\Theta_m(\bx)\\
=&\sum_{j=0}^J \lambda_{\varepsilon,j} \td{f}_j(\bx,1)
+\int_{\Y}\sigma_{\varepsilon}(\by)\td{f}(\bx, 1, \by)\ud\nu(\by)
+\varepsilon\left(\sum_{i=1}^n x_i^m+1\right)\ge 0.
\end{aligned}
\end{equation}
%In fact, if $\bu=0$, then \eqref{eq::nn2} is clearly true as $f(\bx, \by)$ is homogeneous in $\bx$. 
%Otherwise, by \eqref{eq::nn1},
%\[
%\begin{aligned}
%    -\int_{\Y}\sigma_{\varepsilon}(\by)f(\bu, \by)\ud\nu(\by)+\varepsilon\sum_{i=1}^nu_i^m
%=\sum_{i=1}^nu_i^m\left(-\int_{\Y}\sigma_{\varepsilon}(\by)f\left(\frac{\bu}{\sqrt[m]{\sum_{i=1}^nu_i^m}}, \by\right)\ud\nu(\by)+\varepsilon\right)\ge 0.
%\end{aligned}
%\]
From the assumption $\Theta_m(\bx)\in\mJ$,  One obtains (iii) by using Lemma \ref{lem::sos}. 

\vskip 5pt
\noindent[(iii) $\Longrightarrow$ Not(i)] Let (iii) hold. We now see that {\rm (i)} fails.  
Suppose to the contrary that there exists $\bu\in\RR^n$ such that 
$f_j(\bu)<0$ and $f(\bu, \by)<0$ for all $j\in[0: J]$ and $\by\in\Y$. 
As $\Y$ is compact, we have $$ M:=\max\left\{
\max_{j\in[0:J]} f_j(\bu),
\max_{\by\in\Y}f(\bu,\by)
\right\}<0.$$ 
Fix an $\varepsilon>0$ such that $\varepsilon \, \Theta_m(\bu)<-M$.
Then by (iii), there exist $\lambda_{\varepsilon,j}\ge 0$, $j\in[0:J]$, and
$\sigma_{\varepsilon}\in\Sigma[\by]$ with 
$\sum_{j=0}^J\lambda_{\varepsilon,j}+\int_{\Y}\sigma_{\varepsilon}(\by)\ud\nu(\by)=1$ such that 
\[
\begin{aligned}
0&\le \sum_{j=0}^J \lambda_{\varepsilon,j} f_j(\bu)+\int_{\Y} \sigma_{\varepsilon}(\by)f\left(\bu,\by\right)\ud\nu(\by)+\varepsilon \Theta_m(\bu)\\
&\le M \left(\sum_{j=0}^{J}\lambda_{\varepsilon,j}+\int_{\Y}\sigma_{\varepsilon}(\by)\ud\nu(\by)\right)+\varepsilon \Theta_m(\bu)\\
&=M+\varepsilon \Theta_m(\bu)<0.
\end{aligned}
\]
a contradiction. Thus, the conclusion follows.
\hfill \qed

\section{Nonnegativity certificates for Jensen-type polynomials with UQs}\label{sec::certificates}
Throughout this section, we work under Assumption~\ref{assump::0}.
We derive SOS-based certificates for the nonnegativity of a polynomial $f_0$
on a feasible set $\mathcal X$ described by finitely many polynomial inequalities
and an uncountable family of polynomial inequalities: 
\[
\mathcal X
:=
\left\{
\bx\in\mathbb{R}^n
\ \middle|\
f_j(\bx)\leq 0,\ j\in[1:J],\quad
f(\bx,\by)\leq 0,\ \forall\,\by\in\mathcal Y
\right\},
\]
where $f_j \in \mathbb{R}[\bx]$, $f \in \mathbb{R}[\bx,\by]$ and $\Y \subseteq \RR^\ny$ is a compact set. 
We collect the constraint functions in the family
\[
\bff
:=
\big\{
f_1(\bx),\ldots,f_J(\bx),f(\bx,\by)
\big\}.
\]
\begin{definition}[Slater condition]\label{def::slater}
We say that the Slater condition holds for $\bff$ if there exists $\bar{\bu}\in\RR^n$ such that
   $f_j(\bar{\bu})<0$ for all $j\in[1:J]$,  and  $f(\bar{\bu}, \by)<0$ for all $\by\in\Y$. In this case,  $\bar{\bu}$ is called a Slater point of the set $\X$.
\end{definition}

For the constraint family $\bff$, define
\begin{equation}\label{eq::Q}
\mC(\bff):=\left\{
\begin{array}{c|c}
\displaystyle  s(\bx)-\sum_{j=1}^J\lambda_j f_j(\bx)-\int_{\Y} f(\bx,\by)\ud\mu(\by)   &  
\displaystyle  \begin{aligned}
    &\lambda_j\ge 0,\ j\in[1:J],\\
    &\mu\in\mathcal{M}_+(\Y), \ s(\bx)\in\Sigma[\bx]_m
  \end{aligned}
\end{array}
\right\}.
\end{equation}
For a fixed finite and nonnegative measure $\nu\in\mathcal{M}_+(\Y)$ with ${\rm supp}(\nu)=\Y$, let
\begin{equation}\label{eq::Qv}
\mC(\bff,\nu):=\left\{
\begin{array}{c|c}
 \displaystyle s(\bx)-\sum_{j=1}^J\lambda_j f_j(\bx)-\int_{\Y} \sigma(\by)f(\bx,\by)\ud\nu(\by)   &  
\displaystyle  \begin{aligned}
    &\lambda_j\ge 0,\ j\in[1:J],\\
    &\sigma(\by)\in\Sigma[\by],\ s(\bx)\in\Sigma[\bx]_m 
  \end{aligned}
\end{array}
\right\}.
\end{equation}
Direct verification shows that every polynomial belonging to the cones $\mC(\bff)$ or $\mC(\bff,\nu)$ is nonnegative on $\X$.
As an application of Theorem \ref{th::alt}, we now derive the following 
certificates for the nonnegativity of $f_0$ over the set $\X$.
\begin{theorem}\label{th::Pastz}
Suppose that  Assumption \ref{assump::0} holds and   the Slater condition is satisfied for $\bff$. 
%    For a polynomial $f_0\in\RR[\bx]$ of degree $m$, assume that $-f_0$ has ENNCs.
    If $f_0(\bx)\ge 0$ on $\X$, then the following statements hold.
    \begin{enumerate}
        \item[\upshape (i)] $f_0(\bx)\in \mC(\bff)$;
        \item[\upshape (ii)] Let $\nu$ be a fixed finite and nonnegative measure $\nu\in\mathcal{M}_+(\Y)$ with ${\rm supp}(\nu)=\Y$. Suppose further that { $\Theta_m \in \mJ$.} Then,
        $f_0(\bx)+\varepsilon\Theta_m(\bx)\in \mC(\bff,\nu)$ for any $\varepsilon>0$.
    \end{enumerate}
\end{theorem}
\begin{proof}
    As $f_0(\bx)\ge 0$ on $\X$, then 
\begin{equation}\label{eq::infeasible0}
\{\bx\in\RR^n \mid f_0(\bx)<0,\ f_1(\bx)<0,\ldots, f_J(\bx)<0, \  
        f(\bx,\by)<0, \ \forall \by\in\Y\} =\emptyset. 
\end{equation}
Let $\bar{\bu}$ be a Slater point of $\X$. As $\Y$ is compact, one has 
\[ 
M_0:=\max\left\{\max_{j\in[1:J]} f_j(\bar{\bu}),\  \max_{\by\in\Y} f(\bar{\bu},\bv)\right\}<0.
\]
%{\color{blue} Replace $\xi$ by $M_0$}
%\vskip 5pt
\noindent [Proof of (i)] Applying Theorem \ref{th::alt} to \eqref{eq::infeasible0} yields  $\lambda_j\ge 0$, 
$j=0,1,\ldots,J$, and $\mu\in\mathcal{M}_+(\Y)$ with $\sum_{j=0}^J\lambda_j+\mu(\Y)=1$ 
such that 
        \begin{equation}\label{eq::feasible0}
        \lambda_0 f_0(\bx)+\sum_{j=1}^J\lambda_j f_j(\bx)+\int_{\Y} f(\bx,\by)\ud\mu(\by)=s(\bx)
        \end{equation}
for  some $s\in\Sigma[\bx]_m$. 
We see that $\lambda_0\neq 0$  (otherwise we have $\sum_{j=1}^J\lambda_j+\mu(\Y)=1$ and 
\[
\sum_{j=1}^J\lambda_j f_j(\bx)+\int_{\Y} f(\bx,\by)\ud\mu(\by)=s(\bx)\ge 0,
\]
for all $\bx\in\RR^n$. Then, for the Slater point $\bar{\bu}$, we have
\[
0\le s(\bar{\bu})=\sum_{j=1}^J\lambda_j f_j(\bar{\bu})+\int_{\Y} f(\bar{\bu},\by)\ud\mu(\by)
\le M_0 \left(\sum_{j=1}^J\lambda_j+\mu(\Y)\right)=M_0<0,
\]
a contradiction). Therefore, $\lambda_0>0$ and hence, by dividing $\lambda_0$ in \eqref{eq::feasible0}, we see that $f_0(\bx)\in \mC(\bff)$. 

\vskip 5pt
\noindent [Proof of (ii)] Fix an $\varepsilon'>0$, we next prove that 
$f_0(\bx)+\varepsilon'\Theta_m(\bx)\in \mC(\bff,\nu)$.

For any $\varepsilon>0$,  applying Theorem \ref{th::alt} to \eqref{eq::infeasible0} yields
$\lambda_{\varepsilon,j}\ge 0$, $j=0,1,\ldots,J$, and
        $\sigma_{\varepsilon}\in\Sigma[\by]$
        with $\sum_{j=0}^J\lambda_{\varepsilon,j}+\int_{\Y}\sigma_{\varepsilon}(\by)\ud\nu(\by)=1$
        such that 
        \begin{equation}\label{eq::epsilon0}
            \lambda_{\varepsilon,0} f_0(\bx) +\sum_{j=1}^J\lambda_{\varepsilon,j} f_j(\bx)+\int_{\Y} \sigma_{\varepsilon}(\by)f(\bx, \by)\ud\nu(\by)+\varepsilon\Theta_m(\bx)=s(\bx),
        \end{equation}
for some $s\in\Sigma[\bx]_m$.
We claim that there exist $\tau>0$ and $\bar{\varepsilon}>0$ such that \eqref{eq::epsilon0} holds, and 
$\lambda_{\varepsilon,0}>\tau$ for all $\varepsilon<\bar{\varepsilon}$. 
Otherwise, there exists a sequence $\{\varepsilon_k\}_{k\in\N}$ 
with $\lim_{k\to\infty}\varepsilon_k=0$ and $\lambda_{\varepsilon_k,0}<1/k$ such that 
   $\sigma_{\varepsilon_k}\in\Sigma[\by]$, 
 $\sum_{j=0}^J\lambda_{\varepsilon_k,j}+\int_{\Y}\sigma_{\varepsilon_k}(\by)\ud\nu(\by)=1$ and \[s_k(\bx):=\lambda_{\varepsilon_k,0} f_0(\bx) +\sum_{j=1}^J\lambda_{\varepsilon_k,j} f_j(\bx)+\int_{\Y} \sigma_{\varepsilon_k}(\by)f(\bx, \by)\ud\nu(\by)+\varepsilon_k \, \Theta_m(\bx) \in \Sigma[\bx]_m.
            \]
Then, for the Slater point $\bar{\bu}$, using the definition of $M_0$ and the fact that $\sigma_{\varepsilon_k}\in\Sigma[\by]$, 
 $\sum_{j=1}^J\lambda_{\varepsilon_k,j}+\int_{\Y}\sigma_{\varepsilon_k}(\by)\ud\nu(\by)=1-\lambda_{\varepsilon_k,0} \ge 0$, it follows that 
\[
        \lambda_{\varepsilon_k,0} f_0(\bar{\bu})+\varepsilon_k\Theta_m(\bar{\bu}) =s_k(\bar{\bu})-\sum_{j=1}^J\lambda_{\varepsilon_k,j} f_j(\bar{\bu})-\int_{\Y} \sigma_{\varepsilon_k}(\by)f(\bar{\bu}, \by)\ud\nu(\by)\ge -(1-\lambda_{\varepsilon_k,0}) M_0.
        \]
Letting $k\to\infty$ yields a contradiction $0\ge -M_0>0$. 

Now, choose an $\varepsilon>0$ such that $\varepsilon<\min\{\tau\varepsilon', \bar{\varepsilon}\}$, then by \eqref{eq::epsilon0},
\[
\frac{\varepsilon}{\lambda_{\varepsilon,0}}\le\frac{\varepsilon}{\tau}\le \varepsilon'\quad
\text{and}\quad f_0(\bx)+\frac{\varepsilon}{\lambda_{\varepsilon,0}}\Theta_m(\bx)\in \mC(\bff,\nu).
\]
Then, as $\Theta_m(\bx)\in\Sigma[\bx]_m$, by the definition of $\mC(\bff,\nu)$, we have 
\[
f_0(\bx)+\varepsilon'\Theta_m(\bx)=f_0(\bx)+\frac{\varepsilon}{\lambda_{\varepsilon,0}}\Theta_m(\bx)
+\left(\varepsilon'-\frac{\varepsilon}{\lambda_{\varepsilon,0}}\right)\Theta_m(\bx)\in \mC(\bff,\nu).
\]
This completes the proof. 
\end{proof}

%\subsection{Finite dimensional truncations}
Note that the cones $\mathcal{M}_+(\Y)$ and $\Sigma[\by]$ in the representations of 
$\mC(\bff)$ and $\mC(\bff,\nu)$ are both infinite dimensional. As a result, the nonnegative certificates  in Theorem \ref{th::Pastz} are still, in general, hard to check. Next, we construct tractable approximations 
of $\mC(\bff)$ and $\mC(\bff,\nu)$, which enable us to provide SOS approximations for the nonnegativity of
$f_0$ over $\X$.

\vskip 5pt

\noindent{\bf Tractable approximation for $\mC(\bff,\nu)$.} We first consider constructing tractable approximation for the cone $\mC(\bff,\nu)$. 
For the fixed measure $\nu\in\mathcal{M}_+(\Y)$ and any $k\in\N$, bounding the degree of the SOS polynomial $\sigma(\by)$ by $2k$, we let
\begin{equation}\label{eq::Qkv}
\mC_k(\bff,\nu):=\left\{
\begin{array}{c|c}
 \displaystyle s(\bx)-\sum_{j=1}^J\lambda_j f_j(\bx)-\int_{\Y} \sigma(\by)f(\bx,\by)\ud\nu(\by)   &  
\displaystyle  \begin{aligned}
    &\lambda_j\ge 0,\ j\in[1:J],\\
    &\sigma(\by)\in\Sigma[\by]_{2k},\ s(\bx)\in\Sigma[\bx]_m 
  \end{aligned}
\end{array}
\right\}.
\end{equation}
Clearly, $\mC_k(\bff,\nu)$, $k \in \mathbb{N}$, are {\itshape inner} approximations of $\mC(\bff,\nu)$ and 
$\mC(\bff,\nu)=\cup_{k\in\N}\mC_k(\bff,\nu)$. 
More importantly, for given $\varepsilon>0$ and $k\in\N$, 
checking whether $f_0(\bx)+\varepsilon\Theta_m(\bx)\in \mC_k(\bff,\nu)$ holds or not amounts to solving   
an SDP feasibility problem.  

We now have the following tractable approximation on nonnegativity of $f_0$ over $\X$ based on the inner approximation cone $\mC_k(\bff,\nu)$. 
\begin{corollary}\label{cor::Qkv}
  Under Assumption \ref{assump::0}, suppose that 
    the Slater condition holds for $\bff$, $\supp(\nu)=\Y$ and $\Theta_m \in \mJ$. 
  If $f_0(\bx)\ge 0$ on $\X$, 
  then the following statements hold:
  \begin{enumerate}[\upshape (i)]
      \item For any $\varepsilon>0$, there exists $k_{\varepsilon}\in\N$ such that 
  $f_0(\bx)+\varepsilon\Theta_m(\bx)\in \mC_k(\bff,\nu)$ for any $k\ge k_{\varepsilon}$;
  \item Suppose that $\sum_{i=1}^n x_i^m-R\in\bff$ for some $R>0$. Then, for any $\varepsilon>0$, 
  there exists $k_{\varepsilon}\in\N$ such that 
  $f_0(\bx)+\varepsilon\in \mC_k(\bff,\nu)$ for any $k\ge k_{\varepsilon}$.
  \end{enumerate}
\end{corollary}
\begin{proof}
    (i) This part follows directly by  Theorem \ref{th::Pastz} (ii).

    (ii) Let $\varepsilon'=\varepsilon/(R+1)$. By (i), 
    there exists $k_{\varepsilon'}\in\N$ such that 
  $f_0(\bx)+\varepsilon'\Theta_m(\bx)\in \mC_k(\bff,\nu)$ for any $k\ge k_{\varepsilon'}$. Note that
  \[
       f_0(\bx)+\varepsilon'\Theta_m(\bx)=f_0(\bx)+\varepsilon'\left(\sum_{i=1}^n x_i^m+1-R+R\right)
       =f_0(\bx)+\varepsilon+\varepsilon'\left(\sum_{i=1}^n x_i^m-R\right).
  \]
  Since $\sum_{i=1}^n x_i^m-R\in\bff$, we have $-\varepsilon'\left(\sum_{i=1}^n x_i^m-R\right)\in\mC_k(\bff,\nu)$ for any $k\ge k_{\varepsilon'}$. Letting $k_{\varepsilon}=k_{\varepsilon'}$, the conclusion follows.
\end{proof}

\noindent{\bf Tractable approximation for $\mC(\bff)$.}  Next, we focus on constructing tractable approximations of $\mC(\bff)$. To do this, we further assume that the parameter set $\Y$ is a semialgebraic set which takes the form that 
\begin{equation}\label{eq::Y}
\Y=\{\by\in\RR^\ny\mid g_1(\by)\ge 0, \ldots, g_r(\by)\ge 0\},
\end{equation}
for some $g_1,\ldots,g_r\in\RR[\by]$.
Write
\[
f(\bx,\by)=\sum_{\ba\in\supp_{\bx}(f)}f_{\ba}(\by)\bx^{\ba}=\sum_{\bb\in\supp_{\by}(f)}f_{\bb}(\bx)\by^{\bb},
\]
where
\[
\begin{aligned}
\suppx{\bx}{f}&\coloneqq\{\ba\in\N^n \mid \bx^{\ba}\by^{\bb}\ \text{appears in }f(\bx, \by) \text{ for some }\bb\in\N^\ny\},\\
\suppx{\by}{f}&\coloneqq\{\bb\in\N^\ny \mid \bx^{\ba}\by^{\bb}\ \text{appears in }f(\bx, \by) \text{ for some }\ba\in\N^n\}.
\end{aligned}
\]
For a (pseudo-moment) sequence of real numbers 
$\bz=(z_{\bb})_{\bb\in \mathbb{N}^{\ny}_{2k}}$, $k\in\N$, define a linear functional
$L_{\bz}: \RR[\by]_{2k}\to\RR$ by letting $L_{\bz}(\phi)=\sum_{\bb\in\supp(\phi)}\phi_{\bb}z_{\bb}$ for any 
$\phi(\by)=\sum_{\bb\in\supp(\phi)}\phi_{\bb}\by^{\bb}\in\RR[\by]_{2k}$.
Let 
\[
L_{\bz}(f(\bx,\by))=\sum_{\ba\in\supp_{\bx}(f)}L_{\bz}(f_{\ba}(\by))\bx^{\ba}\in\RR[\bx].
\]
Similarly, for a (pseudo-moment) sequence of real numbers 
$\bw=(w_{\ba})_{\ba\in \mathbb{N}^n_m}$, define a linear functional
$L_{\bw}: \RR[\bx]_m\to\RR$ such that $L_{\bw}(h)=\sum_{\ba\in\supp(h)}h_{\ba}w_{\ba}$ for any 
$h(\bx)=\sum_{\ba\in\supp(h)}h_{\ba}\bx^{\ba}\in\RR[\bx]_m$.
Let 
\[
L_{\bw}(f(\bx,\by))=\sum_{\bb\in\supp_{\by}(f)}L_{\bw}(f_{\bb}(\bx))\by^{\bb}\in\RR[\by].
\]
Let $r \ge 1$, 
%$\bg:=\{g_1,\ldots,g_r\}$, $d_w:=\lceil\deg(g_w)/2\rceil$, $w\in[1:r]$, and 
\begin{equation}\label{eq::degree}
  \begin{aligned}
\bg&:=\{g_1,\ldots,g_r\},\ d_w:=\lceil\deg(g_w)/2\rceil,\ w\in[1:r],\\
    d_{\bg}&:=\max\{d_w \, \mid \,  w\in[1:r]\},\
d_0:=\max\{d_{\bg},\ \lceil\deg_{\by}(f)/2\rceil\}.
\end{aligned}  
\end{equation}
For $k\ge d_0$, consider the cone 
\begin{equation}\label{eq::Qk}
\mC_k(\bff):=\left\{
\begin{array}{c|c}
\displaystyle  s(\bx)-\sum_{j=1}^J\lambda_j f_j(\bx)-L_{\bz}(f(\bx,\by))   &  
\displaystyle  \begin{aligned}
    &s(\bx)\in\Sigma[\bx]_m,\ \lambda_j\ge 0, j\in[1:J],\ \bz\in\RR^{\vert\mathbb{N}^{\ny}_{2k}\vert},\\ 
    &\bM_k(\bz)\succeq 0,\ \bM_{k-d_w}(g_w\bz)\succeq 0,\ w\in[1:r]
  \end{aligned}
\end{array}
\right\}.
\end{equation}

For any $\mu\in\mathcal{M}_+(\Y)$, $k\ge d_0$ and $\bb\in\N^\ny_{2k}$, 
let $\bz=(z_{\bb})_{\bb\in \mathbb{N}^{\ny}_{2k}}$ with each $z_\bb=\int_{\Y}\by^{\bb}\ud\mu(\by)$.
Then, it holds that
\[
\bM_k(\bz)\succeq 0,\ \bM_{k-d_w}(g_w\bz)\succeq 0,\ w\in[1:r],\ 
L_{\bz}(f(\bx,\by))=\int_{\Y}f(\bx,\by)\ud\mu(\by).
\]
Thus, $\mC(\bff) \subseteq \mC_k(\bff)$. That is, $\mC_k(\bff)$, $k \ge d_0$, are {\itshape outer} approximations of $\mC(\bff)$. Moreover, checking whether 
$f_0(\bx)\in \mC_k(\bff)$, for a given $k\ge d_0$, can also be equivalently rewritten as an SDP feasibility problem. 

Next, we state the following tractable approximation on nonnegativity of $f_0$ over $\X$ based on the outer approximation cone $\mC_k(\bff)$. 
\begin{theorem}\label{th::Qk}
The following statements hold.
    \begin{enumerate}
        \item[\upshape (i)] Suppose that there exists $k$ with  $k\ge d_0$ such that $f_0\in \mC_k(\bff)$ and the rank condition 
\begin{equation}\label{eq::FEC0}
       {\exists\ t\quad\text{s.t. }  d_0\le t\le k\quad \text{and}\quad \rank\ \bM_t(\bz)=\rank\ \bM_{t-d_{\bg}}(\bz)}
      \end{equation}
    holds for the vector $\bz\in\RR^{\vert\mathbb{N}^{\ny}_{2k}\vert}$ in the representation of $f_0$
    in \eqref{eq::Qk}.  Then $f_0(\bx)\ge 0$ on $\X$;
    \item[\upshape (ii)]
       Under Assumption \ref{assump::0}, suppose that 
    the Slater condition holds for $\bff$. Let $\bg=\{g_1,\ldots,g_r\}$ and suppose that the quadratic module $\mQ(\bg)$ is Archimedean.  
    Then $f_0(\bx)\ge 0$ on $\X$ if and only if $f_0(\bx)\in \mC_k(\bff)$ for all $k\ge d_0$.
    \end{enumerate} 
\end{theorem}
\begin{proof}  

\text{\upshape(i)}
Suppose that $f_0\in\mC_k(\bff)$ for some $k\geq d_0$ and that the
rank condition \eqref{eq::FEC0} holds for some $t\in[d_0:k]$. 
%By the definition of
%$\mC_k(\bff)$,
%\[
%\bM_k(\bz)\succeq 0,
%\qquad
%\bM_{k-d_w}(g_w\bz)\succeq 0,
%\quad w\in[1:r].
%\]
%Consequently, the corresponding truncated moment and localizing
%matrices at level $t$ are also positive semidefinite:
%\[
%\bM_t(\bz)\succeq 0,
%\qquad
%\bM_{t-d_w}(g_w\bz)\succeq %0,
%\quad w\in[1:r].
%\]
By the flat truncation theorem
\cite{CF2000,NieFlatTruncation} and the rank condition \eqref{eq::FEC0},  the truncated sequence
$\bz|_{2t}$ admits a finite nonnegative representing measure $\mu$
supported on $\Y$. Moreover, as $t\geq d_0$, 
all moments required to evaluate $L_{\bz}(f(\bx,\by))$ are contained
in $\bz|_{2t}$, and so, 
$L_{\bz}(f(\bx,\by))
=
\int_{\Y}f(\bx,\by)\,d\mu(\by)$. 
Since $f_0\in\mC_k(\bff)$, there exist
$s\in\Sigma[\bx]_m$ and $\lambda_j\geq 0$, $j\in[1:J]$, such that
\[
f_0(\bx)
=
s(\bx)
-\sum_{j=1}^{J}\lambda_j f_j(\bx)
-L_{\bz}(f(\bx,\by))
=
s(\bx)
-\sum_{j=1}^{J}\lambda_j f_j(\bx)
-\int_{\Y}f(\bx,\by)\,d\mu(\by).
\]
Therefore, for every $\bx\in\X$, 
$f_0(\bx)
\geq 0$,
because $s(\bx)\geq 0$, $\lambda_j\geq 0$, $f_j(\bx)\leq 0$,
and $f(\bx,\by)\leq 0$ for all $\by\in\Y$.

\medskip

\text{\upshape(ii)}
Suppose first that $f_0(\bx)\geq 0$ for every $\bx\in\X$. By
Theorem~\ref{th::Pastz}\textup{(i)}, we have
$f_0\in\mC(\bff)$. For every $k\geq d_0$, let
$\mu\in\mathcal M_+(\Y)$ 
be a measure appearing in a representation of $f_0$ in
$\mC(\bff)$, and define
\[
z_{\bb}:=\int_{\Y}\by^{\bb}\,d\mu(\by),
\qquad
\bb\in\mathbb N^\ny_{2k}.
\]
Then, 
$\bM_k(\bz)\succeq 0,\ \bM_{k-d_w}(g_w\bz)\succeq 0,\ w\in[1:r]$ and 
$L_{\bz}(f(\bx,\by))=\int_{\Y}f(\bx,\by)\ud\mu(\by)$. 
Thus, $f_0\in\mC_k(\bff)$ for every $k\geq d_0$.

Conversely, suppose that
$f_0\in\mC_k(\bff)$ for all $k\geq d_0$. 
We now see that $f_0(\bx)\geq 0$ for every $\bx\in\X$. 
Assume, to the contrary, that there exists
$\hat{\bu}\in\X$ such that
$f_0(\hat{\bu})<0$. 
Let $\bar{\bu}\in\X$ be a Slater point. Since
$f_0(\hat{\bu})<0$, we can choose $\varepsilon>0$ sufficiently small
such that
$\varepsilon f_0(\bar{\bu})+f_0(\hat{\bu})<0$. 
Define the linear functional
\[
L_\varepsilon
:=
\varepsilon L_{\bar{\bu}}+L_{\hat{\bu}},
\]
where $L_{\bu}$ denotes the point-evaluation functional
$L_{\bu}(q):=q(\bu)$ for all $q \in \RR[\bx]_m$. Direct verification shows that $L_{\varepsilon}(q) \ge 0$ for any SOS polynomial $q$, and $L_{\varepsilon}(f_j)=\varepsilon f_j(\bar \bu)+ f_j(\hat{\bu})<0$, $j \in [1:J]$ (because $\bar{\bu}$ is a Slater point and $\hat{\bu}\in\X$). Regarding $f(\bx,\by)$ as a polynomial in $\bx$, we have
$L_\varepsilon(f(\bx,\by))
=
\varepsilon f(\bar{\bu},\by)
+
f(\hat{\bu},\by) \in \mathbb{R}[\by]$. 
%Since $\bar{\bu}$ is a Slater point and $\hat{\bu}\in\X$, one has 
%$f(\bar{\bu},\by)<0$ and 
%$f(\hat{\bu},\by)\leq 0$ for all %$\by\in\Y$.  
%Consequently, 
Similarly, one has 
$L_\varepsilon(f(\bx,\by))<0$ for all 
$\by\in\Y$. In other words, 
$H_{\varepsilon}(\by):=-L_\varepsilon(f(\bx,\by))$ 
is positive over $\Y$. 
By Putinar's Positivstellensatz and the Archimedean assumption on
$\mQ(\bg)$, there exists $k_0\geq d_0$ such that
$H_{\varepsilon}(\by)=-L_\varepsilon(f(\bx,\by))
\in\mQ_{k_0}(\bg)$. 

Since $f_0\in \mC_{k_0}(\bff)$, there exist
    $s(\bx)\in\Sigma[\bx]_m$, $\lambda_j\ge 0$, $j\in[1:J]$ and $\bz\in\RR^{\vert\mathbb{N}^{\ny}_{2k_0}\vert}$ with 
    \begin{equation}\label{eq::moment}
    \bM_{k_0}(\bz)\succeq 0,\ \bM_{k_0-d_w}(g_w\bz)\succeq 0,\ w\in[1:r],
    \end{equation}
    such that 
    $f_0(\bx)=s(\bx)-\sum_{j=1}^J\lambda_j f_j(\bx)-L_{\bz}(f(\bx,\by))$. 
    As $H_{\varepsilon}(\by)\in Q_{k_0}(\bg)$, the condition \eqref{eq::moment} implies that 
    \[
    -L_{\varepsilon}(L_{\bz}(f(\bx,\by)))=L_{\bz}(-L_{\varepsilon}(f(\bx,\by))) =L_{\bz}(H_{\varepsilon})\ge 0.
    \]
    Noting that $L_{\varepsilon}(s(\bx)) \ge 0$ (thanks to the fact that $s$ is an SOS polynomial) and $L_{\varepsilon}(f_j)<0$, it follows that 
    \[
L_{\varepsilon}(f_0(\bx))=
    L_{\varepsilon}(s(\bx))-\sum_{j=1}^J\lambda_j L_{\varepsilon}(f_j(\bx))-L_{\varepsilon}(L_{\bz}(f(\bx,\by)))\ge 0.
    \]
On the other hand, direct verification shows that    $ L_{\varepsilon}(f_0(\bx))=\varepsilon f_0(\bar{\bu})+f_0(\hat{\bu})<0$,     a contradiction. Hence, we have $f_0(\bx)\ge 0$ on $\X$.
\end{proof}

\section{Applications and Numerical Examples}\label{sec::sipop}
{In this section, we provide two applications illustrating our developments in the previous sections, and present some numerical examples. All computations in the sequel were performed and timed on a desktop 
with 8-Core Intel i7 3.8 GHz CPUs and 8 GB RAM.}
\subsection{} \noindent{\bf Application 1: source localization problems.} Consider the source localization problem: 
$$\min_{\bx \in \RR^n}f_{\rm SL}(\bx)= \sum_{\ell=1}^p (\|\bx-\bv_\ell\|_U^2 -r_\ell^2)^{2d},$$
where $d \ge 1$,  $p \in \mathbb{N}$ is the number of sensors, for each $\ell=1,\ldots,p$, ${\bv}_\ell \in \RR^n$ is the given location of the sensor $\ell$, $r_\ell > 0$ is the measured distance to each sensor $\ell$ corrupted by noise, and $\|\cdot\|_U$ is the norm induced by a positive definite symmetric matrix $U$ and is given by $\|\bv\|_U=\sqrt{\bv^TU\bv}$ for all $\bv$. 
Note that  $f_{\rm SL}=g \circ h$ which is a polynomial with degree $m=4d$, where $g(\bxi)=\sum_{\ell=1}^p \xi_\ell^{2d}$ and 
$h(\bx)=(h_1(\bx),\ldots,h_p(\bx))$ with 
\[
h_\ell(\bx)=\|\bx-\bv_\ell\|_U^2 -r_\ell^2=\bx^TU\bx-2(U \bv_\ell)^T\bx+(\bv_\ell^TU\bv_\ell -r_\ell^2), \ \ell=1,\ldots,p.
\]
We now suppose the following {\it  dimension condition} holds:  ${\rm dim} \, {\rm span}\{\bv_2-\bv_1,\ldots,\bv_p-\bv_1\} \le n-1$ (in the case $n=3$, 
this assumption means that the sensor locations $\bv_\ell$, $\ell=1,\ldots,p$ 
are co-planar).  As $U$ is positive definite,  ${\rm dim} \, {\rm span}\{U(\bv_2-\bv_1),\ldots,U(\bv_p-\bv_1)\} \le n-1$, and so, the assumption in Proposition \ref{prop::subsclass4} (c) is satisfied. This implies that $f_{\rm SL}$ is a Jensen-type polynomial with a Jensen map $\Phi_{\mJ^{\csosc}_{h,4d}}$ defined as in Proposition \ref{prop::subsclass4}.

Moreover, under this dimension condition,  the corresponding Jensen map has an explicit form. To see this, let 
$\overline{\bd} \in \mathbb{R}^n \backslash\{\bf 0\}$ be such that $(U\overline{\bd})^T(\bv_\ell-\bv_1)=0$, $\ell=2,\ldots,p$. Recall that, for any $L\in\mathcal{B}[\bx]_{4d}$, 
%then let $\varphi_3(L)$ be an arbitrarily fixed point in $\RR^n$; otherwise, 
 $\Phi_{\mJ^{\csosc}_{h,4d}}(L)=\bu$ where $\bu=(u_i)_{i\in[1:n]}\in\RR^n$ satisfies 
$h_\ell(\bu)=L(h_\ell(\bx))$ for all $\ell \in[1:p]$. Now, let $\overline{\bx}=(L(x_1),\ldots,L(x_n))$. Let $h_0(\bx)=\bx^TU\bx$. Direct verification shows that 
\[
\gamma_L:=h_0(\overline{\bx})- L(h_0)=h_\ell(\overline{\bx})-L(h_\ell), \ \forall \ \ell=1,\ldots,p. 
\]
Clearly, $\gamma_L \le 0$ (as $h_0$ is a convex quadratic function, and so, Proposition \ref{prop::Jensen} applies). 
Consider 
$\bx(t)=\overline{\bx}+t\, \overline{\bd}$ with $t \in \mathbb{R}$. Then, for all $\ell=1,\ldots,p$, 
\[
\begin{aligned}
h_\ell(\overline{\bx}+t\overline{\bd})-L(h_\ell)
&= h_\ell(\overline{\bx})-L(h_\ell)
   + h_\ell(\overline{\bx}+t\overline{\bd})-h_\ell(\overline{\bx})\\
&= \gamma_L+t^2\overline{\bd}^{\,T}U\overline{\bd}
   +2t\bigl(U(\overline{\bx}-\bv_\ell)\bigr)^T\overline{\bd}\\
&= \gamma_L+t^2\overline{\bd}^{\,T}U\overline{\bd}
   +2t\bigl(U(\overline{\bx}-\bv_1)\bigr)^T\overline{\bd},
\end{aligned}
\] 
% \begin{eqnarray*}
% h_\ell(\bx(t))- L(h_\ell) &=&  (\overline{\bx}+t\, \overline{\bd})^TU(\overline{\bx}+t\, \overline{\bd})-2(U^T \bv_\ell)^T(\overline{\bx}+t\, \overline{\bd})+(\bv_\ell^TU\bv_\ell -r_\ell^2)- L(h_\ell) \\
% &= & (\overline{\bx}+t\, \overline{\bd})^TU(\overline{\bx}+t\, \overline{\bd})-2(U^T \bv_1)^T(\overline{\bx}+t\, \overline{\bd})+(\bv_\ell^TU\bv_\ell -r_\ell^2)- L(h_\ell) \\
% &=&  t^2 \overline{\bd}^TU \overline{\bd} + 2t \, (U^T(\overline{\bx}-\bv_1))^T\overline{\bd} +\gamma_L, 
% \end{eqnarray*}
where the last equality holds as $(U\overline{\bd})^T(\bv_\ell-\bv_1)=0$, $\ell=2,\ldots,p$. 
Let $t_L$ be a solution of the quadratic equation  
$t^2 \overline{\bd}^TU \overline{\bd} + 2t \, (U(\overline{\bx}-\bv_1))^T\overline{\bd} +\gamma_L=0$.  
Note that such a solution always exists as $U$ is positive definite and $\gamma_L \le 0$, and has a closed-form. Then, we can set  
\begin{equation}\label{eq:explicit_Jensen}
    \Phi_{\mJ^{\csosc}_{h,4d}}(L)= \bx(t_L)= \overline{\bx}+t_L\, \overline{\bd}.
\end{equation}

Note that $f_{\rm SL}$ is coercive, and so, $f^{\star}:=\min_{\bx \in \RR^n}f_{\rm SL}(\bx)$ is 
well-defined and is finite, and a global minimizer exists. It follows from Proposition~\ref{prop::subsclass4}{(ii)} that 
\begin{equation}\label{eq::sos4f_SL}
   {f^{\star}=}\sup\{\gamma: f_{\rm SL}-\gamma \in \Sigma[{\bx}]_{4d}\}.
\end{equation} 
Note that this SOS {reformulation} problem can be cast as an SDP, and its dual problem reads as 
\begin{equation}\label{eq::moment4f_SL}
    \inf\{ L(f_{\rm SL}):  L\in(\Sigma[\bx]_{4d})^*, \ L(1)=1\},
\end{equation}
which can also be cast as an SDP. Notice that the duality gap between \eqref{eq::sos4f_SL} and its dual \eqref{eq::moment4f_SL} is zero (as Slater condition holds for \eqref{eq::moment4f_SL}).
Using the Jensen inequality (C) for the Jensen-type polynomial $f_{\rm SL}$, we see that {for any solution $L^{\star}$ of \eqref{eq::moment4f_SL}, $\Phi_{\mJ^{\csosc}_{h,4d}}(L^{\star})$ is a global minimizer for the sensor localization problem. } %{\color{blue} Some explanations added.}

Next, we present numerical examples for solving the source localization problem using the corresponding SOS reformulation. 

\smallskip 

\noindent \textbf{Numerical setup.} Following the experimental framework in \cite{beck2017branch}, 
we consider $p=300$ sensors in $\mathbb{R}^3$, and set the true source location to 
$\bx^{\mathrm{true}} = (1,1, -5)\in \mathbb{R}^3$. The sensor locations are generated on the horizontal 
affine plane $\mathcal H:=\{\bx\in\mathbb{R}^3:x_3=5\}$. We generate $ \bv_\ell=(v_{\ell,1},v_{\ell,2},5)^\top$ 
so that $v_{\ell,1}$ and $v_{\ell,2}$ are drawn randomly from the interval $[-10,10]$. The measured distances are 
generated as $r_\ell= \|\bx^{\mathrm{true}}{-} \bv_\ell\|{+}\varepsilon_\ell$, where $\varepsilon_\ell$ are drawn 
from a normal distribution with mean zero and standard deviation parameter $\sigma$ independently for $\ell=1,\ldots,300$. For convenience, we let $U$ be the $(3 \times 3)$ identity matrix and let the degree parameter $d=1$ (and so, the loss function is a non-homogeneous polynomial with degree $4$).  Note that the dimension condition follows by the numerical setup, and so, the loss function $f_{\rm SL}$ is a Jensen-type polynomial associated with the explicit Jensen map given in \eqref{eq:explicit_Jensen} with $\overline{\bd}=(0,0,1) \in \mathbb{R}^3$.  Thus, 
%the SOS relaxation is exact, and 
a solution can be found using the explicit form of Jensen map as described above. 

\smallskip 

\noindent{\bf Numerical Results:} In solving the SOS {reformulation \eqref{eq::sos4f_SL}}
of the source localization problem (together with its dual \eqref{eq::moment4f_SL}), 
we convert the reformulation and its dual into SDPs using CVX and solve them with the SDP solver {\sf MOSEK} \cite{mosek}.

For the purpose of visualization, we first solve an instance where we set the number of sensors $p=300$, standard deviation parameter $\sigma=0.1$, and plot out the location of the sensors, computed source location and the true source location in Figure 1. 
%The Jensen-map equations associated with an optimal dual solution of \eqref{eq::moment4f_SL} 
%have two solutions, one of which is the reflection of the other across
%the sensor plane. 
%For the instance considered here, the equations defining the Jensen map have two solutions, 
For the instance considered here, the Jensen-map equations indeed produce two solutions for the sensor localization problem (one solution $\bx=(x_1,x_2,x_3)$ and the other, reflected across the plane $\bx_{\rm reflected}=(x_1,x_2,10-x_3)$). Here, we only plot out the solution whose third coordinate is smaller than $5$.  CVX outputs the computed source location as $\bx_{\rm output}=(1.045, 0.938, -5.030)$ with the corresponding objective value $f_{\rm SL}(\bx_{\rm output})$  as {$21737.70$} compared with the true location $\bx_{\mathrm{true}} = (1,1, -5)^\top$ and the objective value  $f_{\rm SL}(\bx_{\mathrm{true}})$  as {$22071.93$}. The CPU time used is $0.095$ seconds for solving the SDP relaxation. We also report the values of the quantities: Location error, computed by $\|\bx_{\mathrm{output}}-\bx_{\mathrm{true}}\|$, Relative location error, computed by 
$\frac{\|\bx_{\mathrm{output}}-\bx_{\mathrm{true}}\|}{\|\bx_{\mathrm{true}}\|}$ and Relative objective error, computed by $\frac{|f_{\rm SL}(\bx_{\mathrm{output}})-f_{\rm SL}(\bx_{\mathrm{true}})|}{|f_{\rm SL}(\bx_{\mathrm{true}})|}$. The results for this instance are summarized in Table 1. 
%\includegraphics[scale=0.36]{graph_source_location1.jpg}

%\begin{figure}[htbp]
%    \centering
%
%    \begin{minipage}[t]{0.48\textwidth}
%        \vspace{0pt}
%        \centering
%        \includegraphics[width=\textwidth]{graph_source_location1.jpg}
%        \captionof{figure}{Source localization results.}
%        \label{fig:source-location}
%    \end{minipage}
%    \hfill
%    \begin{minipage}[t]{0.48\textwidth}
%        \vspace{0pt}
%        \centering
%        \captionof{table}{Comparison: True location vs computed source locations.}
%        \label{tab:source-location}
%
%        \begin{tabular}{lcc}
%            \toprule
%            Quantity & True value & Computed value \\
%            \midrule
%            $x_1$ & $1.000$ & $1.045$ \\
%            $x_2$ & $1.000$ & $0.938$ \\
%            $x_3$ & $-5.000$ & $-5.030$ \\
%            \midrule
%            Objective value & $2207.193$ & $2173.770$ \\
%            Location error & --- & $0.0823$ \\
%            Relative location error & --- & $1.58\%$ \\
%            Relative objective error & --- & $1.51\%$ \\
%            CPU time (seconds) & --- & $0.135$ \\
%            \bottomrule
%        \end{tabular}
%    \end{minipage}
%\end{figure}

\begin{figure}[htbp]
    \centering
    \begin{minipage}[t]{0.46\textwidth}
        \vspace{0pt}
        \centering
        \includegraphics[width=\linewidth]{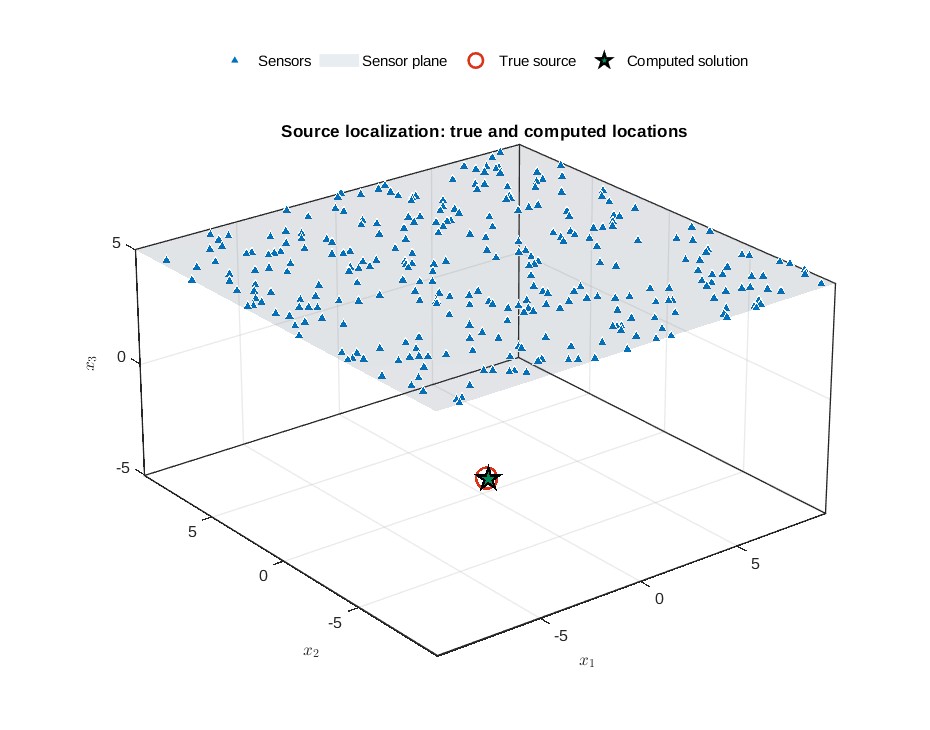}  
        \captionof{figure}{Source localization results.}
        \label{fig:source-location}
    \end{minipage}
    \hfill
    \begin{minipage}[t]{0.50\textwidth}
        \vspace{0pt}
        \centering
        \captionof{table}{Comparison: true vs.\ computed source locations.}
        \label{tab:source-location}

        \small                              
        \setlength{\tabcolsep}{3pt}         

        \resizebox{\linewidth}{!}{%         
        \begin{tabular}{lcc}
            \toprule
            Quantity & True value & Computed value \\
            \midrule
            $x_1$ & $1.000$ & $1.045$ \\
            $x_2$ & $1.000$ & $0.938$ \\
            $x_3$ & $-5.000$ & $-5.030$ \\
            \midrule
            Objective value          & {$22071.93$} & {$21737.70$} \\
            Location error           & ---        & $0.0823$   \\
            Relative location error  & ---        & $1.58\%$   \\
            Relative objective error & ---        & $1.51\%$   \\
            CPU time (seconds)             & ---        & $0.095$    \\
            \bottomrule
        \end{tabular}%
        }
    \end{minipage}
\end{figure}

We also vary the number of sensors $p\in\{200,300,500\}$ and the noise
standard deviation $\sigma\in\{1,0.1,0.01\}$. For each pair $(p,\sigma)$,
we generate 50 independent instances. Table~2 reports, for each pair,
the sample averages of the computed source locations, relative location
errors, relative objective errors, and CPU times (in seconds) across the 50 trials. %The average source location is computed componentwise.

% We also vary the number of sensors $p$ in $\{200,300,500\}$ and the standard deviation parameter $\sigma$ in $\{1,0.1,0.01\}$. For each pair $(p,\sigma)$, we generate 50 independent instances. In Table 2, we report the sample averages over these 50 trials of the computed source location, relative location error, relative objective error, and CPU time (in seconds). In particular, the reported average source location is the componentwise arithmetic mean of the 50 computed locations.

%we run 50 trials and output the computed solutions, relative location error, 
 %           relative objective error and CPU time per instance (measured in seconds), on average. The results are summarized in Table 2. 
% Requires: \usepackage{booktabs}
%          \usepackage{graphicx}

\begin{table}[htbp]
\centering
\caption{Performance with different $(p,\sigma)$}
\label{tab:jensen-map-results}
\scriptsize
\setlength{\tabcolsep}{3pt}
\resizebox{\textwidth}{!}{%
\begin{tabular}{@{}cccccc@{}}
\toprule
\# of Sensors
& Noise std.
& \shortstack{Average computed location\\$(x_1,x_2,x_3)$}
& \shortstack{Relative location\\error (\%)}
& \shortstack{Relative objective\\error (\%)}
& \shortstack{CPU time\\(sec.)} \\
\midrule
200 & 1.000 & $(0.943,\ 1.027,\ -5.045)$ &  4.926 & 1.941 & 0.106 \\
200 & 0.100 & $(1.004,\ 0.996,\ -4.999)$ &  0.449 & 1.639 & 0.100 \\
200 & 0.010 & $(1.000,\ 1.000,\ -5.000)$ &  0.044 & 1.432 & 0.102 \\
300 & 1.000 & $(0.965,\ 1.040,\ -5.049)$ &  3.865 & 1.263 & 0.101 \\
300 & 0.100 & $(1.000,\ 0.999,\ -5.000)$ &  0.365 & 1.109 & 0.102 \\
300 & 0.010 & $(1.000,\ 1.000,\ -5.000)$ &  0.037 & 1.094 & 0.100 \\
500 & 1.000 & $(1.017,\ 0.993,\ -5.045)$ &  3.200 & 0.898 & 0.105 \\
500 & 0.100 & $(1.001,\ 1.001,\ -5.000)$ &  0.266 & 0.553 & 0.103 \\
500 & 0.010 & $(1.000,\ 1.000,\ -5.000)$ &  0.029 & 0.665 & 0.101 \\
\bottomrule
\end{tabular}%
}
\end{table}
%{\color{blue} We can remove the column of the average objective value.}

\subsection{Application 2: Semi-infinite optimization with Jensen-type polynomials}
Under Assumption~\ref{assump::0}, consider the nonconvex semi-infinite polynomial optimization problem 
\begin{equation}\label{eq::sipop}
\left\{
\begin{aligned}
f^{\star}:=\inf_{\bx\in\RR^n}&\ f_0(\bx)\\
\quad\text{s.t. }&\ f_1(\bx)\le 0,\ldots,f_J(\bx)\le 0,\\
&\ f(\bx, \by)\le 0,\ \forall \by\in \Y,\\
\end{aligned}\tag{SIPOP}
\right.
\end{equation}
where $\Y$ is a semialgebraic set defined as in \eqref{eq::Y} with the functions $g_1,\ldots,g_r$. We assume that $f^{\star}$ is a finite number. The class of problems considered here is generally challenging to solve. It includes, for example, SOS-convex polynomial optimization problems with semi-infinite constraint (in particular, semi-infinite convex quadratic optimization problems) and ENPC polynomial optimization problems with semi-infinite constraints.    
In this section,
we aim to establish two hierarchies of SDP relaxations for \eqref{eq::sipop}
%by Theorem \ref{th::Pastz} (i) and (ii), 
which can provide two sequences of convergent upper and lower bounds of $f^{\star}$. 
Recall the number $d_0$ in \eqref{eq::degree}.

%{\color{blue} I imagine, we need to know how $\Y$ looks like (e.g. defined by $g_1...g_r$) in the lower bound approximation? }{\color{red} I add Remark 7.1}

%First note that 
%\begin{proposition}\label{prop::R_+}
%    Under Assumption \ref{assump::0}, we have $\Phi_{\mJ}(\X)\subset\X$ and 
%    there exists a minimizing sequence
%    of \eqref{eq::sipop} in $\Phi_{\mJ}(\X)$. In particular, if $f^{\star}$ is attainable,
%    an optimal solution exists in $\Phi_{\mJ}(\X)$.
%\end{proposition}
%\begin{proof}
%    By the condition (C), it holds
%    \[
%    \begin{aligned}
%    &f_j(\Phi_{\mJ}(L_{\bu}))\le L_{\bu}(f_j)=f_j(\bu),\ \forall j\in[1:J],\ \
%     f(\Phi_{\mJ}(L_{\bu}),\bv)\le L_{\bu}(f(\bx,\bv))=f(\bu,\bv),\ \forall \bv\in\Y, \\
%    &f_0(\Phi_{\mJ}(L_{\bu}))\le L_{\bu}(f_0)=f_0(\bu).
%    \end{aligned}
%    \]
%    for any $\bu\in\X$. The conclusion follows.
%\end{proof}
%
\vskip 5pt
\noindent{\bf SDP relaxation for upper bounds of $f^{\star}$.}  
For $k\ge d_0$, consider the SDP relaxation problem of \eqref{eq::sipop}
%\begin{equation}\label{eq::uppersdp}
%\left\{
%\begin{aligned}
%    \rho^{\psdp}_k:=\sup_{\rho,\lambda_j,\bz,s}&\ \rho\\
%    \text{s.t. }&\ f_0(\bx)-\rho=s(\bx)+\sum_{j=1}^J\lambda_jf_j(\bx)+L_{\bz}(f(\bx,\by)) \\
%    &\ \lambda_j\ge 0, j\in[1:J],\ \bz\in\RR^{\vert\mathbb{N}^{\ell}_{2k}\vert},\ s(\bx)\in\Sigma[\bx]_m, \\ 
%    &\ \bM_k(\bz)\succeq 0,\ \bM_{k-d_w}(g_w\bz)\succeq 0,\ w\in[1:r].
%\end{aligned}\right.
%\end{equation}
\begin{equation}\label{eq::uppersdp}
    \rho^{\psdp}_k:=\sup_{\rho\in\RR}\ \rho\ \text{ s.t. }\ f_0(\bx)-\rho\in \mC_k(\bff),
\end{equation}
where $\mC_k(\bff)$ is given as in \eqref{eq::Qk}. 
%Denote by $[\bx]_{\frac{m}{2}}$ the column vector of monomials in $\RR[\bx]$ of degree up to $m/2$, i.e.,
%\[
%[\bx]_{\frac{m}{2}}:=[1,\ x_1,\ x_2,\ \cdots,\ x_n,\ x_1^2,\ x_1x_2,\ \cdots,\ x_n^{m/2}]^{\intercal}
%\]
%\begin{assumption}\label{assump::R}
%    We have $f_1=$
%\end{assumption}
%
Consider the dual cones of $\mC_k(\bff)$
\begin{equation}\label{eq::QK*}
    (\mC_k(\bff))^*:=\left\{L\in(\RR[\bx]_m)^*\mid L\in(\Sigma[\bx]_m)^*, 
    \ \ L(f_j)\le 0,\ j\in[1:J],\ -L(f(\bx,\by))\in \Qk(\bg)\right\},
\end{equation}
where $\bg=\{g_1,\ldots,g_r\}$ and $\Qk(\bg)$ is the truncated quadratic module generated by $\bg$. 
The Lagrange dual  of \eqref{eq::uppersdp} reads as the following SDP problem
%\begin{equation}\label{eq::uppersdp_dual}
%\left\{
%\begin{aligned}
%    \rho^{\dsdp}_k:=\inf_{\bw\in\RR^{|\N^n_m|}}&\ L_{\bw}(f_0)\\
%    \text{s.t. }&\ L_{\bw}(1)=1,\ \bM_{m/2}(\bw)\succeq 0, \\
%    &\ L_{\bw}(f_j)\ge 0,\ j\in[1:J],\ L_{\bw}(f(\bx,\by))\in Q_k(\bg).
%\end{aligned}\right.
%\end{equation}
\begin{equation}\label{eq::uppersdp_dual}
    \rho^{\dsdp}_k:=\inf_{L\in(\RR[\bx]_m)^*}\ L(f_0)\ \text{ s.t. }\ L\in(\mC_k(\bff))^*,\  L(1)=1. 
\end{equation}
Before we proceed, we note that, by construction, for all $k \in \mathbb{N}$, 
\[
\rho_{k+1}^{\psdp}\le \rho_k^{\psdp},
\mbox{ and }
\rho_{k+1}^{\dsdp}\le \rho_k^{\dsdp}. 
\]

%We write $\be=(\be_1,\ldots,\be_n)$ for the standard basis of $\RR^n$.  
%For any feasible point $\bw=(w_{\ba})_{\ba\in\N^n_m}$ of \eqref{eq::uppersdp_dual}, 
%recall that 
%\[
%\bx_{\bw}:=\left(\sqrt[m]{w_{m\be_1}},\ldots,\sqrt[m]{w_{m\be_n}}\right)\in\RR^n,
%\]

\begin{lemma}\label{lem::feasibility}
    Suppose that Assumption~\ref{assump::0} holds and let $\Phi_{\mJ}$ be the associated Jensen map given in Assumption~\ref{assump::0}. Then,  $\Phi_{\mJ}(L)\in\X$
    for any feasible point $L$ of \eqref{eq::uppersdp_dual}.
\end{lemma}
\begin{proof}
    For any $L\in (\mC_k(\bff))^*$ with $L(1)=1$, we have $L\in\mathcal{B}[\bx]_m$. Thus, 
    by the condition (C), we have 
    \[
    f_j(\Phi_{\mJ}(L))\le L(f_j)\le 0,\  j\in[1:J],\  \ f(\Phi_{\mJ}(L),\by)\le L(f(\bx,\by))\le 0,\ \by\in\Y,
    \]
    which implies $\Phi_{\mJ}(L)\in\X$.
\end{proof}

\begin{theorem}\label{th::convergence_sym}
Suppose that Assumption \ref{assump::0} and the Slater condition for $\bff$ hold. For $k\ge d_0$,
let $L^{(k,\star)}$ be a minimizer of \eqref{eq::uppersdp_dual} 
with order $k$.
%and recall that
%\[
%\bx_{\bw^{(k,\star)}}:=\left(\sqrt[m]{w^{(k,\star)}_{m\be_1}},\ldots,\sqrt[m]{w^{(k,\star)}_{m\be_n}}\right)\in\RR^n.
%\]
Then, $\rho^{\dsdp}_k\ge \rho^{\psdp}_k \ge f^{\star}$ for any $k\ge d_0$.
%\begin{enumerate}
%    \item[\upshape (i)] $\rho^{\dsdp}_k\ge \rho^{\psdp}_k \ge f^{\star}$ for any $k\ge d_0$;
%\end{enumerate}
If, moreover, the quadratic module $\mQ(\bg)$ is Archimedean, then 
\begin{enumerate}
    \item[\upshape (i)] $\rho^{\dsdp}_k\searrow f^{\star}$ and $\rho^{\psdp}_k\searrow f^{\star}$ as $k\rightarrow\infty$;
    \item[\upshape (ii)] Suppose further that the feasible region {$\mathcal{X}$ of \eqref{eq::sipop} is compact}. Then, $\{\Phi_{\mJ}(L^{(k,\star)})\}_k$ is a bounded sequence, and for any convergent subsequence $\{\Phi_{\mJ}(L^{(k_i,\star)})\}_i$ 
of $\{\Phi_{\mJ}(L^{(k,\star)})\}_k$, the limit
$\lim_{i\rightarrow\infty}\Phi_{\mJ}(L^{(k_i,\star)})$ is a global minimizer of \eqref{eq::sipop}.
Consequently, if the set of optimal solutions of \eqref{eq::sipop} 
is a singleton, then $\lim_{k\rightarrow\infty} \Phi_{\mJ}(L^{(k,\star)})$ is the unique 
optimal solution of \eqref{eq::sipop}.
\end{enumerate}
\end{theorem}
\begin{proof}
As $f_0(\bx)-f^{\star}\ge 0$ on $\X$ and $\mC(\bff)\subset\mC_k(\bff)$ for each $k\ge d_0$, 
by Theorem \ref{th::Pastz} one has $\rho^{\psdp}_k\ge f^{\star}$. Noting that 
\eqref{eq::uppersdp_dual} is a minimization problem whose dual problem is \eqref{eq::uppersdp},  
this together with the standard weak 
duality for SDP implies that $\rho^{\dsdp}_k\ge \rho^{\psdp}_k\ge f^{\star}$. 
Now assume that $\mQ(\bg)$ is Archimedean.  
    
\noindent (i) Fixing an arbitrary $\varepsilon\in(0, 1)$, there exists a point $\bu^{(\varepsilon)}\in\X$
such that $f_0(\bu^{(\varepsilon)})<f^{\star}+\varepsilon$.
Let $\bar{\bu}\in\X$ be a Slater point and 
$L^{(\varepsilon)}=(1-\varepsilon)L_{\bu^{(\varepsilon)}}+\varepsilon L_{\bar{\bu}}$,
where $L_{\bu}$ denotes the point-evaluation functional
$L_{\bu}(q):=q(\bu)$ for all $q \in \RR[\bx]_m$.  
Then, it is clear that $L^{(\varepsilon)}\in\mathcal{B}[\bx]_m$.
% $(\Sigma[\bx]_m)^*$ and $L^{(\varepsilon)}(1)=1$. 
For $j\in[1:J]$,
\[
L^{(\varepsilon)}(f_j)=(1-\varepsilon)f_j(\bu^{(\varepsilon)})+\varepsilon f_j(\bar{\bu})<0.
\]
For any $\by\in\Y$, 
\[
L^{(\varepsilon)}(f(\bx,\by))=(1-\varepsilon)f(\bu^{(\varepsilon)},\by)+\varepsilon f(\bar{\bu},\by)<0.
\]
As $\mQ(\bg)$ is Archimedean, by Putinar's Positivstellensatz, 
there exists $k_{\varepsilon}\in\N$ such that $-L^{(\varepsilon)}(f(\bx,\by))\in \mQ_k(\bg)$ 
for all $k\ge k_{\varepsilon}$. Hence, for all $k\ge k_{\varepsilon}$, $L^{(\varepsilon)}$ 
is feasible to \eqref{eq::uppersdp_dual} with order $k$ and thus
\[
\rho^{\dsdp}_k\le L^{(\varepsilon)}(f_0)=
(1-\varepsilon)f_0(\bu^{(\varepsilon)})+\varepsilon f_0(\bar{\bu})
\le (1-\varepsilon)(f^{\star}+\varepsilon)+\varepsilon f_0(\bar{\bu}).
\]
Since $\varepsilon>0$ is arbitrary, letting $\varepsilon\to 0$, 
we have $\lim_{k\to\infty}\rho^{\dsdp}_k\le f^{\star}$.
Thus, by the weak duality,
$\rho^{\dsdp}_k\searrow f^{\star}$ and $\rho^{\psdp}_k\searrow f^{\star}$ as $k\rightarrow\infty$;

\noindent (ii) { Let us further assume that $\mathcal{X}$ is compact.} Then, the preceding lemma implies that $\{\Phi_{\mJ}(L^{(k,\star)})\}_k \subseteq \mathcal{X}$ is a bounded sequence. Let $\bu^{\star}:=\lim_{i\rightarrow\infty}\Phi_{\mJ}(L^{(k_i,\star)})$. 
By Lemma \ref{lem::feasibility}, we have $\bu^{\star}\in\X$. By the condition (C),
\[
f_0(\bu^{\star})=\lim_{i\to\infty}f_0(\Phi_{\mJ}(L^{(k_i,\star)}))\le \lim_{i\to\infty}L^{(k_i,\star)}(f_0)=\lim_{i\to\infty}\rho^{\dsdp}_{k_i}=f^{\star}.
\]
Thus, $\bu^{\star}$ is a global minimizer of \eqref{eq::sipop}.
\end{proof}

\begin{theorem}\label{th::convergence_finite}
      Suppose that Assumption \ref{assump::0} and the Slater condition for $\bff$ hold.
      Let $k\ge d_0$, and let $\bz^{(k)}$ be the vector in the representation \eqref{eq::Qk} of 
      $f_0(\bx)-\rho^{\dsdp}_k$ and 
 %     $(\rho^{(k)}, \lambda^{(k)}_j, j\in[1:J], \bz^{(k)},s^{(k)})$ the vector in the representation of $f_0\in Q_k(\bff)$; and 
      $L^{(k,\star)}$ be an optimal solution of \eqref{eq::uppersdp_dual}. 
      If  the following condition holds
      \begin{equation}\label{eq::FEC}
       {\exists\ t\quad\text{s.t. }  d_0\le t\le k\quad \text{and}\quad  \rank\ \bM_t({\bz^{(k)}})=\rank\ \bM_{t-d_{\bg}}({\bz^{(k)}})},
%          {\color{blue} r_0:=} \, \rank\ \bM_k(\bz^{(k)})=\rank\ \bM_{k-d_{\bg}}(\bz^{(k)}),
      \end{equation}
then $\rho^{\psdp}_k=f^{\star}$.  Moreover, if we further assume that  strong duality holds
between \eqref{eq::uppersdp} and \eqref{eq::uppersdp_dual}, then $\Phi_{\mJ}(L^{(k,\star)})$ is an 
optimal solution of \eqref{eq::sipop}.
      \end{theorem}
\begin{proof}
By Theorem \ref{th::Qk}, 
the condition \eqref{eq::FEC} implies that $\rho^{\psdp}_k\le f^{\star}$. 
By Theorem \ref{th::convergence_sym}, $\rho^{\psdp}_k\ge f^{\star}$, and so,  $\rho_k^{\psdp}=f^{\star}$.  
    If the strong duality holds for \eqref{eq::uppersdp} and \eqref{eq::uppersdp_dual}, 
    then by the condition (C),
    \[
    f_0(\Phi_{\mJ}(L^{(k,\star)}))\le L^{(k,\star)}(f_0)=\rho^{\dsdp}_k=f^{\star}.
    \]
    As $\Phi_{\mJ}(L^{(k,\star)})\in\X$ by Lemma \ref{lem::feasibility}, it implies that 
    $\Phi_{\mJ}(L^{(k,\star)})$ is a global minimizer of \eqref{eq::sipop}.
\end{proof}

%{\color{blue} Add a corollary and Remark. Compare with Nie and also the paper with Shenglong.}

%\begin{theorem}
%Suppose that Assumptions \ref{assump::sipop} and \ref{assump::f0} hold, and the strong duality holds
%for \eqref{eq::uppersdp} and \eqref{eq::uppersdp_dual}. 
%Let $\bw^{(k,\star)}=(w^{(k,\star)}_{\ba})_{\ba\in\N^n_m}$ be a minimizer of \eqref{eq::uppersdp_dual} 
%with order $k\ge d_0$ and let
%\[
%\bx_{\bw^{(k,\star)}}:=\left(\sqrt[m]{w^{(k,\star)}_{m\be_1}},\ldots,\sqrt[m]{w^{(k,\star)}_{m\be_n}}\right)\in\RR^n.
%\]
%\begin{enumerate}
%    \item[\upshape (i)] For any convergent subsequence $\{\bw^{(k_i,\star)}\}_i$ 
%of $\{\bw^{(k,\star)}\}_k$,
%$\lim_{i\rightarrow\infty}\bx_{\bw^{(k_i,\star)}}$ is a global minimizer of \eqref{eq::sipop}.
%Consequently, if the set of optimal solutions of \eqref{eq::sipop} 
%is a singleton, then $\lim_{k\rightarrow\infty}\bx_{\bw^{(k,\star)}}$ is the unique global minimizer;
%\item[\upshape (ii)] If the condition \eqref{eq::FEC} holds for some $k\ge d_0$, then 
%$\bx_{\bw^{(k,\star)}}$ is a global minimizer of \eqref{eq::sipop}.
%\end{enumerate}
%\end{theorem}
%\begin{proof}
%    
%\end{proof}

\noindent{\bf SDP relaxations for lower bounds of $f^{\star}$.} 
Let $\nu$ be an arbitrary, fixed, finite nonnegative Borel measure supported on $\Y$.
For %a fixed $\varepsilon>0$ and 
each $k\in\N$, 
%by restricting the degree of $s(\bx)$ in \eqref{eq::lower},
consider the SDP relaxation problem
\begin{equation}\label{eq::lowersdp}
\gamma^{\psdp}_k:=\sup_{\gamma\in\RR}\ \gamma\ \text{ s.t. }\ f_0(\bx)-\gamma\in \mC_k(\bff, \nu),
\end{equation}
where $\mC_k(\bff, \nu)$ is given as in \eqref{eq::Qkv}. 
Consider the dual cones of $\mC_k(\bff,\nu)$
\begin{equation}\label{eq::Qkv*}
(\mC_k(\bff,\nu))^*:=
\left\{
\begin{array}{c|c}
 \displaystyle  L\in(\RR[\bx]_m)^*  &  
 \displaystyle  \begin{aligned}
    & L\in(\Sigma[\bx]_m)^*, \ \ L(f_j)\le 0,\ j\in[1:J],\\
    & \int_{\Y} \sigma(\by)L(f(\bx,\by))\ud\nu(\by)\le 0, \forall \sigma(\by)\in\Sigma[\by]_{2k}
   \end{aligned}
\end{array}
\right\}.
\end{equation}
The dual SDP relaxation problem of \eqref{eq::lowersdp} reads as
\begin{equation}\label{eq::lowersdp_dual}
\gamma^{\dsdp}_k:=\inf_{L\in(\RR[\bx]_m)^*}\ 
L(f_0)\ \text{ s.t. }\ L\in(\mC_k(\bff,\nu))^*,\  L(1)=1. 
\end{equation}
By construction, it can be directly verified that for each $k \in \mathbb{N}$, 
\[
\gamma_k^{\psdp}\le \gamma_{k+1}^{\psdp},
\mbox{ and }
\gamma_k^{\dsdp}\le \gamma_{k+1}^{\dsdp}.
\]
%Clearly, $\gamma_k^{\psdp}\le \gamma_k^{\dsdp}\le f^{\star}$ for each $k\in\N$. In fact,

%\begin{equation}\label{eq::lowersdp}
%\gamma^{\psdp}_k(\varepsilon):=\sup_{\gamma\in\RR}\ \gamma\ 
%    \text{ s.t. }\ f_0(\bx)-\gamma+\varepsilon\Theta_m(\bx)\in \mC_k(\bff, \nu), 
%%\left\{
%%\begin{aligned}
%%    \gamma^{\star}_k(\varepsilon):=\sup_{\gamma,\lambda_j,\sigma,s}&\ \gamma\\
%%    \text{s.t. }&\ f_0(\bx)-\gamma+\varepsilon\Theta_m(\bx)=s(\bx)+\sum_{j=1}^J\lambda_jf_j(\bx)+\int_{\Y}\sigma(\by)f(\bx,\by)\ud\nu(\by), \\
%%    &\ \lambda_j\ge 0, j\in[1:J], \sigma(\by)\in\Sigma[\by]_{2k}, s(\bx)\in\Sigma[\bx]_m.
%%\end{aligned}\right.
%\end{equation}
%and its dual SDP relaxation problem
%\begin{equation}\label{eq::lowersdp_dual}
%\gamma^{\dsdp}_k(\varepsilon):=\inf_{\bw\in\RR^{|\N^n_m|}}\ 
%L_{\bw}(f_0+\varepsilon \Theta_m)\ \text{ s.t. }\ \bw\in(\mC_k(\bff,\nu))^*,\  w_0=1. 
%\end{equation}
%

\begin{lemma}\label{lem::feasibility2}
    Suppose that Assumption \ref{assump::0} holds and $\supp(\nu)=\Y$. 
    For $k\in\N$, let $L^{(k)}$ be a feasible point of \eqref{eq::lowersdp_dual}. 
    If $\bu=\lim_{k\to\infty}\Phi_{\mJ}(L^{(k)})$, then $\bu\in\X$.
\end{lemma}
\begin{proof}
    For each $k\in\N$, by the condition (C), 
    \begin{equation}\label{eq::sosint}
    \begin{aligned}
    &f_j(\Phi_{\mJ}(L^{(k)}))\le L^{(k)}(f_j(\bx))\le 0,\ \forall j\in[1:J],\\
    &  \int_{\Y} \sigma(\by)f(\Phi_{\mJ}(L^{(k)}),\by)\ud\nu(\by)\le
    \int_{\Y} \sigma(\by)L^{(k)}(f(\bx,\by))\ud\nu(\by)\le 0,\ \ \forall \sigma(\by)\in\Sigma[\by]_{2k}.
    \end{aligned}
    \end{equation}
    Fix an arbitrary $\sigma\in\Sigma[\by]$. Since $\sigma$ has finite
degree, there exists $k_\sigma\in\mathbb{N}$ such that
$\sigma\in\Sigma[\by]_{2k}$ for every $k\ge k_\sigma$. Hence, for all
$k\ge k_\sigma$, feasibility of $L^{(k)}$ gives
\[
f_j(\Phi_{\mJ}(L^{(k)}))\le 0, \ \forall j\in[1:J] \mbox{ and } \int_{\Y}\sigma(\by)
f(\Phi_{\mJ}(L^{(k)}),\by)\,\ud\nu(\by)
\le 0.
\]
Passing to the limit and using the continuity of $f$ and $f_j$, and $\bu=\lim_{k\to\infty}\Phi_{\mJ}(L^{(k)})$, it follows that 
\[
f_j(\bu) \le 0 , \ \forall j\in[1:J] \mbox{ and } \int_{\Y}\sigma(\by)f(\bu,\by)\,\ud\nu(\by)\le0.
\]
Since the latter holds for every $\sigma\in\Sigma[\by]$,  \cite[Theorem 3.2]{LasserreOuter} shows that $f(\bu,\by)\le 0$ on $\Y$ and 
    thus $\bu\in\X$.  
    % Since $\bu=\lim_{k\to\infty}\Phi_{\mJ}(L^{(k)})$, by the continuity of $f_j$, and the condition (C),
    % \[
    % f_j(\bu)=\lim_{k\to\infty}f_j(\Phi_{\mJ}(L^{(k)}))\le 0, \ \forall j\in[1:J], 
    % \]
    % and 
    % \[
    % \int_{\Y} \sigma(\by)f(\bu,\by)\ud\nu(\by)=\lim_{k\to\infty}\int_{\Y} \sigma(\by)f(\Phi_{\mJ}(L^{(k)}),\by)\ud\nu(\by)\le 0,\  \forall \sigma(\by)\in\Sigma[\by].
    % \]
    % By \cite[Theorem 3.2]{LasserreOuter}, $f(\bu,\by)\le 0$ on $\Y$ and 
    % thus $\bu\in\X$.
\end{proof}

\begin{theorem}\label{th::lower}
{For each $k\in\N$, $\gamma_k^{\psdp}\le \gamma_k^{\dsdp}\le f^{\star}$. Moreover,
under Assumption \ref{assump::0}, suppose that 
    the Slater condition holds for $\bff$ and $\supp(\nu)=\Y$, then the following hold.}
%     Suppose that Assumption \ref{assump::0} and the Slater condition for $\bff$ hold, $\supp(\nu)=\Y$ and 
     \begin{enumerate}
%     \item[\upshape (i)]  For each $k\in\N$, $\gamma_k^{\psdp}\le \gamma_k^{\dsdp}\le f^{\star}$. 
         \item[\upshape (i)]  If $\sum_{i=1}^n x_i^m-R\in\bff$ for some $R>0$ and $\Theta_m \in \mJ$, then  $\gamma_k^{\psdp}\nearrow f^{\star}$ and 
         $\gamma_k^{\dsdp}\nearrow f^{\star}$ as $k\to\infty$;
     \item[\upshape (ii)] 
     Let $L^{(k,\star)}$ be an optimal solution of \eqref{eq::lowersdp_dual}. For any convergent subsequence $\{\Phi_{\mJ}(L^{(k_i,\star)})\}_i$ 
of $\{\Phi_{\mJ}(L^{(k,\star)})\}_k$, the limit
$\bu^{\star}:=\lim_{i\rightarrow\infty}\Phi_{\mJ}(L^{(k_i,\star)})$ is a global minimizer of \eqref{eq::sipop}.
     \end{enumerate}
\end{theorem}
\begin{proof}
Let $\bu$ be an arbitrary feasible point of \eqref{eq::sipop}, and let $L_{\bu}$ be the point-evaluation functional at $\bu$, namely,
$L_{\bu}(q):=q(\bu)$ for all $q \in \RR[\bx]$. 
%let $L_{\bu}$ be the linear functional that evaluating at the point $\bu$.  
Then $L_{\bu}$ is feasible to \eqref{eq::lowersdp_dual}. Hence, 
    $\gamma_k^{\dsdp}\le f^{\star}$
    and $\gamma_k^{\psdp}\le \gamma_k^{\dsdp}\le f^{\star}$ holds due to the weak duality. 
    
(i) Note that $f_0-f^\star \ge 0$ on $\mX$. This part follows by applying Corollary \ref{cor::Qkv} (ii) and the weak duality.

(ii) By Lemma \ref{lem::feasibility2}, we have $\bu^{\star}\in\X$. 
By the condition (C) in the definition of Jensen-type polynomials,
\[
f^{\star}\le f_0(\bu^{\star})=\lim_{i\to\infty}f_0(\Phi_{\mJ}(L^{(k_i,\star)}))
\le \lim_{i\to\infty}L^{(k_i,\star)}(f_0)=\lim_{i\to\infty}\gamma_{k_i}^{\dsdp} \le f^{\star}.
\]
Therefore, $f_0(\bu^{\star})=f^{\star}$ and $\bu^{\star}$ is an optimal solution of \eqref{eq::sipop}.
\end{proof}

\begin{remark}
To implement the SDP relaxations \eqref{eq::lowersdp} and
\eqref{eq::lowersdp_dual}, one needs the moments
$\int_{\Y}\by^{\bb}\,\ud\nu(\by)$ for all $\bb\in\N^\ny$. 
For several common choices of $\Y$ and $\nu$, these moments are
available explicitly. Examples include:
\begin{itemize}
    \item a hyperrectangle in $\RR^\ny$ with Lebesgue measure;
    \item the unit sphere in $\RR^\ny$ with
    $(\ny-1)$-dimensional surface measure;
    \item the unit ball in $\RR^\ny$ with Lebesgue measure;
    \item a simplex in $\RR^\ny$ with Lebesgue measure.
\end{itemize}
Thus, for these choices, the moment data required by the SDP
relaxations can be computed in closed form.
\end{remark}

\begin{remark}[Comparison to the literature]
Hu, Klep, and Nie \cite{HKN2024} establish Positivstellensatz and
moment results for sets described by universally quantified
polynomial inequalities. Let $\nu$ be a finite nonnegative Borel measure
with $\operatorname{supp}(\nu)=\Y$. Since $\Y$ is compact, $\nu$ has
compact support. Applying their SOS formulation to
\eqref{eq::sipop}, with the inequalities written in the form
$f_j (\bx)\leq 0$ and $f(\bx,\by)\leq 0$, gives the following SDP 
relaxation with order $k$:
\begin{equation}\label{eq::hkn}
\left\{
\begin{aligned}
    \sup_{\gamma\in\RR}\quad & \gamma\\
    \text{\rm s.t.}\quad
    & f_0(\bx)-\gamma
      =s_0(\bx)-\sum_{j=1}^J s_j(\bx)f_j(\bx)
       -\int_{\Y}\sigma(\bx,\by)f(\bx,\by)\,\ud\nu(\by),\\
    & s_0,s_1,\ldots,s_J\in\Sigma[\bx],
      \qquad \sigma\in\Sigma[\bx,\by],\\
    & \deg(s_0)\leq 2k,\qquad
      \deg(s_jf_j)\leq 2k,\quad j=1,\ldots,J,\\
    & \deg(\sigma f)\leq 2k.
\end{aligned}
\right.
\end{equation}
Under suitable assumptions,  both this SOS
relaxation and its moment dual provide lower bounds on the true optimal value $f^\star$,
and the bounds converge to $f^\star$ as $k\to\infty$.

The lower-bound hierarchy \eqref{eq::lowersdp} shares a key feature with
the hierarchy in \cite{HKN2024}: it yields lower bounds on $f^\star$
that converge to $f^\star$.  On the other hand, the SDP hierarchy \eqref{eq::lowersdp} is simpler and, in general, is cheaper to solve.  In particular, the multipliers for the finite
constraints are nonnegative scalars (instead of SOS polynomials  with increasing degrees), and the multiplier for the 
universally quantified constraint depends only on $\by$ (rather than on both variables $\bx$ and $\by$). Here, under the Jensen-class structural 
assumptions, the convergence result established above shows that a simpler convergent SDP hierarchy can be utilized to solve the original problem.

We also establish an SDP hierarchy that provides upper bounds on
$f^\star$, which was not considered in~\cite{HKN2024}. Together, the upper and lower bounds give an optimality gap
that can be used to assess solution quality and as a stopping criterion,
including when the rank condition for finite convergence is
not satisfied.
\end{remark}

\vskip 5pt

\noindent{\bf Numerical Examples for Semi-infinite Optimization with Jensen-type Polynomials.} 
Below, we provide two numerical examples illustrating how to utilize the lower- and upper-bound hierarchies to solve semi-infinite programming problems involving Jensen-type polynomials. 

We implement the relaxations \ref{eq::uppersdp}, \eqref{eq::uppersdp_dual},
\eqref{eq::lowersdp} and \eqref{eq::lowersdp_dual} using {\sf Yalmip} \cite{YALMIP}
and solve the resulting SDP problems by {\sf MOSEK} \cite{mosek}. 

\begin{example}\label{example:7.1}
Consider the following semi-infinite polynomial optimization problem:
\begin{equation}\label{eq::enpc}
\left\{
\begin{aligned}
\min_{\bx\in\RR^2}\quad
& f_0(\bx)=2x_1^6+2x_2^6-3x_1^2x_2^2\\
\text{s.t.}\quad
& f_1(\bx)=x_1^6+x_2^6-2\leq 0,\\
& f(\bx,\by)
=1-(1+y_1^2)x_1^6-(1+y_2^2)x_2^6\\
&\qquad
-\bigl(2+(y_1-0.5)^2+(y_2+0.5)^2\bigr)x_1^3x_2^3
\leq 0,\qquad \forall\,\by\in\Y,
\end{aligned}
\right.
\end{equation}
where $\Y=[-1,1]^2$. The feasible set is shown in Figure~\ref{fig:enpc_fs}.

We first note that Assumption~\ref{assump::0} holds for \eqref{eq::enpc} with
$\mJ=\mJ^{\enpc}_{\subnm{2}{6}}$ (ENPC polynomials with $2$ variables and degree at most $6$), and the Slater condition holds for the
constraint system.

We then solve the SDP relaxations \eqref{eq::uppersdp} and
\eqref{eq::uppersdp_dual} to
compute upper bounds on the true optimal value $f^\star$, using the relaxation order
$k=1$. Both yield the value $-0.2022$:
 $\rho_1^{\psdp}=\rho_1^{\dsdp}=-0.2022$. The computed moment matrix satisfies the rank condition \eqref{eq::FEC} with $t=1$ (to numerical precision).  Hence,
Theorem~\ref{th::convergence_finite} certifies finite convergence holds numerically, yields a numerically certified candidate global minimizer 
$\bu^\star=(0.7857,\,0.7857)$, subject to machine precision. 
Figure~\ref{fig:enpc_fs} shows the feasible set, the contour curves of
$f_0$ (including the level $f_0(\bx)=-0.2022$), and the computed
minimizer (the red point).

Alternatively, for this example, we solve \eqref{eq::lowersdp} and \eqref{eq::lowersdp_dual} to
compute lower bounds on $f^\star$, where we take $\nu$ to be Lebesgue measure on $\Y$. By Theorem~\ref{th::lower}, the
sequences $\{\gamma_k^{\psdp}\}$ and
$\{\gamma_k^{\dsdp}\}$ converge asymptotically to the true optimal value $f^\star$ as $k \rightarrow \infty$. For
$k=1,\ldots,10$, we report the bounds $\gamma_k^{\dsdp}$, the points
\[
\bu^{(k)}
:=
\Phi_{\mJ^{\enpc}_{\subnm{2}{6}}}\bigl(L^{(k,\star)}\bigr),
\]
where $L^{(k,\star)}$ is an optimal solution of
\eqref{eq::lowersdp_dual}, and the corresponding CPU times (measured in seconds) in
Table~\ref{tab:enpc_results}. At $k=10$, $\gamma_{10}^{\dsdp}=-0.2067$, while the computed point is 
$\bu^{(10)}=(0.7823,\,0.7823)$. Together with the upper bound computed earlier, this shows that the approximate solution’s objective value is within $5\times10^{-3}$ of the true optimal value.

\begin{figure}[htbp]
    \centering
    \begin{minipage}[c]{0.4\textwidth}
        \centering
        \includegraphics[width=\linewidth]{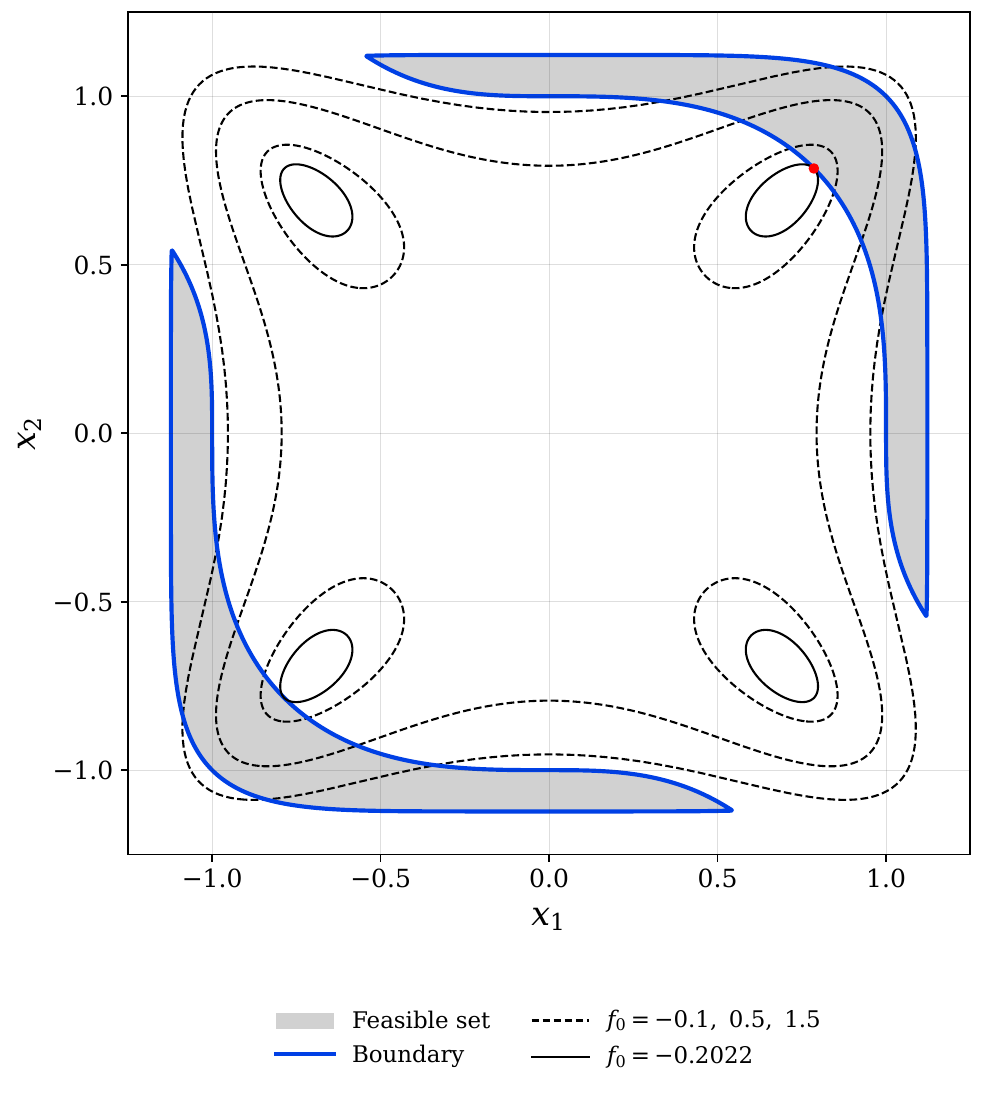}
        \caption{The feasible set of Problem \eqref{eq::enpc}
        and the contour curves of $f_0$.}
        \label{fig:enpc_fs}
    \end{minipage}
    \hfill
    \begin{minipage}[c]{0.5\textwidth}
        \centering
        \captionof{table}{Numerical results for Problem \eqref{eq::enpc}.}
        \label{tab:enpc_results}
        \small
        \setlength{\tabcolsep}{4pt}
        \renewcommand{\arraystretch}{1.1}
        \resizebox{0.7\linewidth}{!}{%
            \begin{tabular}{@{}cccc@{}}
                \toprule
                $k$ & $\gamma_k^{\dsdp}$ & $\bu^{(k)}$ & {\rm time} \\
                \midrule
                 $1$ & $-0.2311$ & $(0.7587,\,0.7587)$ & $0.13\,\mathrm{s}$ \\
                 $2$ & $-0.2258$ & $(0.7650,\,0.7650)$ & $0.14\,\mathrm{s}$ \\
                 $3$ & $-0.2199$ & $(0.7710,\,0.7710)$ & $0.14\,\mathrm{s}$ \\
                 $4$ & $-0.2168$ & $(0.7739,\,0.7739)$ & $0.15\,\mathrm{s}$ \\
                 $5$ & $-0.2135$ & $(0.7768,\,0.7768)$ & $0.15\,\mathrm{s}$ \\
                 $6$ & $-0.2114$ & $(0.7785,\,0.7785)$ & $0.16\,\mathrm{s}$ \\
                 $7$ & $-0.2101$ & $(0.7797,\,0.7797)$ & $0.18\,\mathrm{s}$ \\
                 $8$ & $-0.2084$ & $(0.7810,\,0.7810)$ & $0.20\,\mathrm{s}$ \\
                 $9$ & $-0.2079$ & $(0.7814,\,0.7814)$ & $0.22\,\mathrm{s}$ \\
                $10$ & $-0.2067$ & $(0.7823,\,0.7823)$ & $0.25\,\mathrm{s}$ \\
                \bottomrule
            \end{tabular}%
        }
    \end{minipage}
\end{figure}

\end{example}

\begin{example}\label{example:7.2}
Let $h_i:\mathbb{R} \rightarrow \mathbb{R}$, $i=1,2$ be defined as $h_1(t)=t^3-t$ and $h_2(t)=t^2+t$.
Consider the semi-infinite polynomial optimization problem
\begin{equation}\label{eq::sosc}
\left\{
\begin{aligned}
\min_{\bx\in\RR^2}&\ f_0(\bx) = (h_1(x_1)-1)^2+(h_2(x_2)-1)^2\\
\quad\text{s.t. }&\ f_1(\bx) = (h_1(x_1)-2h_2(x_2)+1)^2-1 \le 0,\\
&\ f(\bx, \by) = \left(h_1(x_1)-\frac{y_1^2+y_1y_2}{2}\right)^2+\left(h_2(x_2)-\frac{y_2^2+y_1y_2}{2}\right)^2-1\le 0,\ \forall \by\in \Y,\\
\end{aligned}
\right.
\end{equation}
where $\Y=[-1,1]^2.$ The feasible set is shown in Figure \ref{fig:sosc_fs} (the gray area). 
%One can check that 
%\[
%f^{\star}=3-2\sqrt{2}\approx
%\]

Assumption \ref{assump::0} holds for this problem with $\mJ=\mJ^{\csosc}_{h,6}$. 
The SOS-convex polynomials corresponding to $f_0$, $f_1$ and $f(\cdot, \by)$, $\by\in\Y$ are
\[
\begin{aligned}
&g_0(\bxi)=(\xi_1-1)^2+(\xi_2-1)^2,\ \  g_1(\bxi)=(\xi_1-2\xi_2+1)^2-1,\\
&g(\bxi,\by)=\left(\xi_1-\frac{y_1^2+y_1y_2}{2}\right)^2+\left(\xi_2-\frac{y_2^2+y_1y_2}{2}\right)^2-1, 
\ \by\in\Y.
\end{aligned}
\]
Note that the Slater condition holds for $\mathcal{F}=\{f_1,f\}$, and none of the 
polynomials $f_0,f_1$ and $f$ is convex in $\bx$.

Solving the SDP relaxations \eqref{eq::uppersdp} and
\eqref{eq::uppersdp_dual} to obtain upper bounds on the true optimal value $f^{\star}$, we obtain 
$\rho_k^{\psdp}=\rho_k^{\dsdp}=0.1716$, $k=2,3,4,5$.  
Although the repeated values suggest that the upper bounds have stabilized numerically, the rank condition~\eqref{eq::FEC} fails at each of these orders. Therefore, finite convergence is not certified. Following the definition of $\Phi_{\mJ^{\csosc}_{h,m}}(L)$ in
Proposition~\ref{prop::subclass2}, we solve the following system of cubic-quadratic  
polynomial equations (which can be efficiently solved and indeed has a closed-form solution):
\begin{equation}\label{eq::equations}
\left\{
\begin{aligned}
x_1^3-x_1
&=L^{(k,\star)}(x_1^3-x_1),\\
x_2^2+x_2
&=L^{(k,\star)}(x_2^2+x_2),
\end{aligned}
\right.
\end{equation}
where $L^{(k,\star)}$ is an optimal solution of
\eqref{eq::uppersdp_dual}, with $k=2$. This system has two solutions
\[
\bu^{(1)}=(1.2510,\,0.4783)
\qquad\text{and}\qquad
\bu^{(2)}=(1.2510,\,-1.4783),
\]
which correspond to two numerical approximate solutions of \eqref{eq::sosc}.
As an illustration, Figure~\ref{fig:sosc_fs}
shows several level curves of $f_0$, including the curve
$f_0(\bx)=0.1716$, and marks $\bu^{(1)}$ and $\bu^{(2)}$ in red.

% Solving the SDP relaxations \eqref{eq::uppersdp} and \eqref{eq::uppersdp_dual} for upper bounds yields $\rho_k^{\psdp} = \rho_k^{\dsdp} = 0.1715$ for $k=2,3,4,5$, which numerically indicates that the finite convergence has occurred. However, the rank condition \eqref{eq::FEC} does not hold for any $k=2,3,4,5$. According to the definition of $\Phi_{\mJ^{\csosc}_{h,m}}(L)$ in Proposition \ref{prop::subclass2}, we solve the system of cubic polynomial equations (which can be efficiently solved, and indeed has a closed-form solution)
% \begin{equation}\label{eq::equations}
%    \left\{
% \begin{aligned}
% x_1^3-x_1&=L^{(k,\star)}(x_1^3-x_1)\\
% x_2^2+x_2&=L^{(k,\star)}(x_2^2+x_2),
% \end{aligned}
% \right. 
% \end{equation}
% with $L^{(k,\star)}$ being an optimal solution of \eqref{eq::uppersdp_dual} and $k=2$. We obtain two points $\bu^{(1)}=(1.2510,  0.4783)$ and $\bu^{(2)}=(1.2510, -1.4783)$ corresponding to the optimal solutions of \eqref{eq::sosc}. 
% To see the accuracy of the results, 
% we plot the contour curves of $f_0$ in Figure \ref{fig:sosc_fs}, among which is the one $f_0(\bx)=0.1715$, and mark the points $\bu^{(1)}$ and $\bu^{(2)}$ in red. 

We next solve the SDP relaxations \eqref{eq::lowersdp} and
\eqref{eq::lowersdp_dual} to obtain lower bounds on the true optimal value $f^\star$, where we take $\nu$ to be Lebesgue measure on $\Y$. We use them to assess the quality of the approximate solution $\bu^{(1)}$ and $\bu^{(2)}$ . For each
order $k=1,\ldots,10$, we compute the lower bound
$\gamma_k^{\dsdp}$ and solve \eqref{eq::equations}, obtaining two points
denoted by $\bu^{(1,k)}$ and $\bu^{(2,k)}$. Table~\ref{tab:sosc_results}
reports these points and the CPU time required to solve
\eqref{eq::lowersdp_dual}. In particular,
$\gamma_{10}^{\dsdp}=0.1711$. Thus, although the rank condition~\eqref{eq::FEC} does not certify the
optimality of the upper bound $\rho_2^{\psdp}=\rho_2^{\dsdp}=0.1716$,
the lower bound $\gamma_{10}^{\dsdp}=0.1711$ shows that  
%$f^\star\in[0.1711,\,0.1715]$. 
% Hence, the optimal value is determined to within $4\times10^{-4}$.
% Next, we solve the SDP relaxations \eqref{eq::lowersdp} and \eqref{eq::lowersdp_dual} to obtain lower bounds on $f^{\star}$.
% For each order $k=1,\ldots,10$,
% we compute the lower bounds $\{\gamma_k^{\dsdp}\}$ and solve the system of equations \eqref{eq::equations}, which yields two points denoted by $\bu^{(1,k)}$ and $\bu^{(2,k)}$. We show these results, as well as the CPU time 
% for solving \eqref{eq::lowersdp_dual}, in Table \ref{tab:sosc_results}. 
% In particular, we have $\gamma_{10}^{\dsdp} = 0.1711$. Although the rank condition \eqref{eq::FEC} does not certify that the upper bound $\rho_2^{\psdp} = \rho_2^{\dsdp} = 0.1715$ is optimal, the computed lower bounds guarantee that, 
%Hence, 
the  approximate solutions found by the upper bound hierarchy,   $\bu^{(1)}$ and $\bu^{(2)}$  are $\varepsilon$-optimal with $\varepsilon<10^{-3}$ (taking into account of numerical rounding error), in the sense that  their objective values differ from the true optimal value by
at most $10^{-3}$.
 %this upper bound is $\varepsilon$-optimal with $\varepsilon=0.0004$. 

\begin{figure}[htbp]
    \centering
    \begin{minipage}[c]{0.45\textwidth}
        \centering
        \includegraphics[width=\linewidth]{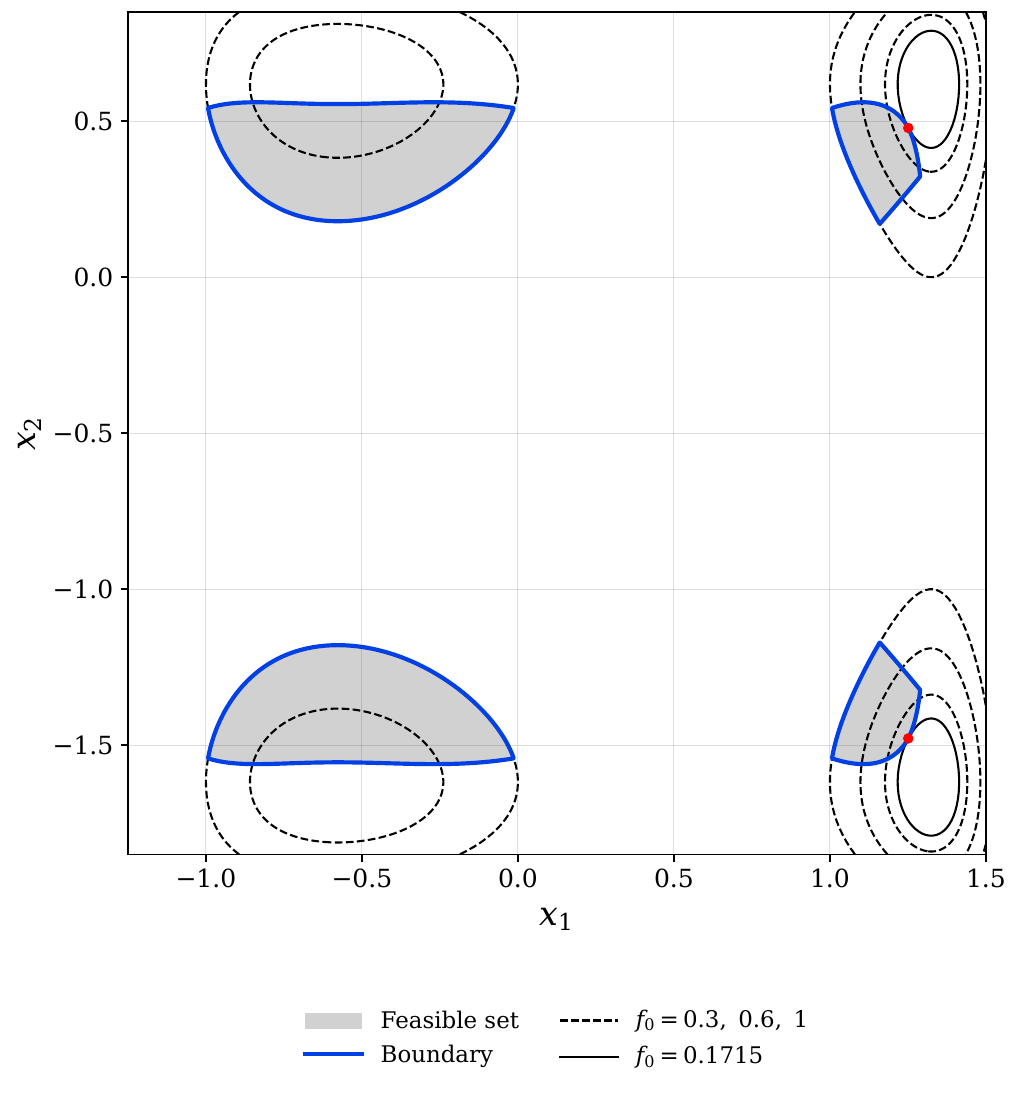}
        \caption{The feasible set of Problem \eqref{eq::sosc}
        and the contour curves of $f_0$.}
        \label{fig:sosc_fs}
    \end{minipage}
    \hfill
    \begin{minipage}[c]{0.48\textwidth}
        \centering
        \captionof{table}{Numerical results for Problem \eqref{eq::sosc}.}
        \label{tab:sosc_results} 
        \small
        \setlength{\tabcolsep}{3pt}
        \renewcommand{\arraystretch}{1.1}
        \resizebox{\linewidth}{!}{%
            \begin{tabular}{@{}ccccc@{}}
                \toprule
                $k$ & $\gamma_k^{\dsdp}$
                    & $\bu^{(1,k)}$ & $\bu^{(2,k)}$ & {\rm time} \\
                \midrule
       $ 1$ & $0.1186$ & $(1.2642,\,0.5032)$ & $(1.2642,\,-1.5032)$ & $0.16${\rm s}\\
       $ 2$ & $0.1480$ & $(1.2566,\,0.4889)$ & $(1.2566,\,-1.4889)$ & $0.16${\rm s}\\
       $ 3$ & $0.1603$ & $(1.2537,\,0.4833)$ & $(1.2537,\,-1.4833)$ & $0.17${\rm s}\\
       $ 4$ & $0.1657$ & $(1.2524,\,0.4809)$ & $(1.2524,\,-1.4809)$ & $0.18${\rm s}\\
       $ 5$ & $0.1683$ & $(1.2518,\,0.4797)$ & $(1.2518,\,-1.4797)$ & $0.18${\rm s}\\
       $ 6$ & $0.1696$ & $(1.2515,\,0.4791)$ & $(1.2515,\,-1.4791)$ & $0.20${\rm s}\\
       $ 7$ & $0.1703$ & $(1.2513,\,0.4788)$ & $(1.2513,\,-1.4788)$ & $0.20${\rm s}\\
       $ 8$ & $0.1708$ & $(1.2512,\,0.4786)$ & $(1.2512,\,-1.4786)$ & $0.23${\rm s}\\
       $ 9$ & $0.1710$ & $(1.2512,\,0.4785)$ & $(1.2512,\,-1.4785)$ & $0.25${\rm s}\\
       $10$ & $0.1711$ & $(1.2511,\,0.4784)$ & $(1.2511,\,-1.4784)$ & $0.30${\rm s}\\
                \bottomrule
            \end{tabular}%
        }
    \end{minipage}
\end{figure}

\end{example}

\section{Conclusions}\label{sec::conclusion}
In this paper, we introduced Jensen classes and their associated Jensen maps
as a framework for extending Jensen-type inequalities beyond
SOS-convexity. This framework includes polynomials with essentially
nonpositive coefficients and, through a composition principle, yields
further classes of possibly nonconvex Jensen-type polynomials. These
include componentwise univariate compositions and compositions with
structured quadratic maps. We then used Jensen-type polynomials to establish theorems of the alternative
for systems of strict polynomial inequalities involving universally
quantified variables. The resulting certificates provide sums-of-squares
conditions for infeasibility and, in turn, nonnegativity certificates
under the assumptions of our results. We also demonstrated the applicability of the framework in two
optimization settings. For source localization, we established exactness
of an SOS relaxation under a dimension condition and gave an explicit
Jensen map for recovering a global minimizer. For semi-infinite
optimization with Jensen-type polynomials, we developed convergent
semidefinite programming hierarchies that provide upper and lower bounds
on the optimal value. We further established finite convergence and
solution-extraction results under the stated assumptions, and illustrated
the methods with numerical examples.

These results suggest several directions for further study, including
identifying additional Jensen classes and composition rules, and
developing more efficient computational methods for the resulting
semidefinite relaxations. It is also worth investigating whether the results can be extended to the matrix setting as in  \cite{Guo1,Guo2}, and used to  study optimality conditions for robust single- or multi-objective optimization \cite{Chuong}.

\smallskip
\noindent \textbf{Statements and Declarations}

\smallskip 
{\footnotesize
 \noindent \textbf{Declaration of Competing Interest}. The authors declare that they have no known competing financial interests or personal relationships that could have appeared to influence the work reported in this paper.

\noindent \textbf{Declaration of AI-Assisted Technologies}. 
% The authors acknowledge that Examples \ref{example:4.1}, \ref{example:4.2} and \ref{example:7.2}, and the proof of Proposition \ref{prop::sigma} in the Appendix were developed with the assistance of DeepSeek.  
% These examples and Proposition \ref{prop::sigma} in the Appendix in the present version were subsequently rewritten and verified by the authors, who assume full
% responsibility for their correctness. 
% "The authors acknowledge DeepSeek’s assistance in developing Examples \ref{example:4.1}, \ref{example:4.2}, and \ref{example:7.2}, and the proof of Proposition \ref{prop::sigma}. The authors subsequently rewrote and verified these materials and take full responsibility for their correctness."
% The authors acknowledge that Deepseek was used for the development in constructing Examples \ref{example:4.1}, \ref{example:4.2} and \ref{example:7.2}, and assist with the proof of Proposition \ref{prop::sigma} in the Appendix. The authors are responsible for the final content and conclusions.  
The authors acknowledge that Examples \ref{example:4.1}, \ref{example:4.2}, and \ref{example:7.2}, as well as the proof of Proposition \ref{prop::sigma} in the Appendix, were developed with the assistance of ChatGPT Edu and DeepSeek. These parts were subsequently rewritten and verified by the authors, who assume full responsibility for their correctness.
}

\appendix

\section{Proof of Lemma \ref{lemma:claim}} \label{sec::proof0}
\begin{proof} We start by first establishing the following claim, which is a restricted version of the desired conclusion with the case when the degree of the even monomial $\bx^{2\boldsymbol{v}}$ is exactly equal to $2s$. 

\noindent\textbf{Claim.} For any $s \in \mathbb{N}$ with $s>0$, $L \in \mathcal{B}[{\bx}]_{2s}$ and any $\boldsymbol{v}\in\mathbb{N}_s^n$ with
$|\boldsymbol{v}|= s$, one has 
\begin{equation}\label{CL}
L(\bx^{2\boldsymbol{v}})
\le
\prod_{i=1}^n
L(x_i^{2s})^{v_i/s}.
\tag{CL}
\end{equation}
Below, we establish  the claim by induction on $n$, where $n$ is the dimension of the vector ${\bx}$. The induction statement is understood to hold simultaneously for every degree parameter 
$s\in\mathbb{N}$ with $s>0$.  

\medskip

To avoid the possible occurrence of zero values in logarithmic
expressions later on, fix $\varepsilon>0$ and define
\[
L_\varepsilon
:=
\frac{L+\varepsilon\delta_{\mathbf{1}}}{1+\varepsilon},
\]
where $\delta_{\mathbf{1}}$ denotes evaluation at
$\mathbf{1}=(1,\ldots,1)$. Then
$L_\varepsilon\in\mathcal{B}[\bx]_{2s}$  and
$L_\varepsilon(\bx^{2\boldsymbol\alpha})>0$ 
for every multi-index $\boldsymbol\alpha$. We first prove the claim \eqref{CL} for
$L_\varepsilon$ and then let $\varepsilon\downarrow0$.

\medskip
\noindent
\textbf{Base case: $n=1$.}  {Direct verification shows that claim \eqref{CL} holds with equality in this case.}

\medskip

\noindent
\textbf{Induction step.} Assume that the claim \eqref{CL} holds in dimension $k-1$ for every degree parameter $s$. We prove it in dimension $k$. Let
\[
\bv
=
(\overline{\bv},v_k),
\qquad
\overline{\bv}
=
(v_1,\ldots,v_{k-1}),
\]
and set
$d:=|\overline{\bv}|=s-v_k$. 

If $d=0$, then $\bv=s \, \be_k$, and the claim \eqref{CL} holds with equality.
Thus, we may assume $d>0$. Define
\[
\widetilde L_\varepsilon(\phi)
:=
L_\varepsilon\big(
x_k^{2v_k}\phi(x_1,\ldots,x_{k-1})
\big), \mbox{ 
for } 
\phi\in\mathbb{R}[x_1,\ldots,x_{k-1}]_{2d}.
\]
If $\phi$ is a sum of squares, then
$\widetilde L_\varepsilon(\phi)\ge0$. 
Moreover,
$\widetilde L_\varepsilon(1)
=
L_\varepsilon(x_k^{2v_k})>0$. 
Hence
$\widehat L_\varepsilon
:=
\frac{\widetilde L_\varepsilon}
{\widetilde L_\varepsilon(1)}$  
is a normalized positive functional on
$\mathbb{R}[x_1,\ldots,x_{k-1}]_{2d}$ and so, $\widehat L_\varepsilon \in \mathcal{B}[x_1,\ldots,x_{k-1}]_{2d}$.

Since the induction hypothesis holds for every degree parameter, it
applies to $\widehat L_\varepsilon$ with degree parameter $d$. Since
$|\overline{\bv}|=d$, we obtain
\[
\widehat L_\varepsilon
\left(
x_1^{2v_1}\cdots x_{k-1}^{2v_{k-1}}
\right)
\le
\prod_{i=1}^{k-1}
\widehat L_\varepsilon(x_i^{2d})^{v_i/d}.
\]
Multiplying by
$\widetilde L_\varepsilon(1)$ gives
\begin{equation}\label{eq:1}
L_\varepsilon(\bx^{2\bv})
\le
\prod_{i=1}^{k-1}
L_\varepsilon(x_i^{2d}x_k^{2v_k})^{v_i/d}.
\end{equation}

We next estimate each factor in \eqref{eq:1}. Fix
$i\in\{1,\ldots,k-1\}$, and define
\[
u_\ell
:=
L_\varepsilon
\left(
x_i^{2\ell}x_k^{2(s-\ell)}
\right),
\qquad
\ell=0,\ldots,s.
\]
All $u_\ell$ are strictly positive. Applying the polynomial Cauchy-Schwarz inequality (CS) to
\[
\phi=x_i^{\ell-1}x_k^{s-\ell+1},
\qquad
q=x_i^{\ell+1}x_k^{s-\ell-1},
\]
gives
$u_\ell^2\le u_{\ell-1}u_{\ell+1}$, 
$\ell=1,\ldots,s-1$.
Thus, for 
$a_\ell:=\log u_\ell$ and $\Delta_\ell:=a_\ell-a_{\ell-1}$, one has 
$\Delta_1\le\Delta_2\le\cdots\le\Delta_s$. 
Since $d+v_k=s$, the average of the first $d$ differences is no
larger than the average of all $s$ differences:
\[
\frac{1}{d}\sum_{\ell=1}^d\Delta_\ell
\le
\frac{1}{s}\sum_{\ell=1}^s\Delta_\ell.
\]
Noting that
$\sum_{\ell=1}^d\Delta_\ell=a_d-a_0$ and 
$\sum_{\ell=1}^s\Delta_\ell=a_s-a_0$, this implies that 
$a_d-a_0
\le
\frac{d}{s}(a_s-a_0)$. 
Since $s-d=v_k$, this gives
$a_d
\le
\frac{v_k}{s}a_0+\frac{d}{s}a_s$. In other words, 
\begin{equation}\label{eq:2}
L_\varepsilon(x_i^{2d}x_k^{2v_k})
\le
L_\varepsilon(x_i^{2s})^{d/s}
L_\varepsilon(x_k^{2s})^{v_k/s}.
\end{equation}
Substituting \eqref{eq:2} into \eqref{eq:1}, and using
$\sum_{i=1}^{k-1}v_i=d$, 
we obtain
\[
\begin{aligned}
L_\varepsilon(\bx^{2\bv})
\le
\prod_{i=1}^{k-1}
\left[
L_\varepsilon(x_i^{2s})^{d/s}
L_\varepsilon(x_k^{2s})^{v_k/s}
\right]^{v_i/d}
&=
\left(\prod_{i=1}^{k-1}
L_\varepsilon(x_i^{2s})^{v_i/s}\right)
L_\varepsilon(x_k^{2s})^{v_k/s}\\
&=
\prod_{i=1}^{k}
L_\varepsilon(x_i^{2s})^{v_i/s}.
\end{aligned}
\]
Therefore, the claim holds for $L_\varepsilon$ in dimension $k$. Thus, by the induction, the claim holds for $L_\varepsilon$ for any dimension $n \in \mathbb{N}$.

Now, since
$L_\varepsilon(\phi)\longrightarrow L(\phi)$
as $\varepsilon\downarrow0$, 
for every polynomial $\phi$ of degree $2s$, passing to the limit
gives
\[
L(\bx^{2\bv})
\le
\prod_{i=1}^{n}
L(x_i^{2s})^{v_i/s}.
\]
This proves the claim \eqref{CL}. 

We next see the desired conclusion. To do this, for any $s \in \mathbb{N}$, $L \in \mathcal{B}[{\bx}]_{2s}$ and any $\boldsymbol{v}\in\mathbb{N}_s^n$ with
$|\boldsymbol{v}| \le s$. From the claim, it suffices to consider the case where $|\boldsymbol{v}| < s$. Define $\widetilde{L}: \mathbb{R}[\bx,z]_{2s} \rightarrow \mathbb{R}$ by 
\[
\widetilde{L}(\psi(\bx,z))= L(\psi(\bx,1)) \mbox{ for all } \psi\in \mathbb{R}[\bx,z]_{2s}
\]
Direct verification shows that $\widetilde{L} \in \mathcal{B}[\bx,z]_{2s}$ 
We now apply the claim \eqref{CL} to the functional $\widetilde L$ and the monomial in the variables $(x_1,\ldots,x_n,z)$ with the multi-index
$\widetilde\bv
=
(v_1,\ldots,v_n,s-|\bv|)$. 
Since $|\widetilde\bv|=s$, we obtain
\[
\widetilde L
\left(
\bx^{2\bv}z^{2(s-|\bv|)}
\right)
\le
\prod_{i=1}^{n}
\widetilde L(x_i^{2s})^{v_i/s}
\,
\widetilde L(z^{2s})^{(s-|\bv|)/s}.
\]
By the definition of $\widetilde L$, one has 
$\widetilde L
\left(
\bx^{2\bv}z^{2(s-|\bv|)}
\right)
=
L(\bx^{2\bv})$, 
\[
\widetilde L(x_i^{2s})
=
L(x_i^{2s}),
\quad i=1,\ldots,n,\quad \mbox{and}\  
\widetilde L(z^{2s})
=
L(1)
=
1.
\]
Consequently,
$L(\bx^{2\bv})
\le
\prod_{i=1}^{n}
L(x_i^{2s})^{v_i/s}$, and hence the conclusion follows. 
\end{proof}

\section{Proof of Proposition \ref{prop::sigma}}\label{sec::proof}

{\begin{proof}
We prove the proposition by using the method of mollification.

If $\mu(\Y)=0$, the assertion follows directly. Hence, we assume that
$\mu(\Y)>0$. Since the family $\{p_t\}_{t\in T}$ is equicontinuous on
the compact set $\Y$, for the given $\varepsilon>0$ there exists
$\delta_\varepsilon>0$ such that
\[
\left|p_t(\bv^{(1)})-p_t(\bv^{(2)})\right|
\leq \frac{\varepsilon}{2\mu(\Y)}
\]
for every $t\in T$ and every
$\bv^{(1)},\bv^{(2)}\in\Y$ satisfying
$\|\bv^{(1)}-\bv^{(2)}\|<\delta_\varepsilon$. Moreover, equicontinuity and pointwise boundedness on the compact set
$\Y$ imply uniform boundedness. Thus,
$N:=\sup_{t\in T}\sup_{\by\in\Y}|p_t(\by)|<+\infty$. 
If $N=0$, then $p_t\equiv 0$ on $\Y$ for every $t\in T$, and the
conclusion is immediate. We now assume that $N>0$.

Define the continuous, nonnegative, compactly supported function
\[
\varphi_\varepsilon(\by):=
\begin{cases}
\displaystyle
\exp\!\left(\frac{1}{\|\by\|-\delta_\varepsilon}\right),
& \|\by\|<\delta_\varepsilon,\\[2mm]
0,& \|\by\|\geq \delta_\varepsilon.
\end{cases}
\]
For $\by\in\Y$, let
$c_\varepsilon(\by)
:=
\int_{\Y}
\varphi_\varepsilon(\bv-\by)\,\ud\nu(\bv)$. 
%Since $\operatorname{supp}(\nu)=\Y$ and
%$\varphi_\varepsilon$ is strictly positive on
%$B(0,\delta_\varepsilon)$, every relatively open neighbourhood of
%$\by\in\Y$ has positive $\nu$-measure. Consequently,
Direct verification shows that $c_\varepsilon(\by)>0$ for all $\by\in\Y$. 
%Moreover, $c_\varepsilon$ is continuous on $\Y$. Since $\Y$ is compact,
%\[
%\min_{\by\in\Y}c_\varepsilon(\by)>0.
%\tag{3}
%\]
Define
\[
g_\varepsilon(\by)
:=
\int_{\Y}
\frac{\varphi_\varepsilon(\by-\bv)}
     {c_\varepsilon(\bv)}
\,\ud\mu(\bv),
\qquad \by\in\mathbb R^\ny.
\]
As $\Y$ is compact, direct verification shows that the function $g_\varepsilon$ is continuous and
nonnegative on $\Y$. Also, its construction depends
only on $\mu$, $\nu$, and $\varepsilon$, and is therefore independent of
$t$. Since $\mu$ and $\nu$ are finite, Fubini's theorem gives, for every
$t\in T$,
\[
\begin{aligned}
\int_{\Y}p_t(\by)g_\varepsilon(\by)\,d\nu(\by)
&=
\int_{\Y}
\left(
\frac{1}{c_\varepsilon(\bv)}
\int_{\Y}
p_t(\by)\varphi_\varepsilon(\by-\bv)\,d\nu(\by)
\right)
d\mu(\bv)\\
&=
\int_{\Y}
\left(
\frac{1}{c_\varepsilon(\by)}
\int_{\Y}
p_t(\bv)\varphi_\varepsilon(\bv-\by)\,d\nu(\bv)
\right)
d\mu(\by).
\end{aligned}
\]
Consequently, by the definition of $c_\varepsilon(\by)$, for all $t \in T$, 
\[
\begin{aligned}
L_{1,t}
&:=
\left|
\int_{\Y}p_t(\by)\,d\mu(\by)
-
\int_{\Y}p_t(\by)g_\varepsilon(\by)\,d\nu(\by)
\right|\\
&\leq
\int_{\Y}
\frac{1}{c_\varepsilon(\by)}
\int_{\Y}
|p_t(\by)-p_t(\bv)|
\varphi_\varepsilon(\bv-\by)\,d\nu(\bv)
\,d\mu(\by)\\
&\leq
\frac{\varepsilon}{2\mu(\Y)}
\int_{\Y}
\frac{c_\varepsilon(\by)}{c_\varepsilon(\by)}
\,d\mu(\by)
=\frac{\varepsilon}{2}.
\end{aligned}
\]
%This estimate is uniform in $t\in T$.

Since $g_\varepsilon$ is continuous and nonnegative on the compact set
$\Y$, the quantity
$M:=\max_{\by\in\Y}\sqrt{g_\varepsilon(\by)}<\infty$. Choose $\varepsilon'>0$ such that
$\varepsilon'(2M+\varepsilon')
<
\frac{\varepsilon}{4N\nu(\Y)}$. 
By the Stone--Weierstrass theorem, there exists a polynomial
$h_\varepsilon\in\mathbb R[\by]$ satisfying
\[
\left|
h_\varepsilon(\by)-\sqrt{g_\varepsilon(\by)}
\right|<\varepsilon', \ 
\qquad \by\in\Y.
\]
Set
$\widetilde{\sigma}_\varepsilon(\by)
:=
h_\varepsilon(\by)^2\in\Sigma[\by]$.
Then, for every $\by\in\Y$,
\begin{equation}\label{eq::ineq}
\begin{aligned}
|g_\varepsilon(\by)-\widetilde{\sigma}_\varepsilon(\by)|
&=
\left|
h_\varepsilon(\by)-\sqrt{g_\varepsilon(\by)}
\right|
\left|
h_\varepsilon(\by)+\sqrt{g_\varepsilon(\by)}
\right|
\leq
\varepsilon'(2M+\varepsilon')
<
\frac{\varepsilon}{4N\nu(\Y)}.
\end{aligned}
\end{equation}
Consequently, for all $t\in T$,
\[
\begin{aligned}
L_{2,t}
&:=
\left|
\int_{\Y}p_t(\by)g_\varepsilon(\by)\,d\nu(\by)
-
\int_{\Y}p_t(\by)\widetilde{\sigma}_\varepsilon(\by)\,d\nu(\by)
\right| \leq
N\int_{\Y}
|g_\varepsilon(\by)-\widetilde{\sigma}_\varepsilon(\by)|
\,d\nu(\by) \leq \frac{\varepsilon}{4}.
\end{aligned}
\]
Now, using Fubini's Theorem again, it holds $\int_{\Y}g_{\varepsilon}(\by)\ud\nu(\by)=\mu(\Y)>0$. 
    Then due to \eqref{eq::ineq}, we can
    assume without loss of generality that $\int_{\Y}\td{\sigma}_{\varepsilon}(\by)\ud\nu(\by)>0$.
    Now let 
    \[
    \sigma_{\varepsilon}(\by)=\frac{\int_{\Y}g_{\varepsilon}(\bv)\ud\nu(\bv)}{\int_{\Y}\td{\sigma}_{\varepsilon}(\bv)\ud\nu(\bv)}\td{\sigma}_{\varepsilon}(\by)
    =\frac{\mu(\Y)}{\int_{\Y}\td{\sigma}_{\varepsilon}(\bv)\ud\nu(\bv)}\td{\sigma}_{\varepsilon}(\by).
    \]
    We have $\int_{\Y}\sigma_{\varepsilon}(\by)\ud\nu(\by)=\mu(\Y)$ and by \eqref{eq::ineq}, for all $t \in T$,
    \[
    \begin{aligned}
    L_{3,t}:=&\left\vert\int_{\Y}p_t(\by)\td{\sigma}_{\varepsilon}(\by)\ud\nu(\by)
    -\int_{\Y}p_t(\by)\sigma_{\varepsilon}(\by)\ud\nu(\by)\right\vert\\
    \le& 
    \int_{\Y}\vert p_t(\by)\vert\frac{\left\vert\int_{\Y}\td{\sigma}_{\varepsilon}(\bv)\ud\nu(\bv)-\int_{\Y}g_{\varepsilon}(\bv)\ud\nu(\bv)\right\vert}{\int_{\Y}\td{\sigma}_{\varepsilon}(\bv)\ud\nu(\bv)}\td{\sigma}_{\varepsilon}(\by)\ud\nu(\by)\\
    \le& \frac{\varepsilon}{4}\int_{\Y}\frac{1}{\int_{\Y}\td{\sigma}_{\varepsilon}(\bv)\ud\nu(\bv)}\td{\sigma}_{\varepsilon}(\by)\ud\nu(\by)=\frac{\varepsilon}{4}.
    \end{aligned}
    \]
Therefore, for all $t \in T$,
    \[
    \left\vert \int_{\Y}p_t(\by)\ud\mu(\by)-\int_{\Y}p_t(\by)\sigma_{\varepsilon}(\by)\ud\nu(\by)\right\vert
    \le L_{1,t}+L_{2,t}+L_{3,t}\le\varepsilon.
    \]
    Finally, note that 
the polynomial $\sigma_\varepsilon$ is constructed independently of
$t$, and so, the proof is complete.
\end{proof}}

\bibliographystyle{abbrv}
\bibliography{main,new}

\end{document}